\documentclass[3p,times]{elsarticle}
\usepackage{graphicx}
\usepackage{sidecap}
\usepackage{amssymb,amsmath,amsthm,amsbsy}
\usepackage{calligra,calrsfs,mathrsfs}
\usepackage{multirow}
\usepackage{ulem,cancel}
\usepackage[usenames,dvipsnames]{color}
\usepackage{copyrightbox}
\usepackage{cases}
\usepackage{pgfplots}
\usepackage{subcaption}
\usepackage{tikz} 
\usepackage{tkz-euclide}
\usepackage{tikz-dimline}
\usetikzlibrary[patterns]
\usepackage{hyperref}
\usepackage[ruled,vlined]{algorithm2e}
\usepackage{algpseudocode}
\usepackage{dsfont}
\usepackage{caption}
\usepackage{epstopdf}
\usepackage{epsfig}
\usepackage{natbib}
\usepackage{grffile}
\usepackage{comment}
\usepackage{enumerate}
\usepackage{booktabs}
\usepackage{accents}
\usepackage{colortbl}
\usepackage{todonotes}
\usepackage{array}
\newcolumntype{P}[1]{>{\centering\arraybackslash}p{#1}}
\usepackage{float}
\usepackage[T1]{fontenc}
\usepackage[latin9]{inputenc}
\usepackage{geometry}
\usepgfplotslibrary{fillbetween}
\usetikzlibrary{patterns}
\usepackage{color}
\usepackage{stmaryrd}
\usepackage{setspace}
\usepackage{amsthm}
\usepackage{textcomp}
\usepackage{hyperref}
\usepackage{soul,xcolor}
\usepackage[inline]{enumitem}
\usepackage{textcomp}
\usepackage{adjustbox}
\usepackage{rotating}

\allowdisplaybreaks

\newcommand{\ifcomment}{\iffalse}

\newcommand{\ssf}[1]{\mbox{\textsf{#1}}}
\newcommand{\bs}[1]{\boldsymbol{#1}}

\newcommand{\bb}[1]{{\mbox{\sffamily \bfseries #1}}}

\newdefinition{rem}{Remark}

\makeatletter
\def\smallunderbrace#1{\mathop{\vtop{\m@th\ialign{##\crcr
   $\hfil\displaystyle{#1}\hfil$\crcr
   \noalign{\kern3\p@\nointerlineskip}%
   \tiny\upbracefill\crcr\noalign{\kern3\p@}}}}\limits}
\makeatother

\newcommand{\ts}{\textsubscript}

\biboptions{sort&compress}

\journal{Computer Methods in Applied Mechanics and Engineering}

\begin{document}

\begin{frontmatter}

\title{High-order Multiscale Preconditioner for Elasticity of Arbitrary Structures}
\author[psu]{Sabit Mahmood Khan}
\ead{skk6071@psu.edu}
\author[psu]{Yashar Mehmani}
\ead{yzm5192@psu.edu}
\address[psu]{Energy and Mineral Engineering Department, The Pennsylvania State University, University Park, Pennsylvania 16802}
\cortext[ca]{Corresponding author: Yashar Mehmani. Email: yzm5192@psu.edu}

\begin{abstract}
We present a two-level preconditioner for solving linear systems arising from the discretization of the elliptic, linear-elastic deformation equation, in displacement unknowns, over domains that have arbitrary geometric and topological complexity and heterogeneity in material properties (including fractures). The preconditioner is an algebraic translation of the high-order pore-level multiscale method (hPLMM) proposed recently by the authors, wherein a domain is decomposed into non-overlapping subdomains, and local basis functions are numerically computed over the subdomains to construct a high-quality coarse space (or prolongation matrix). The term ``high-order'' stands in contrast to the recent low-order PLMM preconditioner, where BCs of local basis problems assume rigidity of all interfaces shared between subdomains. In hPLMM, interfaces are allowed to deform, through the use of suitable mortar spaces, thereby capturing local bending/twisting moments under challenging loading conditions. Benchmarked across a wide range of complex (porous) structures and material heterogeneities, we find hPLMM exhibits superior performance in Krylov solvers than PLMM, as well as state-of-the-art Schwarz and multigrid preconditioners. Applications include risk analysis of subsurface CO$_2$/H$_2$ storage and optimizing porous materials for batteries, prosthetics, and aircraft.
\end{abstract}

\end{frontmatter}

\section{Introduction}\label{sec:intro}

The accurate prediction of mechanical deformation in solids is important to many engineering applications. In underground CO$_2$ or H$_2$ storage, injected fluids can stress the overburden and surrounding reservoir rock, potentially activating faults or leakage pathways through the caprock \cite{krevor2023subsurface}. By contrast, in geothermal energy extraction and unconventional hydrocarbon recovery, injection-induced fracturing is explicitly pursued, but the challenge lies in the high-precision control of such fractures \cite{mcclure2014geothermal}. In material science and manufacturing, porous microstructures that are macroscopically lightweight yet high-strength are desired for building fuel-efficient aircraft  \cite{marx2020airplane}, shock-absorbing armors \cite{jung2011armor}, durable battery electrodes \cite{mukhopadhyay2014electrodes}, and prosthetic implants and scaffolds for osteoporosis \cite{ziaie2024osteo}. A prerequisite to tackling such problems is an ability to solve deformation equations on structures with arbitrary geometric complexity.

Here, we focus on the linear-elastic regime and consider systems of the form $\mathrm{A}x\!=\!b$ that arise from discretizing the governing equations in the displacement unknown, $x$, on a given domain $\Omega$. While the specific discretization scheme used does not limit our methods, classical finite elements (FEM) is a natural choice. We call the mesh defined on $\Omega$, and used to assemble the system, the \textit{fine grid}. When $\Omega$ is geometrically complex (i.e., highly negative Euler characteristic \cite{edels2010CompTopology}), it must often be finely gridded to capture its intricate details. As a result, $\mathrm{A}$ becomes commensurately large and ill-conditioned \cite{saad2003book}. An example is the microstructure of a porous material (e.g., metal foam) characterized by a high-resolution X-ray $\mu$CT image \cite{wildenschild2013x}. Iterative solvers, like GMRES, are the only viable option to solve such linear systems, but require effective preconditioning to achieve rapid convergence. While numerous preconditioners exist for the problem at hand, the most successful are based on domain decomposition (e.g., Schwarz, FETI, GDSW) \cite{toselli2004book, dolean2015schwarz}, multigrid ideas (e.g., AMG, cAMG) \cite{stuben1987algebraic, gustafsson1998cAMGorig, notay2010agmg}, or multiscale methods \cite{mehmani2021striv} translated into algebraic preconditioners (e.g., MsFE, GMsFEM, MoMsFE) \cite{buck2013MsFElin,
buck2014MsFEoscil, castelletto2017mult, yang2022GMsFEM, arbogast2015two}. Almost all are two-level in design and consist of a fine smoother, M\ts{L}, and coarse preconditioner, M\ts{G}, to attenuate high- and low-frequency errors, respectively. The latter consists of a prolongation, P, and a restriction, R, matrix, with $\mathrm{R}\!=\!\mathrm{P}^\top$ often being the case. A key to the success of any two-level preconditioner is for the column space of P to contain a good approximation to the solution $x$. Because then, this approximation would equal $x_{aprx}\!=\!\mathrm{P} x^o$, where $x^o$ is the solution to a much smaller coarse system $\mathrm{A}^o x^o\!=\!b^o$, with $\mathrm{A}^o\!=\!\mathrm{RAP}$ and $b^o\!=\!\mathrm{R}b$. The main difference among the cited two-level preconditioners is in the way P is constructed. If P is ``good,'' errors in $x_{aprx}$ tend to be dominated by high-frequency modes, which are wiped out (iteratively) by M\ts{L}.

Two-level preconditioners based on multiscale methods have a distinct advantage over other, more black-box variants, in that the underlying material heterogeneity (e.g., in stiffness), geometry (e.g., fractures), and physics (e.g., PDE form) of the problem are embedded into the construction of P. 
This is accomplished by restricting the PDE onto subdomains (or coarse grids), $\Omega_i$, and solving it to construct a set of local basis functions. The boundary conditions (BCs) assumed in solving such local problems, which we refer to as \textit{closure BCs}, are the most important determinant of the quality of P. 
In multiscale finite element/volume (MsFE/V) \cite{hou1997multiscale, castelletto2017mult, jenny2003multi, sokolova2019MsFVbiot, xu2025multixfem}, closure BCs are either linearly varying displacements or solutions of the restricted PDE along portions of the subdomain boundary, $\partial\Omega_i$. In multiscale mortar methods (MoMsFE) \cite{bernardi1989new, belgacem1999mortar, arbogast2007multiscale, khattatov2019domain, moretto2024novel}, further flexibility is offered through the use of mortars, a low-dimensional function space defined on $\partial\Omega_i$. In two-level Schwarz, with Generalized Dryja-Smith-Widlund (GDSW) coarse spaces for example \cite{heinlein2019adaptive, dohrmann2008family, heinlein2020fully}, closure BCs are either restrictions of rigid body motions of $\Omega_i$ onto $\partial\Omega_i$ or eigenfunctions computed on $\partial\Omega_i$ (similar to \cite{dolean2012eigenBCs} for Dirichlet-to-Neumann spaces).
Nearly all of the above methods focus on continuum domains with relatively simple exterior geometry, $\partial\Omega$, placing special emphasis on convergence difficulties posed by high-contrast stiffness fields. The decomposition of $\Omega$ is not of primary concern, and assumed given by a software like METIS \cite{karypis1997METISsoftware}.

When $\partial\Omega$ exhibits extreme geometric and topological complexity, such as a $\mu$CT image of a porous microstructure, the decomposition becomes important and constitutes a key step in imposing high-quality closure BCs. The pore-level multiscale method (PLMM), developed recently by the authors \cite{mehmani2021contact, li2023pore} and translated into a two-level preconditioner in \cite{mehmani2023precond, li2024phase}, utilizes this information in building P. The solid $\Omega$ is decomposed into non-overlapping subdomains by cutting across its geometric constrictions using the watershed transform \cite{beucher1979water}. For a grain pack, subdomains would correspond to the grains and interfaces to the contacts between the grains. Local basis functions are then built using closure BCs that assume the displacement field over each interface between neighboring subdomains is uniform (i.e., interfaces are rigid). In \cite{khan2024crack}, it was found such BCs accurately capture deformation under global tension/compression, but incur higher errors under loading conditions that induce significant local bending/torsion moments on $\partial\Omega_i$, e.g., global shear. To address this limitation, a high-order generalization of PLMM, called hPLMM, was developed by \cite{khan2024high}, where special mortar spaces were utilized for arbitrarily complex interfaces that result from decomposing $\Omega$. 

The goal of this paper is to algebraically translate hPLMM into a two-level preconditioner and compare it to the PLMM preconditioner of \cite{li2024phase}, as well as Schwarz (GDSW) and multigrid (cAMG) preconditioners as benchmarks. To distinguish between the geometric method of \cite{li2024phase} and the preconditioner herein, we use the terms \textit{algebraic hPLMM} and \textit{geometric hPLMM} hereafter. If we drop the prefix, ``hPLMM'' refers to the preconditioner. 
We test hPLMM within GMRES on a wide range of geometries (pore-scale and continuum) and stiffness heterogeneities (including fractures) and observe superior convergence in both iteration count and wall-clock time. When the dimension of the mortar space on each interface, $n$, is one, hPLMM reduces to PLMM. While the sum of the offline cost of building hPLMM ($n\! \in\![2,4]$) and the online cost of applying it in GMRES is comparable to PLMM ($n\!=\!1$), hPLMM requires more time offline than online. This is a clear advantage if hPLMM is used to solve multiple systems with differing source terms, BCs, or even slight perturbations in the coefficient matrix A \cite{li2024phase}. Examples include solves over many load/time steps and/or Newton/staggered iterations in computational plasticity, finite-strain mechanics, wave propagation, and phase-field modeling of fracture \cite{neto2011PlasticFiniteStrain, hughes2012finite, ambati2015review}. In such problems, the offline cost is negligible and the speedup of hPLMM over PLMM exceeds a factor of two. In addition to the linear-elastic, elliptic PDE studied here, the algebraic formulation of PLMM for saddle-point (e.g., Stokes) PDEs \cite{mehmani2018fluid, guo2019mult} has recently been proposed \cite{mehmani2025stokes}.

The paper is organized as follows: Section \ref{sec:prob_des} describes the PDE we aim to solve, and for which we propose a preconditioner. Section \ref{sec:hplmm} briefly reviews the geometric formulation of hPLMM by \cite{khan2024high}. In Section \ref{sec:multi_precon}, we detail the algebraic formulation of hPLMM as a two-level preconditioner. Section \ref{sec:prob_set} outlines the problem set we consider to test and validate hPLMM, both as an approximate coarse-scale solver and as a preconditioner. Section \ref{sec:results} presents our results, and Section \ref{sec:discussion} discusses their broader implications and future extensions. Section \ref{sec:conclusion} concludes the paper.

\section{Problem description}\label{sec:prob_des}

\begin{figure} [t!]
  \centering
  \hspace*{-0.6cm}
  \centerline{\includegraphics[scale=0.30,trim={0 0 0 0},clip]{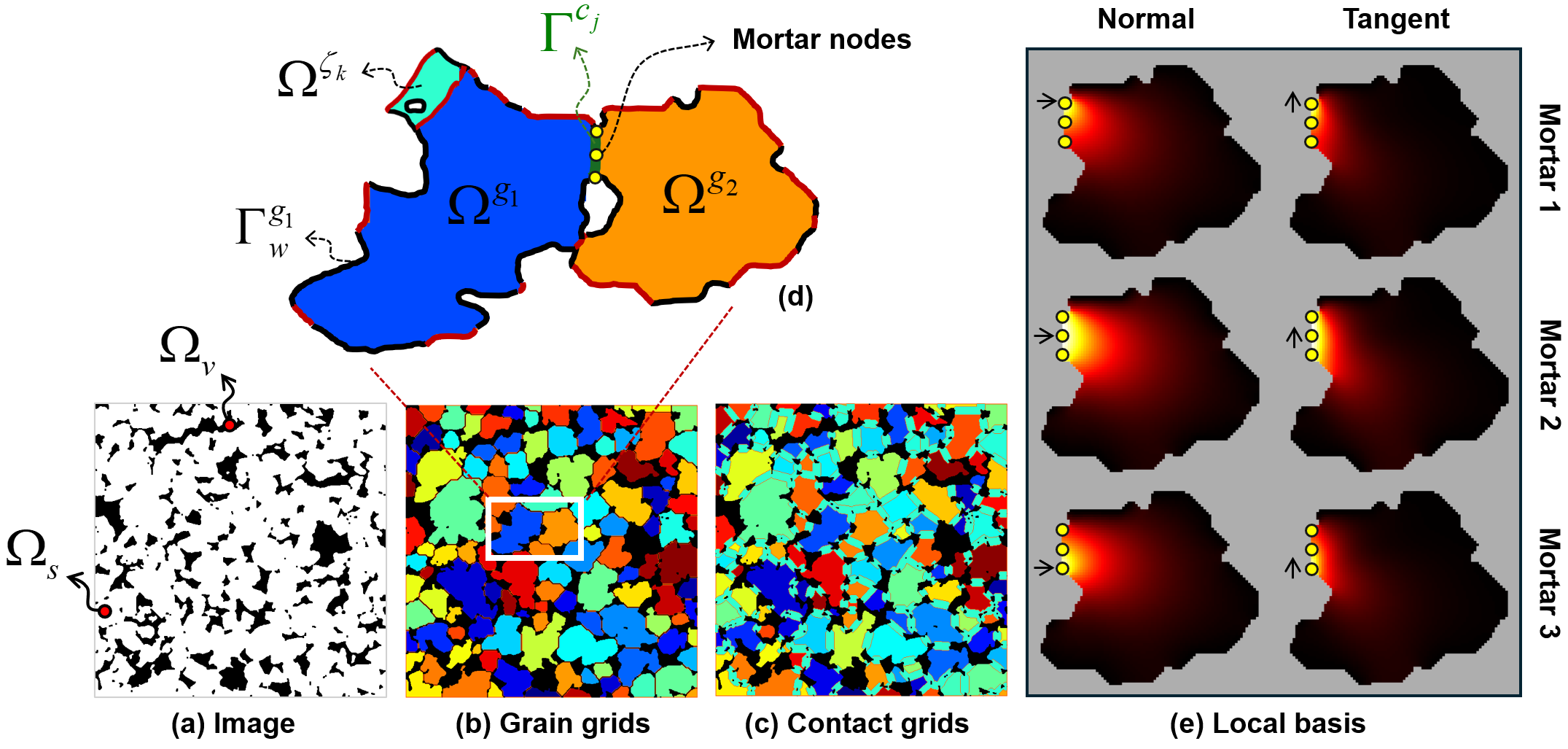}}
	\caption{Schematic of a porous sample captured by a binary image, and its decomposition into grain grids and contact grids. (a) The domain $\Omega$ consists of a solid phase $\Omega_s$ (white) and a void space $\Omega_v$ (black). (b) $\Omega_s$ is decomposed into non-overlapping grain grids $\Omega^{g_i}$ (randomly colored). (c) Contact grids $\Omega^{\zeta_k}$ (cyan) cover a thin region around each contact interface $\Gamma^{c_j}$ (red) and are used to reduce approximation errors. (d) A zoom-in of two adjacent grain grids $\Omega^{g_1}$ and $\Omega^{g_2}$ that share an interface $\Gamma^{c_j}$. The boundary $\partial\Omega^{g_1}$ consists of a fluid-solid interface $\smash{\Gamma^{g_1}_w}$ (black) and contact interfaces shared with other grain grids, like $\smash{\Gamma^{c_j}}$ (green). (e) Displacement-magnitude fields of six shape functions defined on $\Omega^{g_1}$, corresponding to normal (compressive) and tangential (shear) motions imposed at three mortar nodes (yellow dots). These BCs are annotated by black arrows. }
\label{fig:schem}
\end{figure}

Suppose $\Omega\!\subset\!\mathbb{R}^D$ represents a domain with the solid phase $\Omega_s$ and void space $\Omega_v$, where $D$ is the spatial dimension (Fig.\ref{fig:schem}). Let $\partial\Omega_s$ be the Lipschitz boundary of $\Omega_s$ and $\Omega^\circ_s$ its interior. The void-solid interface $\Gamma^w\!=\!\partial\Omega_s\!\cap\!\partial\Omega_v$ and the external boundary $\Gamma^{ex}$ of $\Omega$ comprise $\partial\Omega_s$. The equation we aim to solve on $\Omega^\circ_s$ is the linear-elastic deformation:
\begin{subequations}\label{eq:Linear-elastic}
\begin{equation} \label{eq:Mom}
	\nabla\cdot\bs{\sigma}(\bs{u})=\bs{f}
\end{equation}
\begin{equation} \label{eq:Stress}
	\bs{\sigma}(\bs{u})=\bb{C} : \nabla^s \bs{u}
\end{equation}
\end{subequations}
subject to the following boundary conditions (BCs) on $\partial\Omega_s$:
\begin{subequations}\label{eq:BC}
\begin{align}
	&\left.\bs{u}\,\right|_{\,\Gamma^{D}}= \bs{u}_D
	\; &Dirichlet
	\label{eq:BCb} \\
	&\left.\bs{\sigma}\cdot\bs{n}\,\right|_{\,\Gamma^N}=0 \qquad	
	\left.\bs{\sigma}\cdot\bs{m}\,\right|_{\,\Gamma^N}= 0
	\; &Neumann
\end{align}
\end{subequations}
where $\bs{\sigma}$, $\bs{u}$, $\bb{C}$, and $f$ denote the Cauchy stress tensor, displacement vector, stiffness tensor, and the body force, respectively. We partition the boundary $\partial\Omega_s$ into Dirichlet, $\Gamma^D$, and Neumann, $\Gamma^N$ segments, such that the following hold: $\partial\Omega_s\!=\!\Gamma^D\cup\Gamma^N$, $\Gamma^{ex}\!\subset\!\Gamma^D\cup\Gamma^N$, and $\Gamma^w\!\subset\!\Gamma^N$. The vectors $\bs{n}$ and $\bs{m}$ in Eq.\ref{eq:BC} are the unit normal and unit tangent on $\Gamma^N$, respectively. In this work, we only consider small-strain deformations, where the strain tensor, $\bs{\varepsilon}$, equals the symmetric gradient of the displacement field, $\nabla^{s}\bs{u}\!=\!\left(\nabla\bs{u}+\nabla\bs{u}^{\top}\right)/2$. 

We assume the solid $\Omega_s$ is isotropic with the following stiffness tensor:
\begin{equation}\label{eq:stiff_iso}
	\ssf{C}_{ijkl} = \lambda\delta_{ij}\delta_{kl} + \mu(\delta_{ik}\delta_{jl}+\delta_{il}\delta_{jk}) \; 
\end{equation}
which is described by two Lam\'e parameters: $\lambda$ and $\mu$. We remark that BCs other than Eq.\ref{eq:BC} can also be considered (e.g., roller) without loss of generality to the methods proposed later, but are omitted here for brevity (see \cite{khan2024crack}). Moving forward, we shall adopt bold symbols to denote vectors and tensors and non-bold symbols to denote scalars.

Eq.\ref{eq:Linear-elastic} is discretized using a Galerkin FEM over a Cartesian fine grid on $\Omega_s$. The elements are rectangular/cuboid and the FEM shape functions defined on the elements are bilinear/trilinear in 2D/3D. This yields a linear system: 
\begin{equation} \label{eq:ls_all}
	\hat{\mathrm{A}}\hat{x} = \hat{b}
\end{equation}
where $\hat{\mathrm{A}}$ is the coefficient matrix, $\hat{b}$ is the RHS vector, and $\hat{x}$ is the unknown vector of nodal displacements $\bs{u}$. The Galerkin FEM discretization renders $\hat{\mathrm{A}}$ symmetric. In many practical applications, $\Omega_s$ is very large and finely meshed, resulting in a large $\hat{\mathrm{A}}$. Hence, iterative (e.g., Krylov) methods are the only viable option for solving Eq.\ref{eq:ls_all}. However, rapid convergence depends on the availability of effective preconditioners. In this work, we present such a preconditioner based on an algebraic reformulation of the recent high-order pore-level multiscale method (hPLMM) \cite{khan2024high}. 

\section{Review of the high order pore-level multiscale method (hPLMM)}
\label{sec:hplmm}

Since our hPLMM preconditioner is based on the geometric method of \cite{khan2024high}, we briefly review the latter to develop the reader's intuition about some of the key algebraic operations performed to build the preconditioner. Throughout this section, all references to ``hPLMM'' imply the geometric variant of \cite{khan2024high} unless otherwise stated. In \cite{khan2024high}, hPLMM was proposed to generalize its low-order predecessor PLMM \cite{mehmani2021contact}. Both consist of four steps: (1) decompose $\Omega_s$ into non-overlapping subdomains; (2) build shape/correction functions on each subdomain; (3) couple local functions with a coarse global problem that imposes force balance and kinematic constraints at subdomain interfaces; and (4) interpolate the coarse solution onto the fine grid using the shape/correction functions. The difference between PLMM and hPLMM is that the former assumes the displacement field is uniform over each interface, implying the interfaces are rigid. This limits PLMM's ability to account for local moments. hPLMM removes this drawback by capturing non-uniform displacements on each interface with the use of mortars. Another benefit of hPLMM is that the quality of the approximate solution obtained no longer depends on the quality of the decomposition in Step 1, unlike PLMM where this dependence is strong. The following sections elaborate on each of the above steps for hPLMM.

\subsection{Domain decomposition} \label{sec:domdecomp}
The solid $\Omega_s$ is decomposed into non-overlapping subdomains, $\Omega^{g_i}$, called \textit{grain grids}. They are shown as randomly colored regions in Fig.\ref{fig:schem}b for the domain in Fig.\ref{fig:schem}a. The interfaces shared between adjacent grain grids, $\Gamma^{c_j}$, are referred to as \textit{contact interfaces} (green line in Fig.\ref{fig:schem}d). We consider three ways of decomposing $\Omega_s$: (1) watershed segmentation; (2) spectral partitioning; and (3) Cartesian decomposition. They are detailed below.
\vspace*{5pt}

\textbf{Watershed segmentation:} This morphological operation in image analysis \cite{beucher1979water} decomposes $\Omega_s$, represented by a binary image like Fig.\ref{fig:schem}a, into grain grids by cutting along geometric bottlenecks. Hence, $\Gamma^{c_j}$ coincides with such bottlenecks, and $\Omega^{g_i}$ with local geometric enlargements of $\Omega_s$. For the case of a granular medium, $\Omega^{g_i}$ represents physical grains and $\Gamma^{c_j}$ the grain contacts. Watershed consists of two steps: (1) compute a Euclidean distance map of the image; (2) use local minima of this map as ``seeds'' for the grain grids, which are then grown until shared interfaces form. The number of grain grids is controlled by the number of starting seeds. For details, see \cite{meyer1994marker, khan2024crack}. The limitation lies in that watershed cannot be applied to gray-scale, or non-binary, images corresponding to continuous, or non-discrete, spatial variations in a material property like stiffness $\bb{C}$. This is typical of continuum or Darcy-scale domains. Thus, the low-order PLMM, which heavily relies on watershed, applies to only pore-scale domains. Also, if the input image is noisy or low-resolution, the decomposition quality is poor and PLMM's accuracy deteriorates.

\textbf{Spectral partitioning:} Spectral partitioning is a graph-based decomposition that cuts a graph across its weakest links, or connectivity bottlenecks, and it is the algorithm behind popular software like METIS \cite{karypis1997METISsoftware}. In that sense, spectral partitioning generalizes watershed segmentation as it applies also to gray-scale images, because any such image can be converted into a graph that captures the local connectivity of its pixels. The way this algorithm works is: (1) to obtain $k$ grain grids, find the first $k+1$ eigenvectors of the normalized graph Laplacian, corresponding to its $k+1$ smallest eigenvalues. The first eigenvalue is zero and is discarded; (2) use the eigenvector entries as features for each node in the graph and apply a k-means clustering to partition the nodes into grain grids. This is called a \textit{k-way partitioning}, which can be contrasted with the common alternative of \textit{recursively bisecting} the graph. We only use the former in this work. For details, we refer the reader to \cite{khan2024crack} and the references therein.

\textbf{Cartesian decomposition:} This consists of simply slicing $\Omega_s$ along Cartesian directions into rectangular/cuboid grain grids in 2D/3D. Needless to say, the decomposition is uninformed by geometry or material heterogeneity.
\vspace*{5pt}

To close, hPLMM uses a second set of subdomains $\Omega^{\zeta_k}$, called \textit{contact grids} (cyan colored patches in Fig.\ref{fig:schem}c). Each covers a contact interface and a small region around it. They are created by performing  morphological ``dilations,'' an operation in image analysis, of the pixels comprising a contact interface. The width of $\Omega^{\zeta_k}$ is a user-defined parameter, but 16 pixels is often enough. For details, see \cite{mehmani2021contact}. Unlike \cite{khan2024high}, where overlapping contact grids were merged, here (as in \cite{li2024phase}) we avoid such mergers to ensure the size of each contact grid remains small to permit localized computations.

\subsection{Mathematical notation}\label{sec:note}
Moving forward, we adopt a notation consistent with \cite{khan2024high}. We use $g_i$, $c_j$, and $\zeta_k$ as labels for entities associated with $\Omega^{g_i}$, $\Gamma^{c_j}$, and $\smash{\Omega^{\zeta_k}}$, respectively. The boundary of a grain grid, $\smash{\partial\Omega^{g_i}}$ consists of: (1) contact interfaces shared with other grain grids $\smash{\Gamma^{g_i}_{c_j}\!=\!\partial\Omega^{g_i}\cap\Gamma^{c_j}}$ (green in Fig.\ref{fig:schem}d), and (2) a void-solid interface $\smash{\Gamma^{g_i}_w\! =\!\partial\Omega^{g_i}\cap\Gamma^w}$ (black in Fig.\ref{fig:schem}d). Note that two grain grids may share more than one contact interface (Fig.\ref{fig:schem}d). Similarly, a contact-grid boundary $\smash{\partial\Omega^{\zeta_k}}$ is made up of: (1) a void-solid interface $\smash{\Gamma^{\zeta_k}_w\!=\!\partial\Omega^{\zeta_k}\cap\Gamma^w}$, and (2) segments that intersect the interior of nearby grain grids, $\smash{\Gamma^{\zeta_k}_{g_i}\!=\!\partial\Omega^{\zeta_k}\cap\Omega^{g_i}}$. We note $\partial\Omega^{g_i}$ and $\partial\Omega^{\zeta_k}$ may also intersect the external boundary $\Gamma^{ex}$ of $\Omega_s$, in which case they consist of a third segment where the global BCs in Eq.\ref{eq:BC} are imposed. We define $\smash{C^{g_i}\!=\!\lbrace c_j\,|\,\Gamma^{g_i}_{c_j}\neq\emptyset\rbrace}$ as the set of all $\Gamma^{c_j}$ that intersect $\partial\Omega^{g_i}$, and $\smash{G^{c_j}\!=\!\lbrace g_i\,|\,\Gamma^{g_i}_{c_j}\neq\emptyset\rbrace}$ as the set of all grain grids sharing $\Gamma^{c_j}$. Given $\Gamma^{c_j}$ is shared between two grain grids, $G^{c_j}$ has only two members. We use $N^g$, $N^c$, and $N^{\zeta}$ to denote the number of grain grids, contact interfaces, and contact grids, respectively. We use the term \textit{coarse-scale} for entities associated with $\Omega^{g_i}$ or $\Omega^{\zeta_k}$. Entities associated with the fine grid are termed \textit{fine-scale}. We use the superscripts $f$ and $o$ to indicate fine- and coarse-scale variables, respectively. An entity is \textit{local} if it is confined to a single coarse grid, and \textit{global} if its scope spans $\Omega_s$. For simplicity, we express all hPLMM equations using 2D notation. For example, the unit tangent $\bs{m}$ on a boundary consists of one vector in 2D, but two orthonormal vectors in 3D. We shall use the former to express all local BCs involving $\bs{m}$.

\subsection{Mortar nodes and mortar functions}\label{sec:mortar}
\begin{figure}[t!]
	\centering	
	\includegraphics[trim=15 45 10 65, clip, width=0.8\linewidth]{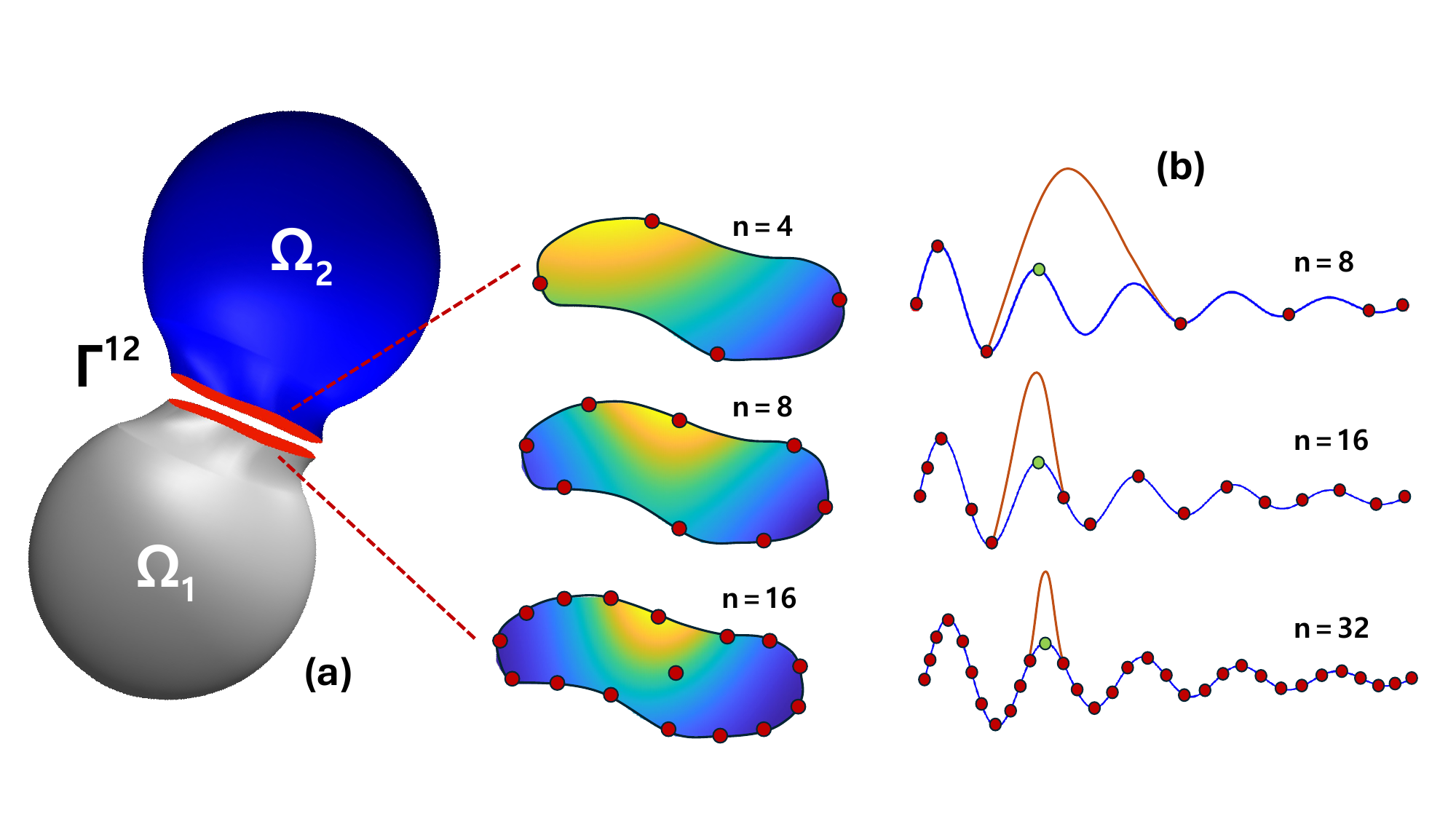}
	\caption{Schematic of mortar nodes and mortar functions. (a) Mortar nodes are placed on a 2D contact interface, $\Gamma_{12}$, between two grain grids, $\Omega_1$ and $\Omega_2$, by treating them as like-charged particles and computing their equilibrium configuration. This configuration is shown for $n\!=\!4$, 8, and 16 mortar nodes. The mortar function associated with one of the nodes is highlighted by the heat map. (b) Same concept as plot (a) but shown for a 1D interface with $n\!=\!8$, 16, and 32. Notice that mortar functions reach a maximum at their corresponding nodes and have localized support.}
\label{fig:mortar}
\end{figure}
The goal of hPLMM is to compute an approximate solution to Eq.\ref{eq:Linear-elastic}. It begins by restricting Eq.\ref{eq:Linear-elastic} onto $\Omega^{g_i}$:
\begin{equation} \label{eq:loc_orig}
	\nabla \cdot \left(\bb{C} :\, \nabla^s \bs{u}^f_{g_i} \right) = \bs{f}
	\qquad s.t. \qquad
	\begin{cases}
	  \bs{u}^f_{g_i} \cdot \bs{n}\:|_{\Gamma^{g_i}_{c_j}}=h_{g_i c_j 1}(\bs{x})  \\
	  \bs{u}^f_{g_i} \cdot \bs{m}\:|_{\Gamma^{g_i}_{c_j}}=h_{g_i c_j 2}(\bs{x})  \\
	\end{cases}
\end{equation}
where $h_{g_i c_j 1}(\bs{x})$ and $h_{g_i c_j 2}(\bs{x})$ correspond to the \textit{unknown} normal and tangential displacements on $\Gamma^{g_i}_{c_j}$, respectively. The position vector $\bs{x}$ indicates these are \textit{functions} defined on $\Gamma^{g_i}_{c_j}$. The BC on the void-solid interface $\Gamma^{g_i}_w$ is zero-stress, i.e., $\smash{\bs{\sigma}^f_{g_i} \cdot \bs{n}\!=\!\bs{0}}$. To make progress, hPLMM approximates $h_{g_i c_j 1}(\bs{x})$ and $h_{g_i c_j 2}(\bs{x})$ with the \textit{localization assumption}:
\begin{equation} \label{eq:closure_hPLMM}
	h_{g_i c_j d}(\bs{x}) \approx 
	    \underset{m_n\in M^{c_j}}{\sum}
	    \overset{D}{\underset{\gamma=1}{\sum}} \,
	    u^o_{g_i m_n \gamma}\, \eta^{(d)}_{m_n \gamma}(\bs{x})
	    \,,\qquad d=1,\ldots,D	
\end{equation}
which is a finite sum of \textit{mortar functions} $\eta^{(d)}_{m_n \gamma}(\bs{x})$, defined below, multiplied by scalar coefficients $u^o_{g_i m_n \gamma}$ called \textit{coarse-scale displacement} unknowns, to be determined from a global problem in Section \ref{sec:local_global}. Each mortar function corresponds to a \textit{mortar node}, a point on the interface that we shall define shortly, labeled with the index $m_n$ (see Fig.\ref{fig:mortar}). The coarse-scale displacement $u^o_{g_i m_n \gamma}$ is associated with node $m_n$ and corresponds to the $\gamma^{th}$ component\footnote{Here and in what follows, ``component'' is with respect to the local coordinate spanned by the unit normal $\bs{n}$ and the unit tangent $\bs{m}$ on $\Gamma^{g_i}_{c_j}$.} of the vector $\bs{u}^o_{g_i m_n}\!=\![u^o_{g_i m_n 1},\,\cdots\,,u^o_{g_i m_n D}]$. In Eq.\ref{eq:closure_hPLMM}, the outer summation is over all mortar nodes positioned on the interface $\smash{\Gamma^{c_j}}$, defined by the index set $\smash{M^{c_j}}$. The inner summation is over the coordinate directions. The symbol $\smash{\eta^{(d)}_{m_n \gamma}(\bs{x})}$ represents the $d^{th}$ component of the vectorial form of the mortar function $\bs{\eta}_{m_n \gamma}(\bs{x})$. The mortar functions form a basis on $\smash{\Gamma^{c_j}}$. 

Let us first discuss how mortar nodes are defined on each interface $\smash{\Gamma^{c_j}}$. Suppose $\Gamma^{c_j}$ consists of $\smash{N^f_{c_j}}$ fine grids (FEM nodes) whose positions $\bs{y}_i$ are contained within the set $\smash{F^{c_j}}$. Let $N^m_{c_j}$ be the number of mortar nodes sought on $\Gamma^{c_j}$ and $\bs{x}_{m_1}$, $\bs{x}_{m_2}$, ..., $\bs{x}_{m_\nu}$ with $
\nu\!=\!N^m_{c_j}$ denoting their positions. These positions are chosen from the finite set $\smash{F^{c_j}}$. The algorithm proposed by \cite{khan2024high} starts from an initial guess for $\bs{x}_{m_i}$, then treats nodes as like-charged particles that repel each other and, consequently, arrange themselves into a configuration that minimizes their collective electric potential. The computations are performed in a sequentially iterative fashion that converges within 5 steps. The starting point is to pick the initial positions of the first two mortar nodes, $\bs{x}_{m_1}$ and $\bs{x}_{m_2}$, as the farthest two points in $\smash{F^{c_j}}$. Then, provided there are more than two mortar nodes sought, the next position $\bs{x}_{m_{i+1}}$ is found from the minimization below:
\begin{equation} \label{eq:node_init}
\bs{x}_{m_{i+1}} = 
		  \underset{\bs{y} \in F^{c_j}}
         {\arg\!\min} \sum_{j=1}^i \frac{1}{\Vert \bs{y} - \bs{x}_{m_j} \Vert}
\end{equation}
Once all nodal positions are determined in this way, we sweep through the nodes again in random order and update their positions. This process is repeated until no further changes in nodal positions is observed. To update a node's position, we hold all other nodes fixed, and solve a minimization similar to Eq.\ref{eq:node_init}, except this time, the summation is over \textit{all} the nodes (upper limit $\nu$ instead of $i$) excluding the one being updated. See \cite{khan2024high} for full details.

We next proceed to defining the mortar functions $\bs{\eta}_{m_n \gamma}(\bs{x})$. Here, we consider two types of mortar functions: (1) Gaussian, and (2) Algebraic. Gaussian mortars were proposed by \cite{khan2024high} and are defined as follows:
\begin{equation} \label{eq:gauss}
	\eta^{(d)}_{m_i \gamma}(\bs{x}) = 
	\begin{cases}
	\frac{\alpha_i(\bs{x})}{\sum_{j=1}^{\nu} \alpha_j(\bs{x})}, \quad d\!=\!\gamma \\
	0, \qquad\qquad d\!\neq\!\gamma \\
	\end{cases}
    \,,\qquad
    \alpha_i(\bs{x}) = \textrm{exp}\left(\dfrac{-\Vert \bs{x} - \bs{x}_{m_i}\Vert^2}
                 {\beta\,\underset{j\in M^{c_j}}
                 {\min} \Vert \bs{x}_{m_i} - \bs{x}_{m_j}\Vert^2}\right)
\end{equation}
where $\bs{x}_{m_i}$ and $\bs{x}_{m_j}$ are mortar-node positions and recall $\nu\!=\!N^m_{c_j}$. The $\gamma^{th}$ component of $\bs{\eta}_{m_i \gamma}(\bs{x})$ attains a maximum of $\sim$1 at node $m_i$ and $\sim$0 at all other nodes. The other components are identically zero. Gaussian mortars satisfy partition of unity on the interface $\Gamma^{c_j}$, i.e., ${\sum}\,\bs{\eta}_{m_i \gamma}(\bs{x})=\bs{e}_\gamma$, where summation is over $m_i\in M^{c_j}$ and $\gamma$, and $\bs{e}_\gamma$ is the unit vector in the $\gamma$ direction. The user-defined parameter $\beta$ controls the spread of the mortar function. We set $\beta\!=\!4$ based on the analysis in \ref{app:gauss_supp}. In \cite{khan2024high}, a similar mortar, called \textit{Fickian}, was proposed that is devoid of certain artifacts Gaussian mortars have in pathological interfaces. They are omitted here as such artifacts are rare and Eq.\ref{eq:gauss} is simpler.

In this work, we propose Algebraic mortars as an alternative to Gaussian (or Fickian) mortars that do not possess a user-defined parameter like $\beta$. They are built by numerically solving the $D-1$-dimensional form of Eq.\ref{eq:Linear-elastic} on $\Gamma^{c_j}$:
\begin{equation} \label{eq:algebra}
	\nabla_{\|}\cdot \left( \bb{C} : \nabla_{\|}^s \bs{\eta}_{m_i \gamma} \right) = \bs{0}
	\qquad s.t. \qquad
	\eta^{(d)}_{m_i \gamma}(\bs{x}_{m_i}) = \delta_{d\gamma}\,, \quad  
	\eta^{(d)}_{m_{i} \gamma}(\bs{x}_{m_{j\neq i}}) = 0
\end{equation}
where $\nabla_{\|}$ and $\nabla_{\|}\cdot$ are the directional gradient and divergence operators tangent to the interface. Note that all components of the mortar function $\bs{\eta}_{m_i\gamma}$ are zero at all mortar nodes except the $\gamma^{th}$ component at node $m_i$, positioned at $\bs{x}_{m_i}$. These constraints along with zero-stress elsewhere on $\partial\Gamma^{c_j}$ constitute the BCs used to solve Eq.\ref{eq:algebra}. In MsFE, similar closure BCs are used to construct basis functions \cite{castelletto2017mult}. Unlike Gaussian mortars, here, $\eta^{(d)}_{m_i\gamma}(\bs{x})\neq\!0$ for $d\!\neq\!\gamma$ in general. Unlike Eq.\ref{eq:gauss}, where partition of unity is imposed explicitly, it is guaranteed automatically in Eq.\ref{eq:algebra} by superposition.

We conclude with a few useful definitions. We denote the span of $\bs{\eta}_{m_n\gamma}(\bs{x})$ on $\Gamma^{c_j}$ by the function space $\mathcal{M}^{c_j}$. The symbol $N^m$ is the total number of mortar nodes defined on \textit{all} interfaces. The set $M^{g_i}\! =\! \{ m_n\!\in\! M^{c_j} \,|\, c_j\!\in\! C^{g_i} \}$ contains the indices of all mortar nodes positioned on the grain-grid boundary $\partial\Omega^{g_i}$. Given a mortar node, $m_n$, the corresponding contact interface hosting the node is given by the mapping $c_j\!=\!\Lambda(m_n)$. A corollary follows that $M^{c_j}\!=\!\{ m_n \, \vert \, \Lambda\left(m_n\right)\!=\!c_j\}$. In the next section, we show how mortars are used by hPLMM to construct an approximate solution to Eq.\ref{eq:Linear-elastic}.

\subsection{Local and global problems} \label{sec:local_global}
The starting point in hPLMM is to approximate the local solution $\smash{\bs{u}^f_{g_i}}$ in Eq.\ref{eq:loc_orig} as the superposition of a number of numerically constructed \textit{shape functions},\footnote{In prior works, these were termed \textit{basis functions}. We deviate from this nomenclature to avoid conflict with \textit{basis vectors} defined later.} $\smash{\bs{\varphi}^f_{g_i m_n d}}$, and one \textit{correction function}, $\smash{\widetilde{\bs{\varphi}}^f_{g_i}}$, both defined on $\Omega^{g_i}$:
\begin{equation} \label{eq:recon}
	\bs{u}^f_{g_i} = \widetilde{\bs{\varphi}}^f_{g_i} + 
	                \underset{m_n\in M^{c_j}}
	                {\sum}\overset{D}{\underset{d=1}{\sum}} 
	                u^o_{g_i m_n d}\, \bs{\varphi}^f_{g_i m_n d}
\end{equation}
The shape function, $\smash{\bs{\varphi}^f_{g_i m_n d}}$ is computed on $\Omega^{g_i}$ from solving: 
\begin{equation} \label{eq:basis}
	\nabla \cdot \left(\bb{C} :\, \nabla^s \bs{\varphi}^f_{g_i m_n d} \right) = 0
	\qquad s.t. \qquad
	\begin{cases}
	  \bs{\varphi}^f_{g_i m_n d} \cdot \bs{n}\:|_{\Gamma^{g_i}_{c_j}}
	              =\delta_{c_k c_j}\,\bs{\eta}_{m_n d}(\bs{x})\cdot\bs{n}
	              =\delta_{c_k c_j}\,\eta^{(1)}_{m_n d}(\bs{x})  \\
	  \bs{\varphi}^f_{g_i m_n d} \cdot \bs{m}\:|_{\Gamma^{g_i}_{c_j}}
	              =\delta_{c_k c_j}\,\bs{\eta}_{m_n d}(\bs{x})\cdot\bs{m}
	              =\delta_{c_k c_j}\,\eta^{(2)}_{m_n d}(\bs{x})  \\
	\end{cases}
\end{equation}
and the correction function is computed on $\Omega^{g_i}$ from solving:
\begin{equation} \label{eq:cor}
	\nabla \cdot \left(\bb{C} :\, \nabla^s \widetilde{\bs{\varphi}}^f_{g_i} \right) = \bs{f}
	\qquad s.t. \qquad
	\begin{cases}
	  \widetilde{\bs{\varphi}}^f_{g_i} \cdot \bs{n}\:|_{\Gamma^{g_i}_{c_j}}=0 \\
	  \widetilde{\bs{\varphi}}^f_{g_i} \cdot \bs{m}\:|_{\Gamma^{g_i}_{c_j}}=0 \\
	\end{cases}
\end{equation}

In Eq.\ref{eq:basis}, $\delta_{c_k c_j}$ is the Kronecker delta and $c_k\!=\!\Lambda(m_n)$. For the more intuitive Gaussian mortars described by Eq.\ref{eq:gauss}, if $d\!=\!1$, $\smash{\bs{\varphi}^f_{g_i m_n d}}$ is obtained by setting the normal displacement on $\Gamma^{g_i}_{c_k}$ equal to $\smash{\eta^{(1)}_{m_n d}(\bs{x})}$, while setting the tangential displacement on $\smash{\Gamma^{g_i}_{c_k}}$ and all normal/tangential displacements on $\Gamma^{g_i}_{c_j} \, \forall c_{j\neq k} \!\in\! C^{g_i}$ equal to zero. If $d\!=\!2$, only the tangential displacement on $\Gamma^{g_i}_{c_k}$ is set to $\smash{\eta^{(2)}_{m_n d}(\bs{x})}$, while all other displacements on all interfaces are set to zero. This is schematized by Fig.\ref{fig:schem}e, where three shape functions for the grain grid in Fig.\ref{fig:schem}d are shown, corresponding to the three mortar nodes (yellow dots; $N^m_{c_j}\!=\!3$) on the highlighted interface (green). For the Algebraic mortars described by Eq.\ref{eq:algebra}, non-zero contributions exist for both the normal and tangential BCs of Eq.\ref{eq:basis}. We remark that Eq.\ref{eq:basis} requires solving $\sum_{c_j\in C^{g_i}}\! N^m_{c_j}\!\times\!D$ local problems on $\Omega^{g_i}$, whereas Eq.\ref{eq:cor} requires solving only one local problem. These computations are fully decoupled spatially across all grain grids and are, thus, amenable to parallelism.

After the shape and correction functions are computed in a pre-processing step, the coarse displacements $u^o_{g_i m_n d}$ are found by solving a global problem that consists of a force balance and a kinematic constraint on all $\Gamma^{c_j}$:
\begin{subequations}\label{eq:global}
\begin{equation} \label{eq:mom_bal}
	\langle\bs{t}_{g_1}\,,\bs{\eta}_{m_n \gamma}\rangle = 
	\langle\bs{t}_{g_2}\,,\bs{\eta}_{m_n \gamma}\rangle \qquad \forall m_n, \gamma
\end{equation}
\begin{equation} \label{eq:stick}
	u^o_{g_1 m_n \gamma} = u^o_{g_2 m_n \gamma} \qquad\quad\qquad \forall m_n, \gamma
\end{equation}
\end{subequations}
Eq.\ref{eq:mom_bal} states that the L$_2$-inner product between the fine-grid traction $\bs{t}$ computed on one side of the interface, namely $\Gamma^{g_{1}}_{c_j}\!\subset\!\partial\Omega^{g_1}$, is equal to that computed on the other side, namely on $\Gamma^{g_{2}}_{c_j}\!\subset\!\partial\Omega^{g_2}$. Eq.\ref{eq:stick} states that the coarse displacement $u^o_{g_1 m_n \gamma}$ defined on $\Gamma^{g_{1}}_{c_j}$ is equal to $u^o_{g_2 m_n \gamma}$ defined on $\Gamma^{g_{2}}_{c_j}$. Using Eq.\ref{eq:recon} and Hooke's law to compute the fine-grid tractions $\bs{t}_{g_1}$ and $\bs{t}_{g_2}$, followed by their substitution into Eq.\ref{eq:global}, yields the following linear system for $u^o_{g_2 m_n d}$:
\begin{equation} \label{eq:residual}
  \Re(\bb{U})=\Re\left([u^o_{g_i m_n d}]_{2N^m D\times 1}\right)=0
\end{equation}
Eq.\ref{eq:residual} is the \textit{global problem}, where $\Re(\bb{U})$ is the residual of size $2N^mD\!\times\!1$; much smaller than Eq.\ref{eq:ls_all} defined on the fine grid. We refer to the approximate fine-scale solution obtained by substituting the $u^o_{g_i m_n \gamma}$ from Eq.\ref{eq:residual} into Eq.\ref{eq:recon} as the \textit{first-pass} solution of hPLMM. We denote the maximum number of mortar nodes used per interface by $n\!=\!\max_{c_j} N^m_{c_j}$.

\subsection{Error control} \label{sec:err_control}

The localization assumption in Eq.\ref{eq:closure_hPLMM} introduces errors in the first-pass solution concentrated near $\Gamma^{c_j}$. To remove these errors, one may either increase the number of mortar nodes per interface, $n$, or adopt an iterative scheme proposed by \cite{khan2024high}, where the BCs of the correction problem (Eq.\ref{eq:cor}) are successively updated using the fine-scale solution from the previous iteration. The procedure requires solving the following \textit{contact problem} on each contact grid $\Omega^{\zeta_k}$:
\begin{equation} \label{eq:dual}
	\nabla \cdot \left(\bb{C} :\, \nabla^s \bs{u}^{f,\omega}_{\zeta_k} \right) = \bs{f}
	\qquad s.t. \qquad
	\begin{cases}
	  \bs{u}^{f,\omega}_{\zeta_k} \cdot \bs{n}\:|_{\Gamma^{\zeta_k}_g}=R_{\Gamma^{\zeta_k}_g}\left[\bs{u}^{f,\omega-1}\right]\cdot\bs{n} \\
	  \bs{u}^{f,\omega}_{\zeta_k} \cdot \bs{m}\:|_{\Gamma^{\zeta_k}_g}=R_{\Gamma^{\zeta_k}_g}\left[\bs{u}^{f,\omega-1}\right]\cdot\bs{m}
	\end{cases}
\end{equation}
for all $\zeta_k$. The symbol $\omega$ denotes the iterations index, and $\smash{R_{\Gamma^{\zeta_k}_g}[\cdot]}$ is an operator that restricts the previous fine-scale solution $\bs{u}^{f, \omega-1}$ onto $\smash{\Gamma^{\zeta_k}_g}$ ($\smash{=\!\cup_i \Gamma^{\zeta_k}_{g_i}}$). The initial guess is chosen to be the first-pass solution from Section \ref{sec:local_global}. In Eq.\ref{eq:dual}, the local BC on the void-solid interface $\Gamma^{\zeta_k}_w$ (not shown) is zero-stress. If $\partial\Omega^{\zeta_k}$ intersects $\Gamma^{ex}$, then Eq.\ref{eq:dual} inherits the global BCs in Eq.\ref{eq:BC}. Notice Eq.\ref{eq:dual} is fully decoupled across all contact grids and thus, is amenable to parallelism. For details, see \cite{khan2024high}. The smoother described in Section \ref{sec:local_smoother} is an algebraic translation of the above iterative scheme. If errors are reduced by increasing $n$, instead of $\omega$, then the cost of offline (pre-processing) computations grows as more shape functions are built via Eq.\ref{eq:basis}. By contrast, if $\omega$ is increased, only online (run-time) computations grow in cost.

\section{Multiscale preconditioner} \label{sec:multi_precon}
The hPLMM preconditioner, $\mathrm{M}$, proposed below is an algebraic translation of the steps outlined in Section \ref{sec:hplmm} and generalizes the low-order PLMM preconditioner by \cite{mehmani2023precond, li2024phase}. Section \ref{sec:precon_structure} discusses the overall structure of $\mathrm{M}$, followed by Sections \ref{sec:global_precon}-\ref{sec:local_smoother} that discuss its two building blocks $\mathrm{M_G}$ and $\mathrm{M_L}$. Below ``hPLMM'' refers to the preconditioner.
 
\subsection{Preconditioner structure} \label{sec:precon_structure}
The hPLMM preconditioner, $\mathrm{M}$, is a multiplicative combination of a global (or coarse) preconditioner, $\mathrm{M_G}$, and a local (of fine-grid) smoother, $\mathrm{M_L}$, as follows. Notice that $\mathrm{M_G}$ is applied first, followed by $\mathrm{M_L}$:
\begin{equation} \label{eq:multi_precond}
	\mathrm{M}^{-1} = \mathrm{M_G^{-1}} + \mathrm{M_L^{-1}}(\mathrm{I} - \hat{\mathrm{A}} \mathrm{M_G^{-1}})
\end{equation}
M\ts{G} attenuates low-frequency error modes, while M\ts{L} removes high-frequency errors. The smoother itself is formulated here as the repeated application of a \textit{base smoother}, M\ts{$l$}, in multiplicative fashion, over $n_{st}$ stages:
\begin{equation} \label{eq:local_smoother}
 	\mathrm{M_L^{-1}} = \sum_{i=1}^{n_{st}} \text{M}_l^{-1} \prod_{j=1}^{i-1}\,(\mathrm{I} - \hat{\mathrm{A}}\text{M}_l^{-1})
 \end{equation}
Any black-box base smoother, like Gauss-Seidel (M\ts{GS}) or incomplete LU-factorization (M\ts{ILU($k$)}), can be used. However, as we will find in Section \ref{sec:results}, black-box smoothers tend to cause the solvers to converge slowly or even stagnate. A better smoother, compatible with the proposed $\mathrm{M_G}$ below, is the additive-Schwarz preconditioner by \cite{li2024phase} referred to as the \textit{contact-grain} smoother, M\ts{CG}, reviewed in Section \ref{sec:local_smoother}. In Section \ref{sec:results}, we benchmark M\ts{CG} against M\ts{ILU($k$)}.

\subsection{Global preconditioner}\label{sec:global_precon}
The global preconditioner M\ts{G} is formulated as follows:
 \begin{equation}\label{eq:mg_struct}
 	\mathrm{M}^{-1}_{\mathrm{G}} = \mathrm{\hat{P}}\, (\mathrm{\hat{R}} \hat{\mathrm{A}} \mathrm{\hat{P}})^{-1} \mathrm{\hat{R}} 
 \end{equation}
where $\mathrm{\hat{P}}$ and $\mathrm{\hat{R}}$ are referred to as the \textit{prolongation} and \textit{restriction} matrices, respectively, which satisfy $\mathrm{\hat{R}}\!=\!\mathrm{\hat{P}}^{\top}$. The matrix $\mathrm{\hat{P}}$ is itself a product of three matrices:
\begin{equation} \label{eq:WQP}
 	\mathrm{\hat{P}} = \mathrm{W} \mathrm{Q} \mathrm{P}
\end{equation}
where W is a \textit{permutation matrix}, Q is a \textit{reduction matrix}, and P a \textit{reduced-prolongation matrix}. In \cite{li2024phase}, these matrices were formulated for PLMM. In what follows, we generalize those formulations to hPLMM.
\\

\textbf{Permutation (W).} The matrix, W, is square and consists of only 1s and 0s. It is unitary (i.e., $\mathrm{W}\mathrm{W}^\top\! =\! \mathrm{I}$) and upon right-multiplying a matrix, shuffles the columns of that matrix. The shuffling is done in such a way that the fine-grid entries associated with each grain grid $\Omega^{g_i}$ and contact interface $\Gamma^{c_j}$ are grouped together, in accordance with the domain decomposition described in Section \ref{sec:domdecomp}. Concretely, applying W to Eq.\ref{eq:ls_all} yields:
\begin{subequations} \label{eq:ls_trans}
\begin{equation} \label{eq:permute}
	\underbrace{\mathrm{W}^\top \hat{\mathrm{A}} \mathrm{W}}_{\mathrm{A}} \underbrace{\mathrm{W}^\top \hat{x}}_{x} = \underbrace{\mathrm{W}^\top \hat{b}}_{b} \quad \Rightarrow \quad \mathrm{A}x = b
\end{equation}
where the permuted $\mathrm{A}$, $b$, and $x$ have the following block structures:
\begin{equation} \label{eq:permute_blk_1}
	\mathrm{A} = 
	\begin{bmatrix}
	\mathrm{A}^g_g & \mathrm{A}^g_c \\
	\mathrm{A}^c_g & \mathrm{A}^c_c
	\end{bmatrix}
	\qquad
	b =
	\begin{bmatrix}
	b^g \\
	b^c
	\end{bmatrix}
	\qquad
	x =
	\begin{bmatrix}
	x^g \\
	x^c
	\end{bmatrix}
\end{equation}
\begin{equation} \label{eq:permute_blk_2}
	\mathrm{A}^g_g =
	\begin{bmatrix}
	\mathrm{A}^{g_1}_{g_1} & \cdots & \mathrm{O} \\
	\vdots & \ddots & \vdots \\
	\mathrm{O} & \cdots & \mathrm{A}^{g_{N^g}}_{g_{N^g}}
	\end{bmatrix}_{N_g^f \times N_g^f}
	\qquad
	\begin{matrix}
	\mathrm{A}^c_g = [ \mathrm{A}^{c_i}_{g_j} ]_{N_c^f \times N_g^f} \\
	\mathrm{A}^g_c = [ \mathrm{A}^{g_i}_{c_j} ]_{N_g^f \times N_c^f} \\
	\mathrm{A}^c_c = [ \mathrm{A}^{c_i}_{c_j} ]_{N_c^f \times N_c^f}
	\end{matrix}
	\qquad\quad
	\begin{matrix}
	b^g = [b^{g_i}]_{N_g^f \times 1} \\
	b^c = [b^{c_i}]_{N_c^f \times 1}
	\end{matrix}
\end{equation}
\end{subequations}
The super/subscripts $g_i$ and $c_j$ represent entries or blocks that belong to either $\Omega^{g_i}$ or $\Gamma^{c_j}$, respectively. $N_{g_i}^f$ and $N_{c_j}^f$ denote the number of fine-scale unknowns associated with $\Omega^{g_i}$ and $\Gamma^{c_j}$, respectively. By construction, $\smash{N_g^f \!=\! \sum_i N_{g_i}^f}$ and $\smash{N_c^f \!=\! \sum_j N_{c_j}^f}$, where $N^g$ is the total number of grain grids.
The matrix $\mathrm{A}^g_g$ is square and block-diagonal, with square blocks $\mathrm{A}^{g_i}_{g_i}$, while $\mathrm{A}^g_c$ and $\mathrm{A}^c_g$ are thin and rectangular. Given the Galerkin FEM discretization used herein on the fine grid, and the self-adjoint weak form of Eq.\ref{eq:Linear-elastic}, $\mathrm{A}^g_c\! =\! (\mathrm{A}^c_g)^\top$ and $\mathrm{A}^{g_i}_{c_j}\! =\! (\mathrm{A}^{c_j}_{g_i})^\top$ hold. Building W is trivial and thus cheap.\\
\\
\textbf{Reduction matrix (Q).}
The reduction matrix Q constitutes the key difference between PLMM and hPLMM, as it is where mortar functions become encoded algebraically. The job of Q is to perform a \textit{weighted} column-sum, when right-multiplying a matrix, over all the entries associated with each contact interface $\Gamma^{c_j}$ separately. For the linear-elasticity PDE given by Eq.\ref{eq:Linear-elastic}, this summation is performed on a per-coordinate-direction basis. The weights correspond precisely to the mortar functions $\eta^{(d)}_{m_n \gamma}(\bs{x})$ defined in Section \ref{sec:mortar}. Concretely, Q is given by:
\begin{equation} \label{eq:reduced_defs}
	\mathrm{Q} =
	\begin{bmatrix}
		\mathrm{I}_{N_g^f \times N_g^f} & \mathrm{O} \\
		\mathrm{O} & \mathrm{Q}^o
	\end{bmatrix}
	\qquad
	\mathrm{Q^o} =
	\begin{bmatrix}
		\mathrm{M}^{m_1} & & \mathrm{O} \\
		& \ddots & \\
		\mathrm{O} & & \mathrm{M}^{m_N}
	\end{bmatrix}_{N_c^f\times N^o_m}
	\;
	\mathrm{M}^{m_i} = 
	\begin{bmatrix}
		\eta^{(1)}_{m_i 1} & \ldots & \eta^{(1)}_{m_i D} \\
		\vdots & \ddots & \vdots \\
		\eta^{(D)}_{m_i 1} & \ldots & \eta^{(D)}_{m_i D}
	\end{bmatrix}_{\,N^f_{c_k} \times D}
\end{equation}
where $N^o_m\!=\!N^mD$ and $c_k\!=\!\Lambda(m_i)$. 
Recall, $N^m$ is the total number of mortar nodes and $D$ the problem dimension. Notice Q is rectangular\footnote{In \cite{li2024machine}, it was mistakenly stated that Q for PLMM is square. It is not, as is obvious from the Q$^o$ block.} and block-diagonal. The (1,1)-block is the $N_g^f\! \times\! N_g^f$ identity matrix, and the (2,2)-block, $\mathrm{Q}^o$, a block-diagonal matrix itself. Each block $\mathrm{M}^{m_i}$ of $\mathrm{Q}^o$ is made up of mortar functions $\eta^{(d)}_{m_n \gamma}(\bs{x})$, understood to represent column vectors defined on the fine grids of their corresponding contact interfaces. Another way to think about $\mathrm{M}^{m_i}$ is a $D$-column matrix, with the $\gamma^{th}$ column representing $\bs{\eta}_{m_i \gamma}(\bs{x})$ (a $D$-dimensional vector) stretched into a column vector.

Using Q, the permuted system in Eq.\ref{eq:permute} can be reduced as follows:
\begin{subequations} \label{eq:reduced}
\begin{equation} \label{eq:reduced_sys}
	\mathrm{A}x = b\;,\quad x\simeq \mathrm{Q}x_M \quad \Rightarrow \quad \mathrm{Q}^\top \mathrm{A} \mathrm{Q} x_M = \mathrm{Q}^\top b \quad \Rightarrow \quad \mathrm{A_M} x_M  = b_M
\end{equation}
where the reduced matrix A\ts{M} now has the following block structure:
\begin{equation} \label{eq:redu_orig}
	\mathrm{A_M} =
	\begin{bmatrix}
		\mathrm{A}_g^g & \bar{\mathrm{A}}_m^g \\
		\bar{\mathrm{A}}_g^m & \bar{\mathrm{A}}_m^m
	\end{bmatrix}
	\qquad\qquad
	\begin{matrix}
	\bar{\mathrm{A}}^m_g = [ \bar{\mathrm{A}}^{m_i}_{g_j} ]_{N_m^o \times N_g^f} \\
	\bar{\mathrm{A}}^g_m = [ \bar{\mathrm{A}}^{g_i}_{m_j} ]_{N_g^f \times N^o_m} \\
	\bar{\mathrm{A}}^m_m = [ \bar{\mathrm{A}}^{m_i}_{m_j} ]_{N^o_m \times N^o_m}
	\end{matrix}
\end{equation}
\end{subequations}
The reduced system in Eq.\ref{eq:reduced_sys} is slightly smaller than Eq.\ref{eq:permute}, provided $N^o_m\!<\! N^f_c$, which implies that the number of mortar nodes per interface is smaller than the number of fine grids per interface. Applying Q in Eq.\ref{eq:reduced_sys} simultaneously imposes the localization assumption in Eq.\ref{eq:closure_hPLMM} and the force balance and kinematic constraint in Eq.\ref{eq:global} across all $\Gamma^{c_j}$. Building Q requires constructing $\smash{\bs{\eta}_{m_i \gamma}}$, which is inexpensive, especially for Gaussian mortars. If the number of mortar nodes per interface is $n\!=\!1$, Eq.\ref{eq:reduced_defs} reduces to PLMM with $\smash{\eta^{(d)}_{m_i \gamma}\!=\!\bs{\delta}_{d\gamma}}$, an all-one (if $d\!=\!\gamma$) or all-zero (if $d\!\neq\!\gamma$) vector.\footnote{In \cite{mehmani2023precond, li2024phase, li2024machine} focused on PLMM, the block $\mathrm{M}^{m_i}$ was concatenated out of a number of $D\!\times\! D$ identity matrices. The tacit assumption there was the fine-scale unknowns in the permuted system are ordered like $x$, $y$, $z$ component for FEM node 1, then $x$, $y$, $z$ for FEM node 2, and so on. Here, we assumed a different ordering that makes the presentation of $\mathrm{M}^{m_i}$  easier, namely, the $x$ components of \textit{all} nodes comes first, then the $y$ component of all nodes come next, and so on. This ordering has no impact on any of the calculations and results that follow. It just clarifies the presentation.}
\\

\textbf{Remark.} The block $\mathrm{M}^{m_i}$ can be built fully algebraically, without resorting to solving Eq.\ref{eq:algebra} geometrically on $\Gamma^{c_j}$, with $c_j\!=\!\Lambda(m_i)$. Concretely, to compute the $\gamma^{th}$ column of $\mathrm{M}^{m_i}$, one must solve a local problem involving the matrix $\mathrm{A}^{c_j}_{c_j}$ subject to Dirichlet constraints on the mortar nodes. Namely, the $\gamma^{th}$ component of the displacement at node $m_i$ is one, and all other nodal displacements are zero. The approach involves four steps: (1) Remove all rows corresponding to the 0/1-constraints on mortar nodes from $\mathrm{A}^{c_j}_{c_j}$; (2) Extract the column corresponding to the 1-constraint from the resulting matrix in Step 2, and call it $e_1$; (3) Remove all columns corresponding to the 0/1-constraints from the reduced matrix in Step 2, and call it $\mathrm{\tilde{A}}^{c_j}_{c_j}$; (4) Compute $-(\mathrm{\tilde{A}}^{c_j}_{c_j})^{-1}e_1$ and expand it by adding back the 0/1 values at mortar nodes.
\\

\textbf{Reduced prolongation matrix (P).}
The reduced prolongation, P, defines a coarse space spanned by its columns, wherein a good approximation to the solution, $x_M$, to the reduced system in Eq.\ref{eq:reduced_sys} exists. It is defined as follows:
\begin{subequations} \label{eq:all_prolong}
 \begin{equation} \label{eq:bc_matrix}
 	\mathrm{P} =
 	\begin{bmatrix}
		\mathrm{B} & \mathrm{C} \\
		\mathrm{I} & \mathrm{O}
	\end{bmatrix}_{(N^f_g+N_m^o)\times(N_m^o+N^g)}
	\quad
 	\mathrm{B} =
	\begin{bmatrix}
		p_1^{g_1} & p_2^{g_1} & \cdots & p_{N}^{g_1} \\
		p_1^{g_2} & p_2^{g_2} & \cdots & p_{N}^{g_2} \\
		\vdots & \vdots & \ddots & \vdots \\
		p_1^{g_M} & p_2^{g_M} & \cdots & p_{N}^{g_M}
	\end{bmatrix}_{N^f_{g} \times N^o_m}
	\quad 
	\mathrm{C} = 
	\begin{bmatrix}
		c^{g_1} & & & \mathrm{O}\\
		 & c^{g_2} & &\\
		 & & \ddots &\\
		\mathrm{O} & & & c^{g_M}
	\end{bmatrix}_{N^f_g \times N^g}
 \end{equation}
with B referred to as the \textit{basis matrix}, and C as the \textit{correction matrix}, composed of the following vectors:
 \begin{align}
    &c^{g_i} = (\mathrm{A}_{g_i}^{g_i})^{-1} b^{g_i} \label{eq:col_cg} \\
 	&p_k^{g_i} =
 	\begin{cases}
 	          -(\mathrm{A}_{g_i}^{g_i})^{-1} 
 	           \bar{\mathrm{A}}_{m_j}^{g_i} \mathrm{R}_m^{m_j} e_k,
  	               &g_i \in G^{c_t} \;\; 
 	                c_t\!=\!\Lambda(m_j) \;\; \;\;
 	                m_j=\lceil k/D \rceil
              \\
 	          O, 
 	               &g_i \notin G^{c_t} \;\;
 	                c_t\!=\!\Lambda(m_j) \;\; \;\;
 	                m_j=\lceil k/D \rceil
 	\end{cases}  \label{eq:col_pk} \\ 
 	&e_k = 
 	       [ 0,\, \cdots, 
			      \underbrace{0, 1, 0,}_{k-1,\, k,\, k+1} 
				   \cdots, 0]^\top_{\,N^o_m \times 1}
				   \label{eq:col_ek}
 \end{align}
\end{subequations}

For sake of brevity in Eq.\ref{eq:bc_matrix}, the dummy subscripts $N$ and $M$ have been used to equal $N^o_m$ and $N^g$ respectively. The column vectors $p^{g_i}_k$ and $c^{g_i}$ denote a \textit{shape vector} and a \textit{correction vector}, respectively, both defined on $\Omega^{g_i}$. The index \textit{k} is mapped to its corresponding mortar node via $m_j\!=\!\lceil k/D \rceil$, and $m_j$ is mapped to its corresponding contact interface hosting the node via $c_t\!=\!\Lambda(m_j)$. Finally, recall the set $G^{c_t}$ contains the only two grain grids that share $\Gamma^{c_t}$. In Eq.\ref{eq:col_ek}, $e_k$ is the unit vector that contains 1 in its $k^{th}$ entry and zero elsewhere. $\smash{\mathrm{R}_m^{m_j}}$ is a \textit{contraction matrix} defined as:
\begin{equation} \label{eq:contract}
   \mathrm{R}_m^{m_j} =
	\begin{bmatrix}
		\Delta^{m_j}_{m_1}, \Delta^{m_j}_{m_2}, \cdots, \Delta^{m_j}_{m_{N}}
	\end{bmatrix}_{D \times N^o_m}
	\qquad
	\Delta^{m_i}_{m_j} =
	\begin{cases}
		\mathrm{I}_{D \times D} \quad \text{  if }\; i=j \\
		\mathrm{O}_{D \times D} \quad \text{if }\; i\neq j
	\end{cases}
\end{equation}
Left-multiplying a $N^o_m \!\times\! 1$ vector defined on all mortar nodes (e.g., $e_k$) by $\mathrm{R}_m^{m_j}$ restricts it to a $D \!\times\! 1$ vector defined on the mortar node $m_j$. Eq.\ref{eq:col_cg} is the algebraic equivalent of solving Eq.\ref{eq:cor} for the the correction function, and Eq.\ref{eq:col_pk} is the algebraic equivalent of solving Eq.\ref{eq:basis} for the shape functions. Both local systems have the coefficient matrix $\smash{\mathrm{A}^{g_i}_{g_i}}$.

The cost of building $\mathrm{P}$ is dominated by $\mathrm{B}$, as multiple shape vectors per grain grid must be computed, whereas only one correction vector per grain grid is needed to build $\mathrm{C}$. We call each column of $\mathrm{P}$ a \textit{basis vector}, which consists of only two non-zero shape vectors. The latter is because the $k^{th}$ coarse-scale unknown, defined on node $m_j$, is shared between only two grain grids. Hence, B is sparse. And given $\smash{N^f_g\!\gg\!N^o_m}$, P is tall and skinny. To build P, $2N^o_m$ shape vectors and $N^g$ correction vectors must be computed, all spatially decoupled across grain grids and thus, amenable to parallelism. Moreover, $c^{g_i}$ is non-zero only on grain grids with a source term, $\bs{f}$ in Eq.\ref{eq:Linear-elastic}, or non-homogeneous BCs.

\subsection{Local smoother}\label{sec:local_smoother}
We next review the \textit{contact-grain} smoother M\ts{CG} proposed by \cite{li2024phase}, which we show later to provide the best pairing with the $\mathrm{M_G}$ from the last section. It is made of two additive-Schwarz preconditioners, $\mathrm{M}_\zeta$ and $\mathrm{M}_g$, as follows:
\begin{equation} \label{eq:local_smoother_CG}
	\text{M}_{\mathrm{CG}}^{-1} = \text{M}_\zeta^{-1} + \text{M}_g^{-1}(\mathrm{I} - 	\hat{\mathrm{A}} \text{M}_\zeta^{-1})
\end{equation}
We refer to $\mathrm{M}_\zeta$ as the \textit{contact-grid} smoother and to $\mathrm{M}_g$ as the \textit{grain-grid} smoother. The job of $\mathrm{M}_\zeta$ is to wipe out high-frequency errors in the interior of contact grids, which overlap $\Gamma^{c_j}$, and the job of $\mathrm{M}_g$ is to remove errors in the interior of grain grids. Both have the standard form of any additive-Schwarz preconditioner \cite{saad2003book} as follows:
\begin{equation} \label{eq:Mg_Mz}
	\text{M}_g^{-1}
	= \sum_{i=1}^{N^g} \mathrm{E}_f^{g_i} 
	(\mathrm{R}_f^{g_i} \mathrm{\hat{A}} \mathrm{E}_f^{g_i})^{-1}
	\mathrm{R}_f^{g_i}
	\quad\qquad
	\text{M}_\zeta^{-1}
	= \sum_{i=1}^{N^{\zeta}} \mathrm{E}_f^{\zeta_i} 
	(\mathrm{R}_f^{\zeta_i} \mathrm{\hat{A}} \mathrm{E}_f^{\zeta_i})^{-1}
	\mathrm{R}_f^{\zeta_i}\
\end{equation}
where matrices $\smash{\mathrm{R}_f^{g_i}}$ and $\smash{\mathrm{R}_f^{\zeta_i}}$ have 0 or 1 entries and restrict any vector they left-multiply onto $\Omega^{g_i}$ and $\Omega^{\zeta_i}$, respectively. The matrices $\smash{\mathrm{E}_f^{g_i}}$ and $\smash{\mathrm{E}_f^{\zeta_i}}$ are mere transposes of these restrictions and extend any vector they left-multiply from $\Omega^{g_i}$ and $\Omega^{\zeta_i}$ to $\Omega$, respectively.
When applied in an iterative solver, $\mathrm{M}_g$ entails solving $N^g$ decoupled systems on the grain grids, and $\mathrm{M}_\zeta$ entails solving $N^\zeta$ decoupled systems on the contact grids; all fully parallelizable. For details, see \cite{li2024phase}.
\section{Problem set} \label{sec:prob_set}

To test the proposed hPLMM preconditioner, we consider the domains in Figs.\ref{fig:problem_set1}-\ref{fig:problem_set2}. They include a small 2D disk pack (P2D) taken from \cite{mehmani2021contact}, a large 2D disk pack (P2DL) from \cite{li2024phase}, a 2D cross-section of a sandstone (S2D) from \cite{data2018sand}, a 2D continuum domain (DARCY) from \cite{khan2024high}, and a 3D bone specimen (BONE) from \cite{data2021bone}. Except DARCY, the foregoing domains correspond to solids described by binary images at the pore scale with homogeneous stiffness $\bb{C}_0$, obtained from Eq.\ref{eq:stiff_iso} with $\lambda\! =\! 8.3$ GPa and $\mu\! =\! 44.3$ GPa. The DARCY domain is described by a gray-scale image, with pixel values $\xi$ serving as stiffness multipliers, i.e., $\bb{C}\!=\!\xi\,\bb{C}_0$. The distribution of $\mathrm{ln}(\xi)$ has a mean of zero, a variance of one, and a Gaussian covariance with a spatial correlation length of 0.1. This results in $\xi\!\in$(0.05, 15) with variations in $\bb{C}$ by up to a factor of $\sim$300. To test for heterogeneity at the pore scale, we consider two variants of S2D shown in Fig.\ref{fig:problem_set1}: S2DC taken from \cite{mehmani2023precond}, and S2DH from \cite{mehmani2021contact}. In S2DC, pre-existing cracks are drawn synthetically so that some fall entirely within grain grids and others intersect contact interfaces. The effect of these cracks on stiffness is modeled through a degradation function $g$ as $g\bb{C}$, where $g\!\in$(0,1) is a continuously varying scalar function. The regularization of sharp cracks as continuous damage fields follows the approach by \cite{borden2012phase} and is detailed in \cite{mehmani2023precond}. In S2DH, $\bb{C}$ assumes one of two values: $\bb{C}_1$ representing a hard material and $\bb{C}_2$ representing a soft material. They contrast by a factor of $10^6$, and are depicted by the dark blue (hard) and light blue (soft) regions in Fig.\ref{fig:problem_set1}. The Lam\'e parameters of $\bb{C}_1$ are $\lambda\! =\! 49\times 10^3$ GPa and $\mu\! =\! 1\times 10^3$ GPa, and those of $\bb{C}_2$ are $\lambda\! =\! 49\times 10^{-3}$ GPa and $\mu\! =\! 1\times 10^{-3}$ GPa.

We decompose the pore-scale domains P2D, P2DL, S2D, and BONE using the watershed algorithm described in Section \ref{sec:domdecomp}. The decompositions of S2DC and S2DH are identical to that of S2D, which conforms to stiffness discontinuities in S2DH but not to the cracks in S2DC. For DARCY, we consider both spectral and Cartesian decompositions in Section \ref{sec:domdecomp}. The resulting grain grids for each domain are depicted as randomly colored regions in Figs.\ref{fig:problem_set1}-\ref{fig:problem_set2}, and contact grids as cyan-colored patches. Contact grids are 16 fine grids (FEM elements) wide in all domains. We subject each domain to a shear load and, separately, to a tensile load. With reference to coordinate axes in Figs.\ref{fig:problem_set1}-\ref{fig:problem_set2}, the BCs follow: In all 2D domains, shear is imposed with a downward unit displacement (negative $y$) on the left boundary ($x\!=\!0$), while keeping the right boundary ($x\!=\!L_x$) clamped. The remaining top/bottom boundaries ($y\!=\!0$ and $L_y$) are stress free. Tensile loading is identical except the left boundary is pulled leftward (negative $x$) by a unit displacement. In the 3D BONE domain, shear is imposed by pulling the left boundary ($y\!=\!L_y$) downward by a unit displacement (negative $z$) while keeping the right boundary ($y\!=\!0$) fixed. The four remaining lateral boundaries are stress free. Tensile loading is again identical, except the left boundary is pulled leftward (positive $y$) by a unit displacement.

\begin{figure} [t!]
  \centering
  \centerline{\includegraphics[scale=0.55,trim={50 20 40 30},clip]{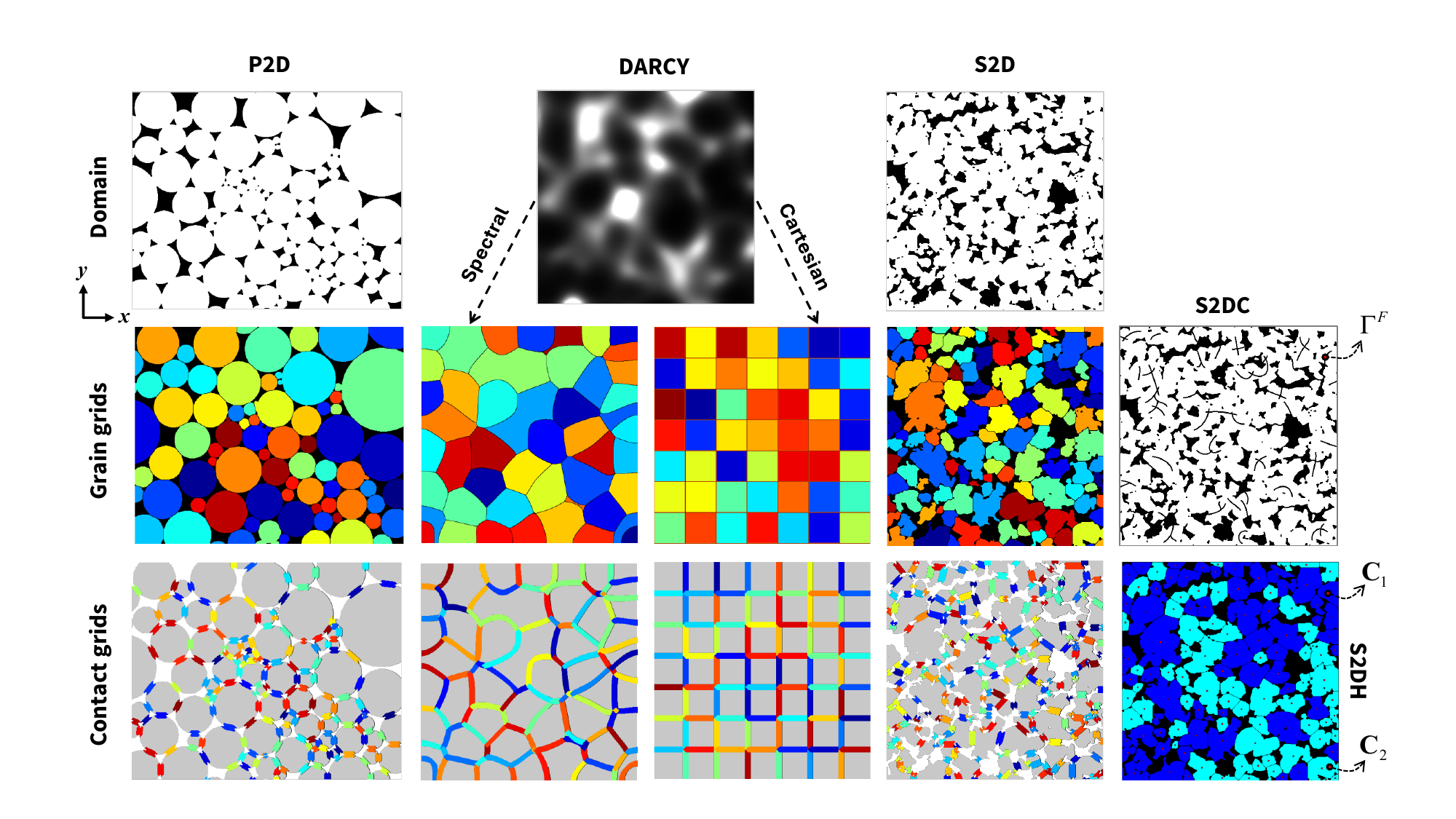}}
	\caption{Domains with complex geometries and material heterogeneity considered to test the hPLMM preconditioner. Binary images correspond to a 2D disk pack (P2D) and a 2D sandstone (S2D), whereas the gray-scale image depicts a continuum domain (DARCY) whose (exponentiated) pixel values represent a stiffness multiplier. The randomly colored regions (middle row) represent grain grids obtained from the decomposition of each domain, and the thin colored strips (bottom row) represent contact grids. The S2D domain has two variants (rightmost column): (1) S2DC with pre-existing cracks annotated by $\Gamma^F$, and (2) S2DH composed of hard (dark blue) and soft (light blue) material with high-contrasting stiffness.}
\label{fig:problem_set1}
\end{figure}

\begin{figure} [t!]
  \centering
  \centerline{\includegraphics[scale=0.55,trim={75 5 60 40},clip]{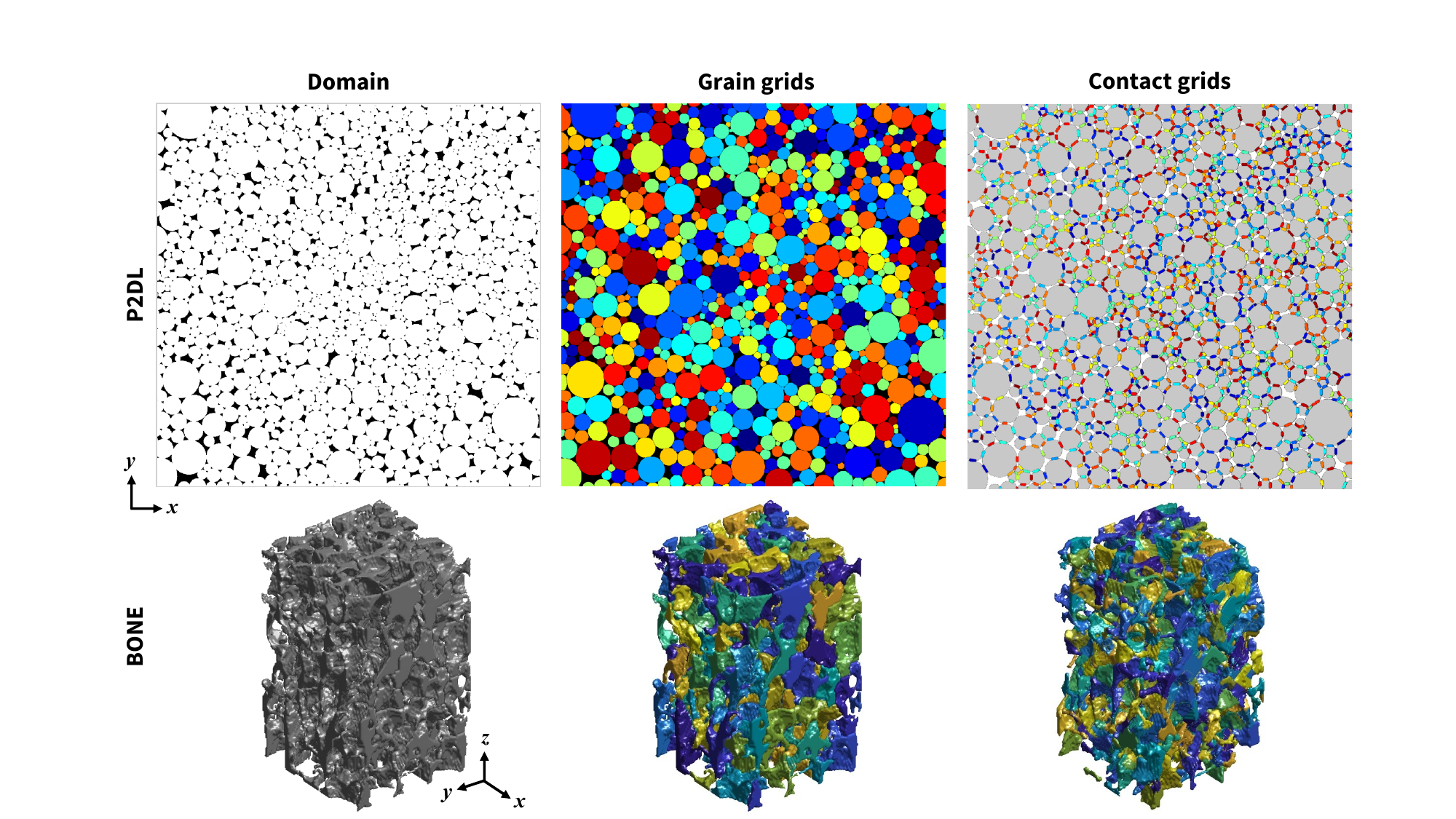}}
	\caption{Continued from Fig.\ref{fig:problem_set1}: (top row) A large disk pack (P2DL) and (bottom row) a 3D-rendered bone image (BONE) considered for testing the hPLMM preconditioner. The randomly colored regions (middle column) are grain grids, and the colored strips (right column) are contact grids.}
\label{fig:problem_set2}
\end{figure}

\begin{table}[h!]\centering
\caption{Geometric information and fine/coarse-grid properties of the domains shown in Figs.\ref{fig:problem_set1}-\ref{fig:problem_set2}.}
\begin{tabular} { P{2.6cm} P{2.6cm} P{2.5cm} P{2.3cm} P{2.1cm} P{2.1cm} P{2.1cm}}
\toprule
\multicolumn{1}{l}{Domains}  & Image size (pixels)   & Domain size (mm) & FEM elements   & FEM nodes   & Grain grids ($N^g$)  & Contact grids ($N^{\zeta}$)  \\
\midrule
\multicolumn{1}{l}{P2D} & $716\times 576$ & $7.16\times 5.76$ & 376,083 & 387,146 & 76 & 127   \\
\multicolumn{1}{l}{P2DL} & $2400\times 2400$ & $24\times 24$ & 5,357,211 & 5,425,144 & 904 & 1541   \\
\multicolumn{1}{l}{S2D, S2DH, S2DC} & $541\times 546$ & $5.41\times 5.46$ & 248,565 & 262,779 & 121 & 145 \\
\multicolumn{1}{l}{DARCY, Spectral} & $640\times 640$ & $0.8\times 0.8$ & 416,082 & 425,826 & 45 & 1  \\
\multicolumn{1}{l}{DARCY, Cartesian} & $640\times 640$ & $0.8\times 0.8$ & 417,120 & 426,213 & 49 & 1  \\
\multicolumn{1}{l}{BONE} & $100\times 100 \times 160$ & $1\times 1\times 1.6$ & 632,877 & 906,643 & 169 & 199  \\
\bottomrule
\end{tabular}
\label{tab:settings}
\end{table}

Table \ref{tab:settings} summarizes each domain's image size, physical dimensions, number of FEM elements and nodes, number of grain grids $N^g$, and number of contact grids $N^{\zeta}$ ($=\!N^c$). For each domain, we solve the linear system in Eq.\ref{eq:ls_all} using a right-preconditioned GMRES solver. We say the solver has converged if the normalized residual satisfies $\| \hat{A}\hat{x} - \hat{b} \| / \| \hat{b} \| \!<\! 10^{-9}$ or the number of iterations reaches 300. All simulations are run serially.  The machine specs for all domains except P2DL and BONE are Intel(R) Core(TM) i7-10700 CPU (8 cores, 2.90 GHz) with 32 GB of RAM, and for P2DL and BONE they are Intel(R) Xeon(R) Gold 6342 CPU (48 cores, 2.80 GHz) with 512 GB of RAM. 

To test the impact of smoother choice in hPLMM, we pair the coarse preconditioner, M\ts{G}, with two smoothers of the form given by Eq.\ref{eq:local_smoother}. They differ only by their base smoother M$_{l}$, for which we select: (1) M\ts{CG} in Section \ref{sec:local_smoother} with $n_{st}\!=\!1$ smoothing stages, and (2) M\ts{ILU(\textit{k})} with $n_{st}\!=\!6$, where the fill-level $k$ of incomplete LU-factorization is 0, unless stated otherwise. The chosen $n_{st}$ for each base smoother is optimal and was determined in \cite{li2024phase}. To construct M\ts{G}, we only use Gaussian mortars with $\beta\!=\!4$ in Eq.\ref{eq:gauss}, because Algebraic mortars yielded similar results as shown in \ref{app:alg_vs_gauss}. To benchmark hPLMM, we compare it against: (1) the low-order PLMM \cite{li2024phase}, obtained by setting the number of mortar nodes per interface, $n$, to one; (2) component-wise AMG (cAMG) proposed by \cite{gustafsson1998cAMGorig,white2019twostage}; and (3) Generalized Dryja-Smith-Widlund (GDSW) preconditioner by \cite{dohrmann2008family, heinlein2016parallel,heinlein2020fully}. The last two are detailed below: 
\\
\\
\textbf{cAMG}. Also known as the \textit{separate displacement component approximation}, the matrix $\hat{\mathrm{A}}$ is first permuted to obtain $\tilde{\mathrm{A}}$, where rows and columns associated with each coordinate direction are grouped together, as follows in 3D:
\begin{equation}
\tilde{\mathrm{A}} =
        \begin{bmatrix}
			\mathrm{A}_{xx} &\mathrm{A}_{xy} &\mathrm{A}_{xz} \\
			\mathrm{A}_{yx} &\mathrm{A}_{yy} &\mathrm{A}_{yz} \\
			\mathrm{A}_{zx} &\mathrm{A}_{zy} &\mathrm{A}_{zz} \\
		\end{bmatrix}
		\approx 
		\begin{bmatrix}
			\mathrm{A}_{xx} & & \\
			 &\mathrm{A}_{yy} & \\
			 & &\mathrm{A}_{zz} \\
		\end{bmatrix}		
\end{equation}
Then, the off-diagonal blocks are removed and the resulting matrix is used to build the cAMG preconditioner. The latter consists of building an AMG preconditioner for each diagonal block separately. Applying cAMG entails solving subsystems associated with these blocks independently. To build and apply the componentwise AMG preconditioners, we adapted the \texttt{amg.m} code from the iFEM GitHub repository \cite{chen2008ifem}, wherein a modified Ruge-Stuben \cite{stuben1987algebraic} coarsening and a two-point interpolation is used. The algorithm performs a single multigrid V-cycle, with the number of levels determined automatically (often 5-8), accompanied by one pre- and one post-smoothing per level via Gauss Seidel.
\\
\\
\textbf{GDSW}. This two-level Schwarz preconditioner bears similarities to the low-order PLMM \cite{li2024phase}. The key differences are: (1) In GDSW, the way $\Omega_s$ is decomposed is \textit{not} prescribed, and left arbitrary. Whereas PLMM insists on watershed segmentation as \textit{the} way to decompose a (pore-scale) structure because of the physical insight that stresses tend to localize at geometric choke points; (2) Shape functions in GDSW are built from harmonic extensions of rigid-body modes of each subdomain, restricted onto its interfaces, to the subdomain interior. This is algebraically equivalent to the closure BCs in PLMM to compute shape functions; (3) There are no contact grids in GDSW. Instead, a second set of subdomains is built by dilating the original non-overlapping ones to create overlaps with some prescribed thickness (here 16 FEM elements, which is comparable to the width of contact grids in hPLMM). These dilated subdomains are then used to build an additive-Schwarz smoother similar to Eq.\ref{eq:Mg_Mz}; (4) The $\mathrm{M_G}$ and $\mathrm{M_L}$ of GDSW are combined additively, not multiplicatively as is done in Eq.\ref{eq:multi_precond} for PLMM; and (5) In GDSW, the option to partition an interface into subregions is left open, in which case the method resembles hPLMM with piecewise constant mortar functions.

Since the arbitrariness of decomposition (Item 1) represents the most significant difference between PLMM and GDSW, we focus our testing on this aspect. We employ a Cartesian decomposition to construct $\mathrm{M_G}$ via the PLMM algorithm, then dilate the resulting rectangular grain grids to build $\mathrm{M_L}$ per Item 3. Fig.\ref{fig:GDSW} in the appendix illustrates this approach for the S2D domain. While we do not implement the additive (instead, multiplicative) combination from Item 4, we retain the ``GDSW'' designation for brevity. The number of grain grids in GDSW is comparable to those in Table \ref{tab:settings} for hPLMM. Specifically, 81, 900, 121, 49, and 175 for P2D, P2DL, S2D (S2DC and S2DH), DARCY, and BONE domains, respectively. The grain grids are roughly square (2D) or cubic (3D) in shape.

\section{Results} \label{sec:results}
We present our results in three parts. In Section \ref{sec:first_pass}, we measure the accuracy of the approximate, or \textit{first-pass}, solution obtained from a single application of M\ts{G} in Section \ref{sec:global_precon}, namely, $\hat{x}_{aprx}\!=\!\mathrm{M_G^{-1}}\hat{b}$. We show this approximation is usable in many applications as a standalone result. In Section \ref{sec:perf_precond}, we pair M\ts{G} with M\ts{L} and discuss the convergence rate of GMRES preconditioned by hPLMM. The associated computational costs are discussed next in Section \ref{sec:comp_cost}.

\subsection{Global preconditioner as an approximate solver} \label{sec:first_pass}
A high-quality first-pass solution $\hat{x}_{aprx}$ indicates the utility of M\ts{G} as a practical standalone approximation, and it generally correlates with faster GMRES convergence in Section \ref{sec:perf_precond}. In what follows, let $n$ denote the maximum number of mortar nodes placed on each contact interface of a domain. We say ``maximum'' because if $n$ is larger than the number of fine grids (FEM nodes) at an interface, then $n$ for that interface is set equal to its number of fine grids. We shall also use the term \textit{geometric hPLMM} to refer to the geometric formulation in \cite{khan2024high}, and \textit{algebraic hPLMM} to refer to the preconditioner formulation herein. Recall $n\!=\!1$ yields the low-order PLMM approximation in \cite{li2024phase}. 

To quantify $L_2$-errors associated with $\hat{x}_{aprx}$, we use the equation below following \cite{khan2024high}:
\begin{equation} \label{eq:err_metric}
	E^{\chi}_{2} = \left(\dfrac{1}{\vert\Omega_s\vert}
	     \int (E^{\chi}_{p})^2\, d\Omega\right)^{1/2}  
	     \,,\qquad 
	E^{\chi}_{p} = \dfrac{\Vert\chi_M-\chi_S\Vert}
	               {\textrm{sup}_{\Omega_s}\Vert \chi_S\Vert}
\end{equation}
where $\chi$ is a placeholder for either the fine-scale displacement, $\bs{u}$, or the fine-scale maximum shear stress, $\sigma^t$. The latter is defined as the radius of the largest Mohr circle, i.e., $(\sigma_{\max}-\sigma_{\min})/2$. The pointwise error $E^{\chi}_{p}$ is calculated by subtracting the approximate hPLMM solution, $\chi_M$, from the exact solution, $\chi_S$. This difference is then normalized by the supremum of $\chi_S$ over the solid domain $\Omega_s$. The $L_2$-error $E^{\chi}_{2}$ is the integral of the square of $E^{\chi}_{p}$ over $\Omega_s$. 
 
Table \ref{tab:err_hPLMM} summarizes the $L_2$-errors of both $\bs{u}$ and $\sigma^t$ on all domains under tensile and shear loading for different $n$. Included are also errors computed in \cite{khan2024high} for the geometric hPLMM on a subset of the domains. A comparison between the algebraic and geometric hPLMM shows that the two errors are of similar magnitude, indicating we have successfully algebraized the geometric hPLMM of \cite{khan2024high} into a coarse preconditioner M\ts{G}. The agreement is not exact because of $O(h)$ differences, where $h$ is the fine-grid size, in labeling the fine grids adjacent to the contact interfaces. As $n$ increases, we see that the errors of $\hat{x}_{aprx}$ also decrease, but more rapidly under shear than under tension. This is attributed to the fact that local bending/torsion moments are more dominant under shear than under tension.

\begin{figure} [t!]
  \centering
  \centerline{\includegraphics[scale=0.5,trim={55 192 55 45},clip]{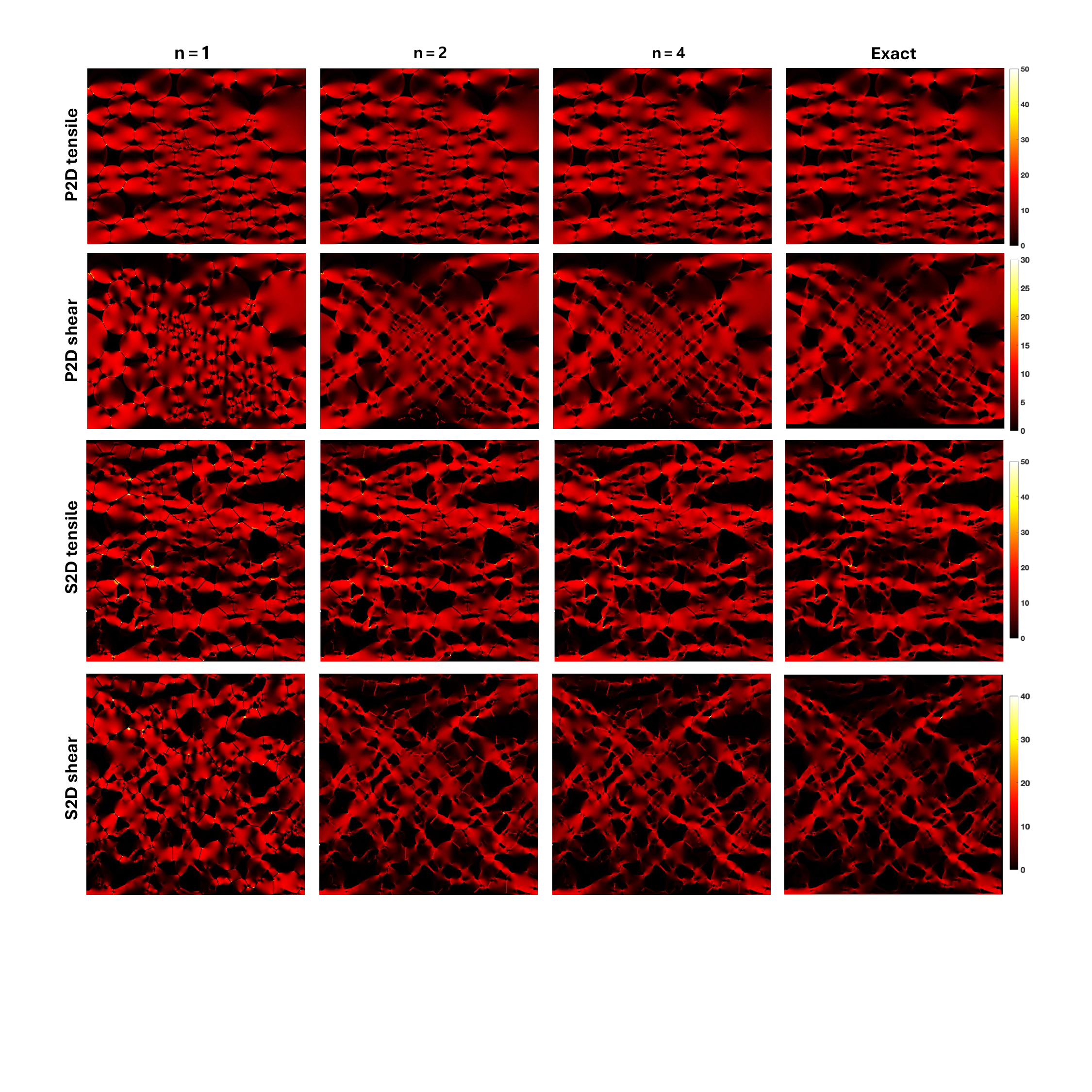}}
	\caption{Comparison of the spatial distributions of maximum shear stress ($\sigma^t$) obtained from a single application of the global preconditioner, M\ts{G}, in hPLMM versus the exact solution. Results correspond to the P2D and S2D domains with different number of mortar nodes per interface, $n$. Plots for both tensile and shear loading are shown, with the former constituting a less stringent test for hPLMM, i.e., good agreement at $n\!=\!1$.}
\label{fig:P2D_S2D_fp}
\end{figure}

\begin{figure} [t!]
  \centering
  \centerline{\includegraphics[scale=0.5,trim={240 50 220 20},clip]{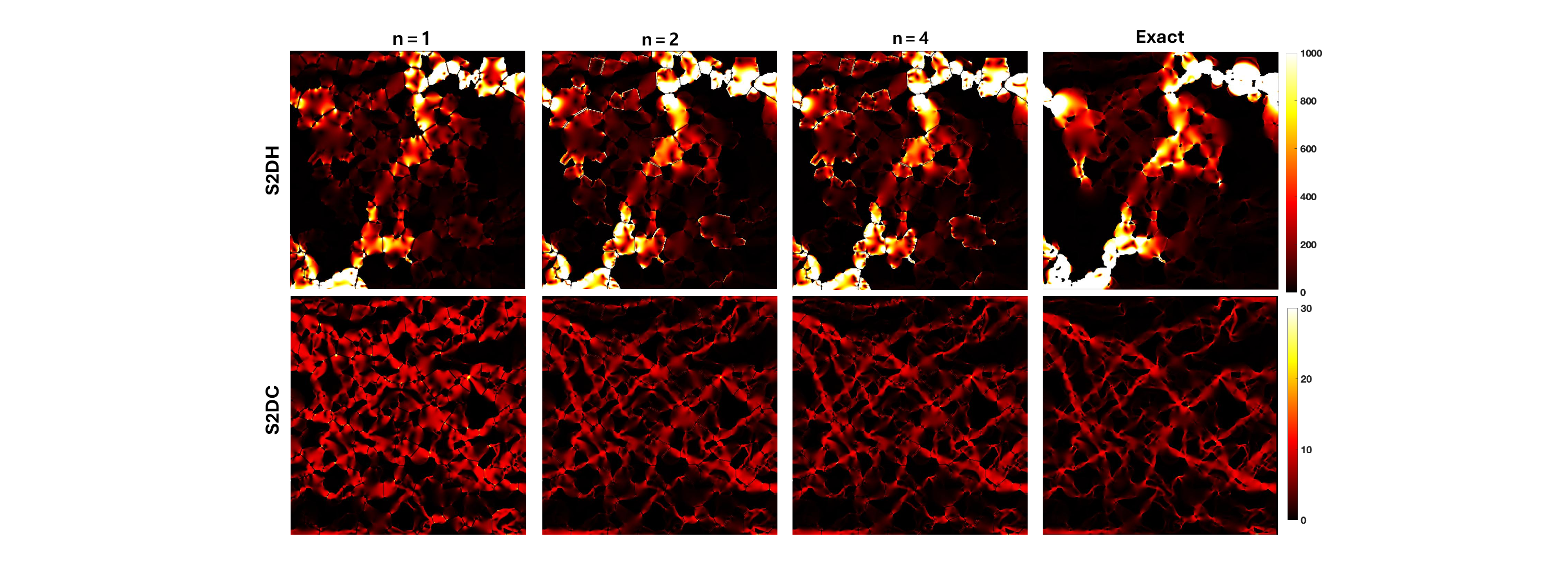}}
	\caption{Comparison of the spatial distributions of maximum shear stress ($\sigma^t$) obtained from a single application of the global preconditioner, M\ts{G}, in hPLMM versus the exact solution. Results correspond to the S2DH and S2DC domains with different number of mortar nodes per interface, $n$.}
\label{fig:S2DH_fp}
\end{figure}

\begin{figure} [t!]
  \centering
  \centerline{\includegraphics[scale=0.5,trim={0 50 7 0},clip]{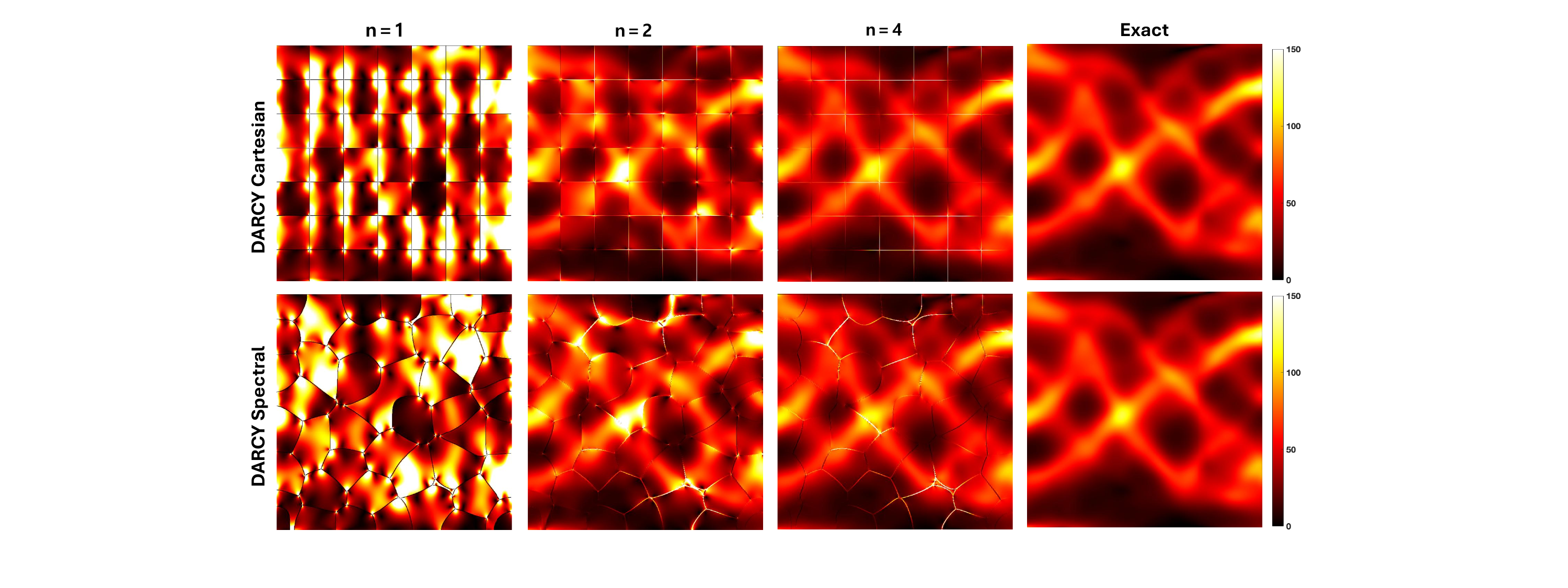}}
	\caption{Comparison of maximum shear stress ($\sigma^t$) from a single application of the global preconditioner, M\ts{G}, in hPLMM versus the exact solution. Results belong to the DARCY domain with Cartesian and spectral decomposition and different number of mortar nodes per interface, $n$.}
\label{fig:Darcy_fp}
\end{figure}

\begin{figure} [t!]
  \centering
  \centerline{\includegraphics[scale=0.5,trim={190 15 250 5},clip]{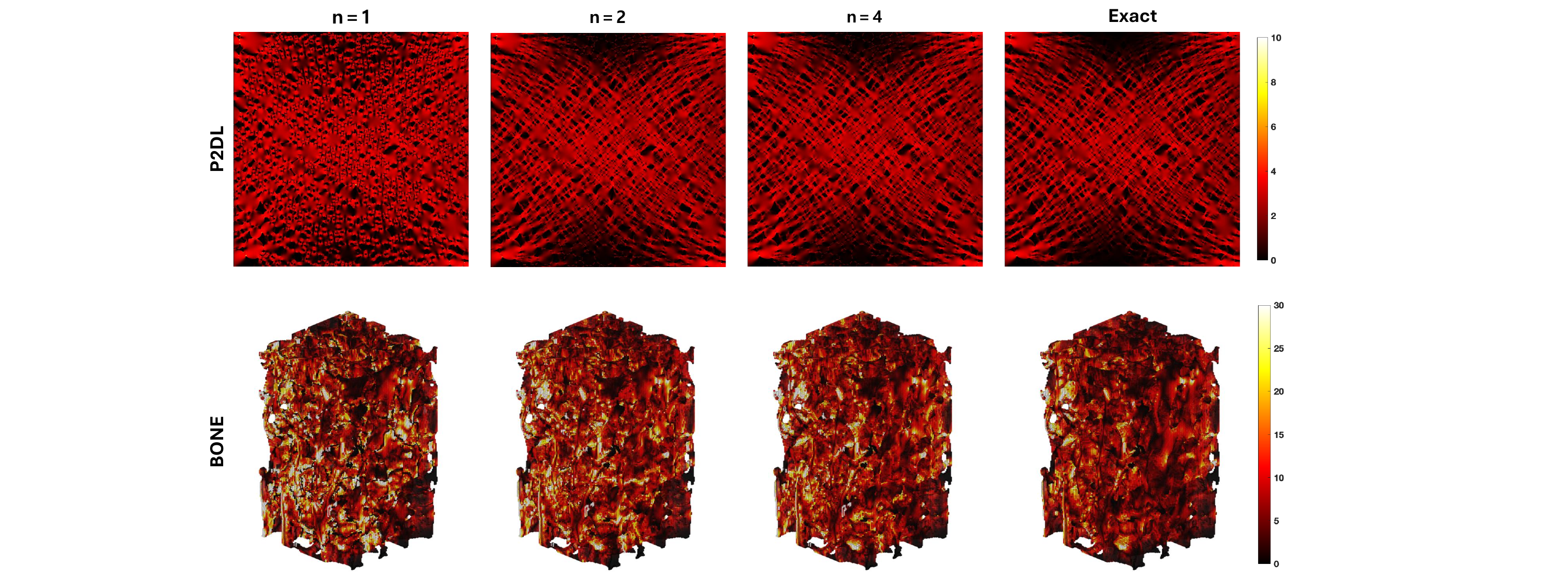}}
	\caption{Comparison of the spatial distributions of maximum shear stress ($\sigma^t$) obtained from a single application of the global preconditioner, M\ts{G}, in hPLMM versus the exact solution. Results correspond to the P2DL and BONE domains with different number of mortar nodes per interface, $n$.}
\label{fig:P2DL_BONE_fp}
\end{figure}

Figs.\ref{fig:P2D_S2D_fp}-\ref{fig:P2DL_BONE_fp} compare $\sigma^t$ from the algebraic hPLMM (i.e., $\hat{x}_{aprx}$) for different $n$ against the exact solution. In Fig.\ref{fig:P2D_S2D_fp} for P2D and S2D, we see that under tensile loading, PLMM ($n\!=\!1$) is already in good agreement with the exact solution, and increasing $n$ results in minor improvement. We therefore omit tensile data from all subsequent figures and focus only on the more stringent tests posed by shear loading. Similarly, we omit plots of displacement, as they are more forgiving in masking errors than $\sigma^t$. Incidentally, the omitted plots look identical to those in \cite{khan2024high} for the geometric hPLMM. Focusing on shear, Figs.\ref{fig:P2D_S2D_fp}-\ref{fig:P2DL_BONE_fp} show that $\hat{x}_{aprx}$ for $n\!=\!1$, namely PLMM, is rather error prone, but improves significantly as $n$ increases. Specifically, the overall force-chain pattern under shear is not captured well with $n\!=\!1$, because $\sigma^t$ appears choppy across contact interfaces. Errors are much larger for DARCY (Fig.\ref{fig:Darcy_fp}), with force chains for $n\!=\!1$ completely incorrect. The errors reduce substantially with $n\!=\!2$, and $\hat{x}_{aprx}$ becomes almost identical to the exact solution with $n\!=\!4$. The reason for the large errors of $n\!=\!1$ is that the localization assumption in Eq.\ref{eq:closure_hPLMM} reduces to a uniform displacement over each interface. Physically, this means interfaces are assumed to be rigid, thereby ignoring local bending/torsion moments. Under shear, such moments dominate microscale deformation \cite{khan2024crack}, rendering $\hat{x}_{aprx}$ a poor approximation. One takeaway here is that $n\!=\!1$ is often sufficient for tension, but $n\!=\!2$ to 4 is needed for shear.

\begin{table}[t!]
\centering
\small
\caption{Summary of $E^{\chi}_{2}$ errors ($\%$) in maximum shear stress ($\chi = \sigma^t$) and displacement ($\chi = \bs{u}$) corresponding to the first-pass solutions from hPLMM's coarse preconditioner M\ts{G} with different number of mortar nodes per interface, $n$. The $E^{\bs{u}}_{2}$ errors are given in parentheses.}
\label{tab:err_hPLMM}
\vspace{0.5em}
\hspace*{-0.2cm}
\begin{tabular}{
    c @{\hspace{1em}} | @{\hspace{1em}} 
    l
    @{\hspace{1.5em}} c @{\hspace{1.5em}} c @{\hspace{1.5em}} c
    @{\hspace{2em}} c @{\hspace{1.5em}} c @{\hspace{1.5em}} c
}
\toprule
\multicolumn{1}{c}{} & \multicolumn{1}{c}{} & \multicolumn{3}{c}{\textbf{Algebraic}} & \multicolumn{3}{c}{\textbf{Geometric}} \\
\cmidrule(lr){3-5} \cmidrule(lr){6-8}
\multicolumn{1}{c}{} & \textbf{Domain} & $n=1$ & $n=2$ & $n=4$ & $n=1$ & $n=2$ & $n=4$ \\
\midrule
\multirow{8}{*}{\rotatebox{90}{\textbf{Tension}}} 
 & P2D               & 3.51 (0.58)  & 2.22 (0.42)  & 2.00 (0.33)  & 2.61 (0.55)  & 1.20 (0.38)  & 1.06 (0.34) \\
 & S2D               & 3.90 (3.77)  & 2.11 (2.05)  & 1.64 (1.83)  & 2.03 (4.19)  & 0.81 (1.59)  & 0.70 (1.63) \\
 & S2DH              & 2.21 (1.23)  & 2.09 (1.09)  & 2.09 (1.12)  & -  & -  & - \\
 & S2DC              & 4.83 (5.14)  & 3.15 (4.63)  & 2.52 (3.90)   & -  & -  & - \\
 & DARCY, Spectral   & 25.5 (4.37)  & 6.60 (0.98)  & 1.96 (0.24)  & 28.2 (4.79)  & 5.27 (1.12)  & 1.69 (0.76) \\
 & DARCY, Cartesian  & 50.8 (5.90)  & 9.41 (1.20) & 0.97 (0.14)   & 41.5 (4.22)  & 8.25 (1.01)  & 0.84 (0.09) \\
 & P2DL              & 1.32 (0.37)  & 0.74 (0.26)  & 0.67 (0.24)  & -  & -  & - \\
 & BONE              & 9.84 (4.33)  & 5.55 (3.06)  & 4.02 (2.09)  & -  & -  & - \\
\midrule
\multirow{8}{*}{\rotatebox{90}{\textbf{Shear}}} 
 & P2D               & 5.44 (1.88)  & 1.70 (0.24)  & 1.58 (0.17)  & 5.49 (2.15)  & 1.11 (0.22)  & 0.98 (0.21) \\
 & S2D               & 5.44 (2.79)  & 1.93 (0.84)  & 1.56 (0.70)  & 6.47 (3.42)  & 1.25 (0.71)  & 1.10 (0.76) \\
 & S2DH              & 2.54 (4.19)  & 2.30 (4.10)  & 2.28 (4.14)  & -  & -  & - \\
 & S2DC              & 7.62 (4.39)  & 3.48 (2.23)  & 2.75 (1.77)  & -  & -  & - \\           
 & DARCY, Spectral   & 55.3 (10.4)  & 13.1 (1.41)  & 5.18 (0.39)  & 59.2 (10.2)  & 7.89 (1.16)  & 3.57 (0.68) \\
 & DARCY, Cartesian  & 97.3 (10.9) & 10.2 (1.05) & 1.41 (0.31)  & 81.4 (10.6)  & 9.18 (0.91)  & 1.10 (0.28) \\
 & P2DL              & 3.82 (7.25)  & 1.12 (4.63)  & 1.04 (3.90)   & -  & -  & - \\
 & BONE              & 7.95 (4.75)  & 5.31 (3.21)  & 3.79 (2.27)  & -  & -  & - \\    
\bottomrule
\end{tabular}
\end{table}

\subsection{Convergence rate of the preconditioned Krylov solver} \label{sec:perf_precond}

\begin{figure} [t!]
  \centering
  \hspace*{0.5cm}
  \centerline{\includegraphics[scale=0.6,trim={165 260 135 190},clip]{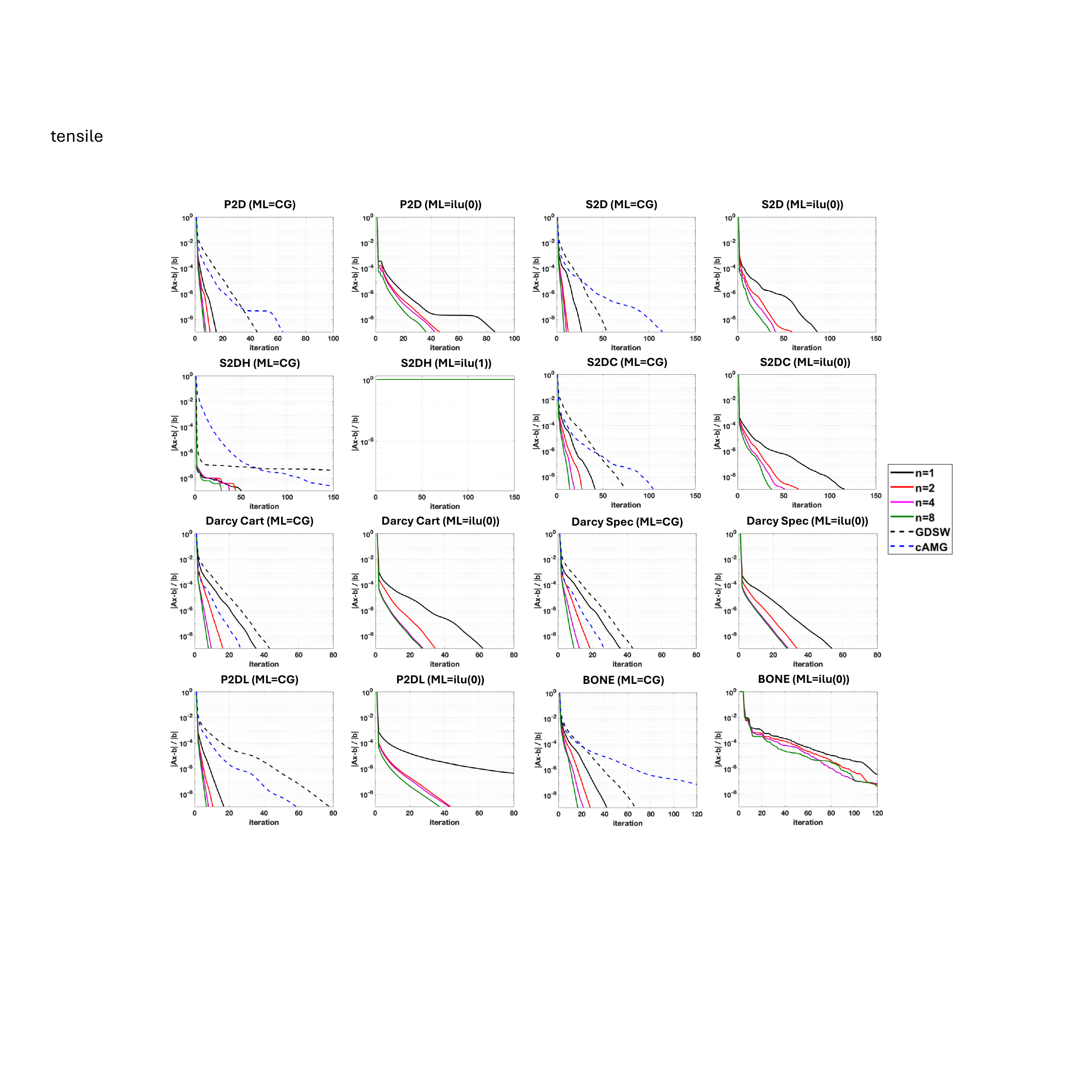}}
	\caption{Normalized residual versus number of GMRES iterations preconditioned by hPLMM, cAMG, and GDSW for all domains under \textit{tensile loading}. hPLMM results correspond to M\ts{G} combined with one of two smoothers, either M\ts{CG} with the number of stages in Eq.\ref{eq:local_smoother} equal to $n_{st}\!=\!1$ or M\ts{ILU($k$)} with $n_{st}\!=\!6$. M\ts{G} is built with different number of mortar nodes $n$ per interface. Recall $n\!=\!1$ corresponds to the low-order PLMM \cite{li2024phase}.}
\label{fig:iter_tense}
\end{figure}

\begin{figure} [t!]
  \centering
  \hspace*{0.5cm}
  \centerline{\includegraphics[scale=0.6,trim={165 260 135 190},clip]{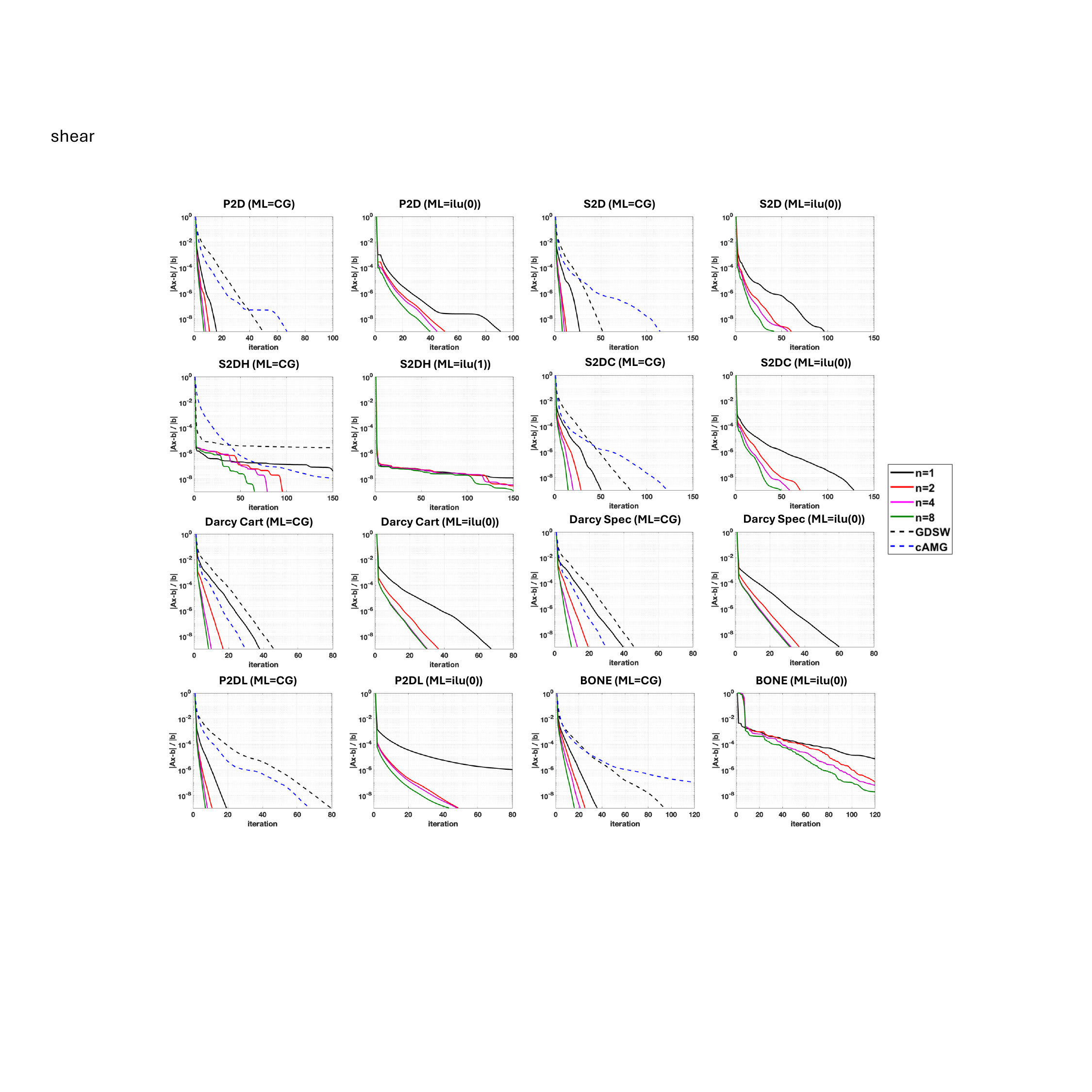}}
	\caption{Normalized residual versus number of GMRES iterations preconditioned by hPLMM, cAMG, and GDSW for all domains under \textit{shear loading}. hPLMM results correspond to M\ts{G} combined with one of two smoothers, either M\ts{CG} with the number of stages in Eq.\ref{eq:local_smoother} equal to $n_{st}\!=\!1$ or M\ts{ILU($k$)} with $n_{st}\!=\!6$. M\ts{G} is built with different number of mortar nodes $n$ per interface. Recall $n\!=\!1$ corresponds to the low-order PLMM \cite{li2024phase}.}
\label{fig:iter_shear}
\end{figure}

Here, we pair the M\ts{G} with a smoother M\ts{L} using Eq.\ref{eq:multi_precond}, then apply the combined hPLMM preconditioner within GMRES to solve Eq.\ref{eq:ls_all}. Table \ref{tab:wct} summarizes the number of GMRES iterations and associated wall-clock times (WCTs) required to converge ($\| \hat{A}\hat{x} - \hat{b} \| / \| \hat{b} \| \!<\! 10^{-9}$) in all domains under shear loading. Results for tensile loading are almost identical and are thus presented in \ref{app:Krylov_tensile}. Table \ref{tab:wct} includes the use of one of two smoothers, M\ts{CG} and M\ts{ILU($k$)}, with M\ts{G} in hPLMM and different number of mortar nodes per interface, namely, $n\!=\!1$, 2, 4, and 8. The ILU fill-level, $k$, is set to 0 in all domains except for S2DH, where it is set to 1. This is because GMRES failed to converge with $k\!=\!0$ for S2DH, which is otherwise the optimal choice according to previous work \cite{li2024phase}. The performances of the cAMG and GDSW preconditioners are also included in Table \ref{tab:wct} as benchmarks. Figs.\ref{fig:iter_tense}-\ref{fig:iter_shear} plot the corresponding GMRES convergence patterns in terms of normalized residual ($\| \hat{A}\hat{x} - \hat{b} \| / \| \hat{b} \|$) versus iterations. Fig.\ref{fig:iter_tense} is for tensile loading and Fig.\ref{fig:iter_shear} for shear, a comparison of which confirms that the convergence rates, thus WCTs, are comparable.

The following key observations stand out: (1) In all cases, the higher $n$ is, the faster GMRES converges, which is consistent with the improvements seen in the first-pass solutions of Figs.\ref{fig:P2D_S2D_fp}-\ref{fig:P2DL_BONE_fp} versus $n$; (2) The largest difference in convergence rate is seen between $n\!=\!1$ and $2$, suggesting that hPLMM is indeed superior to the low-order PLMM. However, gains in convergence rate diminish for larger $n$ ($>\!2$); (3) In all cases, the M\ts{CG} smoother paired with M\ts{G} performs significantly better than M\ts{ILU($k$)}. For the highly heterogeneous S2DH domain, M\ts{ILU($k$)} even struggles to converge regardless of $n$, whereas convergence rate is more than double with M\ts{CG} for $n\!>\!1$. Notice GMRES fails to converge with M\ts{ILU(0)} for S2DH under shear, which required the use of M\ts{ILU(1)} instead. But under tension, with all else held equal, M\ts{ILU(1)} could not converge as seen from Fig.\ref{fig:iter_tense}. The takeaway here is that black-box smoothers are not ideal pairs for the proposed M\ts{G}, whereas M\ts{CG} is; Lastly, (4) GMRES iterations are up to 5 times fewer for hPLMM than either cAMG or GDSW. In S2DH, neither GDSW nor cAMG converges under shear, highlighting the importance of preconditioners like hPLMM, which are informed by the structure and material heterogeneity of the domain.

\begin{figure} [t!]
  \centering
  \centerline{\includegraphics[scale=0.6,trim={90 55 170 15},clip]{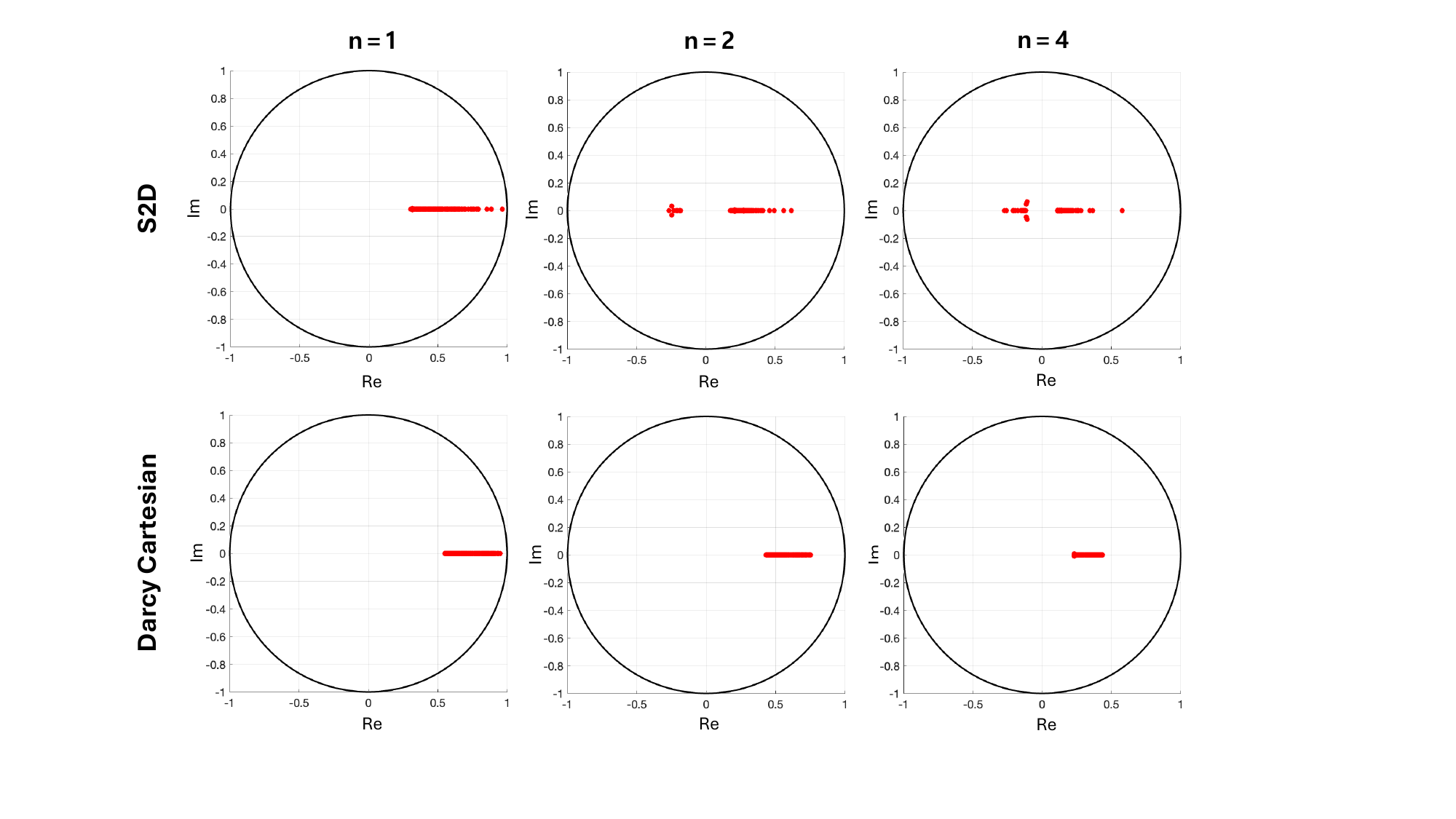}}
	\caption{Spectra (200 largest eigenvalues) of the error propagation matrix $\mathrm{E}$ in Eq.\ref{eq:error_prop} associated with the hPLMM preconditioner for different number of mortar nodes $n$ per interface. Results correspond to the S2D and DARCY (with Cartesian decomposition) domains under shear loading.}
\label{fig:spectra}
\end{figure}

A useful metric for understanding the rapid convergence of the hPLMM preconditioner is the spectral radius of its error propagation matrix $\mathrm{I}-\mathrm{M}^{-1}\hat{\mathrm{A}}$, which can be split into the actions of M\ts{G} and M\ts{L} as follows:
\begin{equation} \label{eq:error_prop}
\mathrm{E} = (\mathrm{I}-\mathrm{M}^{-1}_{\mathrm{L}} \hat{\mathrm{A}})(\mathrm{I}-\mathrm{M}^{-1}_{\mathrm{G}} \hat{\mathrm{A}})
\end{equation}
The split is a consequence of M being a multiplicative combination of M\ts{G} and M\ts{L} in Eq.\ref{eq:multi_precond}. A smaller than unity spectral radius for E implies that a basic iterative solver of the form $\hat{x}_{k+1}\!=\!\hat{x}_{k} + \mathrm{M}^{-1}(\hat{b} - \mathrm{\hat{A}}\hat{x}_k)$ converges. The smaller the spectral radius, the faster the convergence. While GMRES iterations are more involved, the same logic applies. Fig.\ref{fig:spectra} plots the 200 largest eigenvalues of E on the complex plane with $n\!=\!1$, 2, and 4 for the S2D and DARCY (with Cartesian decomposition) domains under shear. Here, M\ts{G} is combined with $\mathrm{M_L}\!=\!\mathrm{M_{CG}}$. We see that as $n$ increases, the eigenvalues of E become more clustered near the origin, resulting in smaller spectral radius. For S2D, the clustering is more dramatic from $n\!=\!1$ to $2$, than from $n\!=\!2$ to 4, which echoes earlier observations of first-pass solutions in Fig.\ref{fig:P2D_S2D_fp}, with minor improvement seen for $n\!>\!2$. For DARCY, clustering plateaus only after $n\!>\!4$, consistent with Fig.\ref{fig:Darcy_fp}.

\begin{table}[t!]
\centering
\caption{Summary of the number of iterations and wall-clock times (WCTs) in seconds required by GMRES to converge ($\| \hat{\mathrm{A}} \hat{x} - \hat{b} \| / \| \hat{b} \| < 10^{-9}$) to the solution of Eq.\ref{eq:ls_all} using different preconditioners. This includes hPLMM, where M\ts{G} is combined with either the M\ts{CG} ($n_{st}\!=\!1$) or M\ts{ILU($k$)} ($n_{st}\!=\!6$) smoother, and the number of mortar nodes per interface is set to $n\!=\!1$, 2, 4, and 8. The $n\!=\!1$ case corresponds to the low-order PLMM. Results for GDSW and cAMG are included as benchmarks. The total WCT, T\ts{tot}, consists of the time needed to build the smoother, T\ts{ML}, time to build the coarse preconditioner, T\ts{MG}, and the time spent by GMRES itself, T\ts{sol}. In all domains, the ILU fill-level is $k\!=\!0$ except in S2DH, where it is $k\!=\!1$. All cases correspond to \textit{shear loading} and \textcolor{red}{red} means ``diverged'' (i.e., not converged within 300 iterations).}
\footnotesize
\setlength{\tabcolsep}{3.5pt}
\begin{tabular}{l l *{4}{r} @{\hspace{12pt}} *{4}{r} @{\hspace{12pt}} r r}
\toprule
  & & \multicolumn{4}{c}{\textbf{Contact-Grain (CG)}} 
  & \multicolumn{4}{c}{\textbf{ILU(0) / ILU(1)}} 
  &  &  \\
\cmidrule(lr){3-6} \cmidrule(lr){7-10} 
 \textbf{Domain} & \textbf{Metric} & $n=1$ & $n=2$ & $n=4$ & $n=8$ & $n=1$ & $n=2$ & $n=4$ & $n=8$ & \textbf{GDSW} & \textbf{cAMG} \\
\midrule

\multirow{5}{*}{P2D} 
 & T\textsubscript{ML}  & 7.60 & 8.08 & 7.93 & 7.39 & 0.11 & 0.14 & 0.12 & 0.11 & 10.92 & -- \\
 & T\textsubscript{MG}  & 10.84 & 14.68 & 22.62 & 39.75 & 11.13 & 14.43 & 22.47 & 40.87 & 10.30 & 1.75 \\
 & T\textsubscript{sol} & 6.17 & 4.25 & 3.69 & 4.82 & 46.44 & 22.37 & 22.39 & 27.94 & 28.08 & 170.9 \\
 & \textbf{T\textsubscript{tot}} & \textbf{24.61} & \textbf{27.01} & \textbf{34.24} & \textbf{51.96} & \textbf{57.68} & \textbf{36.94} & \textbf{44.98} & \textbf{68.92} & \textbf{49.3} & \textbf{172.6} \\
 & \textbf{Iter.}  & \textbf{17} & \textbf{12} & \textbf{9} & \textbf{8} & \textbf{91} & \textbf{51} & \textbf{45} & \textbf{40} & \textbf{50} & \textbf{77} \\
\midrule

\multirow{5}{*}{S2D} 
 & T\textsubscript{ML}  & 4.15 & 4.37 & 4.16 & 4.23 & 0.07 & 0.09 & 0.07 & 0.09 & 6.45 & -- \\
 & T\textsubscript{MG}  & 5.52 & 8.23 & 11.86 & 22.67 & 5.45 & 7.60 & 11.47 & 22.49 & 5.26 & 1.56 \\
 & T\textsubscript{sol} & 6.07 & 2.88 & 2.86 & 5.23 & 31.83 & 18.42 & 19.88 & 26.76 & 18.00 & 205.3 \\
 & \textbf{T\textsubscript{tot}} & \textbf{15.74} & \textbf{15.48} & \textbf{18.88} & \textbf{32.13} & \textbf{37.35} & \textbf{26.11} & \textbf{31.42} & \textbf{49.34} & \textbf{29.71} & \textbf{206.9} \\
 & \textbf{Iter.}  & \textbf{28} & \textbf{13} & \textbf{11} & \textbf{8} & \textbf{97} & \textbf{61} & \textbf{57} & \textbf{42} & \textbf{53} & \textbf{138} \\
\midrule

\multirow{5}{*}{S2DH} 
 & T\textsubscript{ML}  & 4.10 & 4.19 & 4.16 & 4.09 & 0.73 & 1.09 & 0.72 & 0.74 & \textcolor{red}{7.39} & \textcolor{red}{--} \\
 & T\textsubscript{MG}  & 5.40 & 7.37 & 11.18 & 21.55 & 5.52 & 7.28 & 11.18 & 23.39 & \textcolor{red}{5.60} & \textcolor{red}{1.38} \\
 & T\textsubscript{sol} & 106.8 & 31.96 & 27.66 & 43.50 & 202.0 & 93.36 & 106.2 & 135.8 & \textcolor{red}{218.0} & \textcolor{red}{627.2} \\
 & \textbf{T\textsubscript{tot}} & \textbf{116.3} & \textbf{43.52} & \textbf{43.00} & \textbf{69.14} & \textbf{208.2} & \textbf{101.7} & \textbf{118.1} & \textbf{159.9} & \textcolor{red}{\textbf{231.0}} & \textcolor{red}{\textbf{628.6}} \\
 & \textbf{Iter.}  & \textbf{219} & \textbf{96} & \textbf{80} & \textbf{66} & \textbf{296} & \textbf{174} & \textbf{181} & \textbf{151} & \textcolor{red}{\textbf{300}} & \textcolor{red}{\textbf{300}} \\
\midrule

\multirow{5}{*}{S2DC} 
 & T\textsubscript{ML}  & 3.52 & 3.52 & 3.66 & 3.67 & 0.11 & 0.11 & 0.11 & 0.12 & 4.59 & -- \\
 & T\textsubscript{MG}  & 4.99 & 6.84 & 11.63 & 21.42 & 4.99 & 6.87 & 11.23 & 21.70 & 5.39 & 1.33 \\
 & T\textsubscript{sol} & 12.73 & 6.91 & 5.94 & 9.13 & 53.04 & 25.27 & 23.33 & 35.07 & 30.87 & 232.2 \\
 & \textbf{T\textsubscript{tot}} & \textbf{21.24} & \textbf{17.27} & \textbf{21.23} & \textbf{34.22} & \textbf{58.14} & \textbf{32.25} & \textbf{34.67} & \textbf{56.89} & \textbf{40.85} & \textbf{233.5} \\
 & \textbf{Iter.}  & \textbf{51} & \textbf{29} & \textbf{21} & \textbf{15} & \textbf{129} & \textbf{71} & \textbf{59} & \textbf{50} & \textbf{82} & \textbf{123} \\
\midrule

\multirow{5}{*}{DARCY Cart} 
 & T\textsubscript{ML}  & 6.99 & 6.79 & 6.76 & 6.93 & 0.19 & 0.23 & 0.19 & 0.20 & 9.94 & -- \\
 & T\textsubscript{MG}  & 9.79 & 13.17 & 19.46 & 37.07 & 9.69 & 13.20 & 20.29 & 36.54 & 9.96 & 1.98 \\
 & T\textsubscript{sol} & 16.84 & 6.83 & 4.19 & 4.43 & 41.34 & 19.90 & 16.80 & 18.36 & 24.19 & 85.37 \\
 & \textbf{T\textsubscript{tot}} & \textbf{33.62} & \textbf{26.79} & \textbf{30.41} & \textbf{48.43} & \textbf{51.22} & \textbf{33.33} & \textbf{37.28} & \textbf{55.10} & \textbf{44.09} & \textbf{87.35} \\
 & \textbf{Iter.}  & \textbf{38} & \textbf{17} & \textbf{10} & \textbf{9} & \textbf{68} & \textbf{37} & \textbf{31} & \textbf{30} & \textbf{46} & \textbf{30} \\
\midrule

\multirow{5}{*}{DARCY Spec} 
 & T\textsubscript{ML}  & 7.83 & 7.83 & 7.40 & 7.22 & 0.21 & 0.20 & 0.20 & 0.19 & 10.17 & -- \\
 & T\textsubscript{MG}  & 11.41 & 15.73 & 25.35 & 49.19 & 11.24 & 15.51 & 25.59 & 50.51 & 10.12 & 1.97 \\
 & T\textsubscript{sol} & 18.13 & 8.28 & 6.35 & 5.39 & 33.59 & 20.52 & 18.87 & 21.43 & 24.25 & 86.48 \\
 & \textbf{T\textsubscript{tot}} & \textbf{37.37} & \textbf{31.84} & \textbf{39.10} & \textbf{61.80} & \textbf{45.04} & \textbf{36.23} & \textbf{44.66} & \textbf{72.13} & \textbf{44.54} & \textbf{88.45} \\
 & \textbf{Iter.}  & \textbf{40} & \textbf{20} & \textbf{14} & \textbf{10} & \textbf{60} & \textbf{37} & \textbf{33} & \textbf{32} & \textbf{46} & \textbf{30} \\
\midrule

\multirow{5}{*}{P2DL} 
 & T\textsubscript{ML}  & 230.3 & 175.4 & 229.5 & 215.3 & \textcolor{red}{1.90} & 1.84 & 2.12 & 2.35 & 143.4 & -- \\
 & T\textsubscript{MG}  & 216.0 & 285.6 & 478.0 & 963.6 & \textcolor{red}{217.1} & 290.7 & 482.5 & 961.6 & 181.9 & 36.25 \\
 & T\textsubscript{sol} & 147.4 & 85.64 & 81.48 & 142.8 & \textcolor{red}{5106.2} & 448.6 & 543.2 & 851.2 & 612.0 & 3083.1 \\
 & \textbf{T\textsubscript{tot}} & \textbf{593.7} & \textbf{546.6} & \textbf{789.0} & \textbf{1321.6} & \textcolor{red}{\textbf{5325.2}} & \textbf{741.1} & \textbf{1027.8} & \textbf{1815.1} & \textbf{937.3} & \textbf{3119.3} \\
 & \textbf{Iter.}  & \textbf{20} & \textbf{11} & \textbf{9} & \textbf{8} & \textcolor{red}{\textbf{300}} & \textbf{49} & \textbf{49} & \textbf{44} & \textbf{80} & \textbf{68} \\
\midrule

\multirow{5}{*}{BONE}
 & T\textsubscript{ML}  & 126.1 & 126.0 & 142.6 & 140.3 & \textcolor{red}{5.58} & 5.60 & 6.17 & 6.15 & 109.2 & -- \\
 & T\textsubscript{MG}  & 183.1 & 335.2 & 757.9 & 2038.3 & \textcolor{red}{184.0} & 339.7 & 745.4 & 1892.2 & 117.1 & 11.18 \\
 & T\textsubscript{sol} & 147.2 & 124.6 & 463.8 & 2123.7 & \textcolor{red}{1437.9} & 1612.3 & 4987.6 & 23775.0 & 723.4 & 3310.0 \\
 & \textbf{T\textsubscript{tot}} & \textbf{456.4} & \textbf{585.9} & \textbf{1364.3} & \textbf{4302.3} & \textbf{\textcolor{red}{1627.5}} & \textbf{1957.6} & \textbf{5739.1} & \textbf{25673.3} & \textbf{949.7} & \textbf{3321.2} \\
 & \textbf{Iter.}  & \textbf{36} & \textbf{26} & \textbf{22} & \textbf{17} & \textcolor{red}{\textbf{300}} & \textbf{289} & \textbf{271} & \textbf{243} & \textbf{56} & \textbf{233} \\

\bottomrule
\end{tabular}
\label{tab:wct}
\end{table}


\subsection{Computational cost of building and using the preconditioner} \label{sec:comp_cost}
\begin{figure} [t!]
  \centering
  \centerline{\includegraphics[scale=0.7,trim={30 135 30 60},clip]{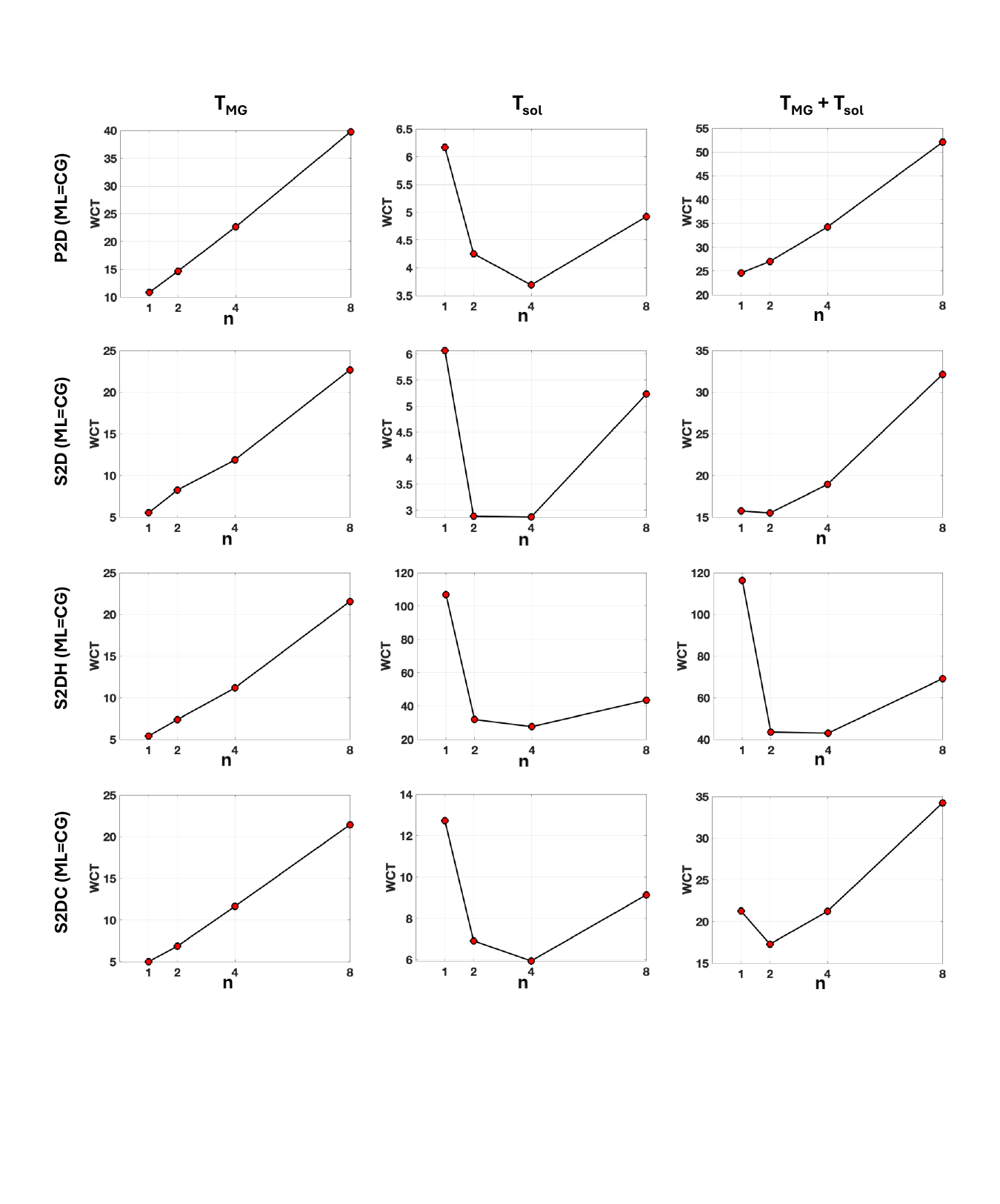}}
	\caption{Wall clock times (WCTs) spent solving Eq.\ref{eq:ls_all} via GMRES preconditioned by hPLMM versus the number of mortar nodes per interface $n$ in the P2D, S2D, S2DH, and S2DC domains under shear loading. Depicted are the WCTs of building the coarse preconditioner, T\ts{MG}, and the time spent by GMRES itself, T\ts{sol}. The last column is the sum of T\ts{MG} and T\ts{sol}, which excludes the smoother's build-time (= T\ts{tot} - T\ts{ML} in Table \ref{tab:wct}).}
\label{fig:wct_1}
\end{figure}

\begin{figure} [t!]
  \centering
  \centerline{\includegraphics[scale=0.7,trim={30 105 30 90},clip]{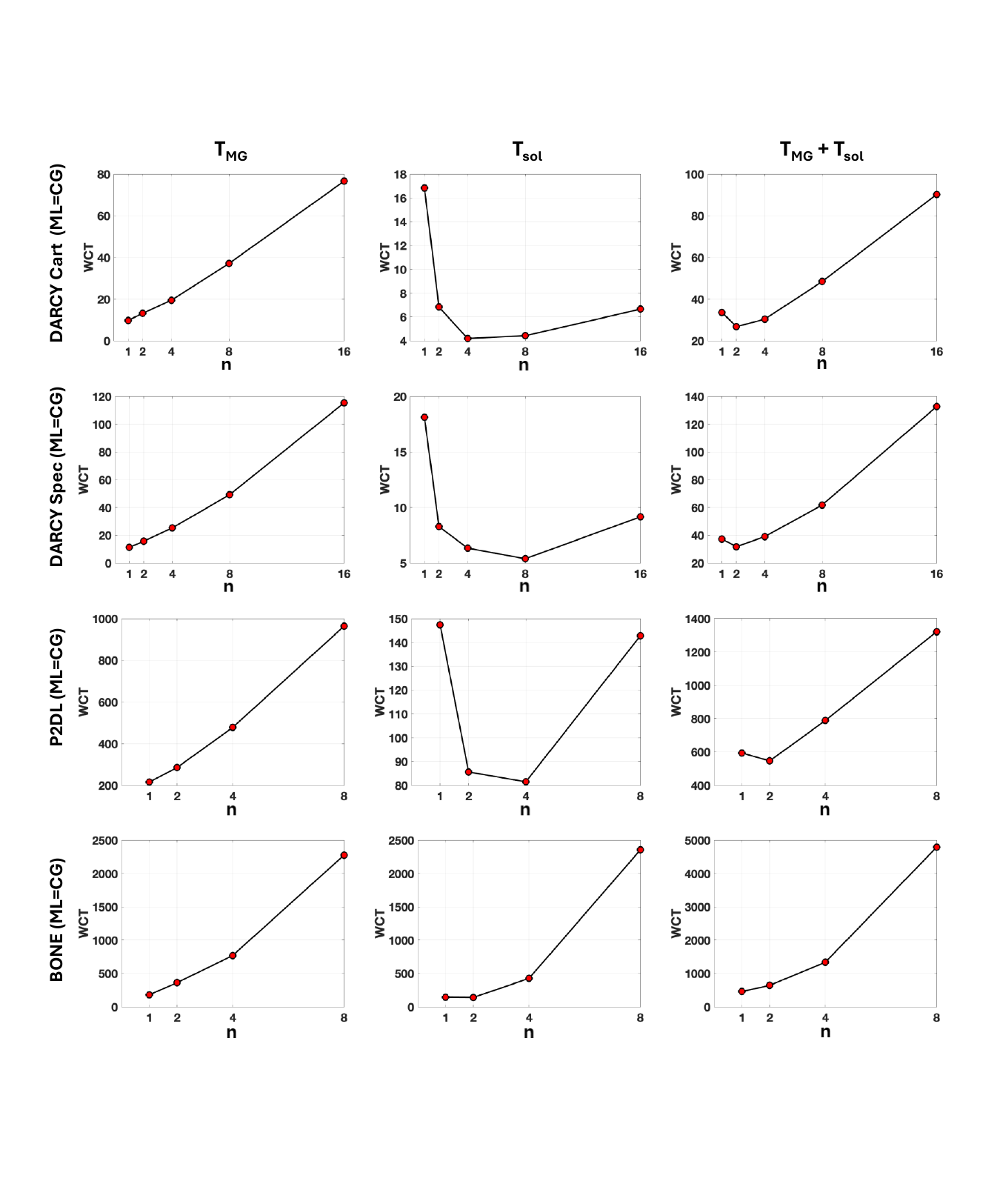}}
	\caption{Wall clock times (WCTs) spent solving Eq.\ref{eq:ls_all} via GMRES preconditioned by hPLMM versus the number of mortar nodes per interface $n$ in the  DARCY Cartesian, DARCY Spectral, P2DL, and BONE domains under shear loading. Depicted are the WCTs of building the coarse preconditioner, T\ts{MG}, and the time spent by GMRES itself, T\ts{sol}. The last column is the sum of T\ts{MG} and T\ts{sol} (= T\ts{tot} - T\ts{ML} in Table \ref{tab:wct}).}
\label{fig:wct_2}
\end{figure}

Table \ref{tab:wct} includes wall-clock times (WCTs) associated with solving Eq.\ref{eq:ls_all} via GMRES for all domains under shear loading. The total WCT, denoted by T\ts{tot} and highlighted in bold font, is comprised of the time required to build the coarse preconditioner, T\ts{MG}, the time needed to build the smoother, T\ts{ML}, and the time spent by GMRES itself, T\ts{sol}. The results include hPLMM with two different smoothers (M\ts{CG} and M\ts{ILU($k$)}) and different number of mortar nodes per interface ($n$), cAMG, and GDSW. In hPLMM and GDSW, ``building the smoother'' refers to performing LU-decompositions of the local systems defined over each grain grid to speed up their repeated solves during GMRES iterations. For cAMG, no such cost is listed separately, as T\ts{MG} includes the overall cost of building the preconditioner. Figs.\ref{fig:wct_1}-\ref{fig:wct_2} depict the WCTs for hPLMM versus $n$ for all domains under shear, when the M\ts{CG} smoother is used.

\begin{figure} [t!]
  \centering
  \centerline{\includegraphics[scale=0.7,trim={120 118 115 105},clip]{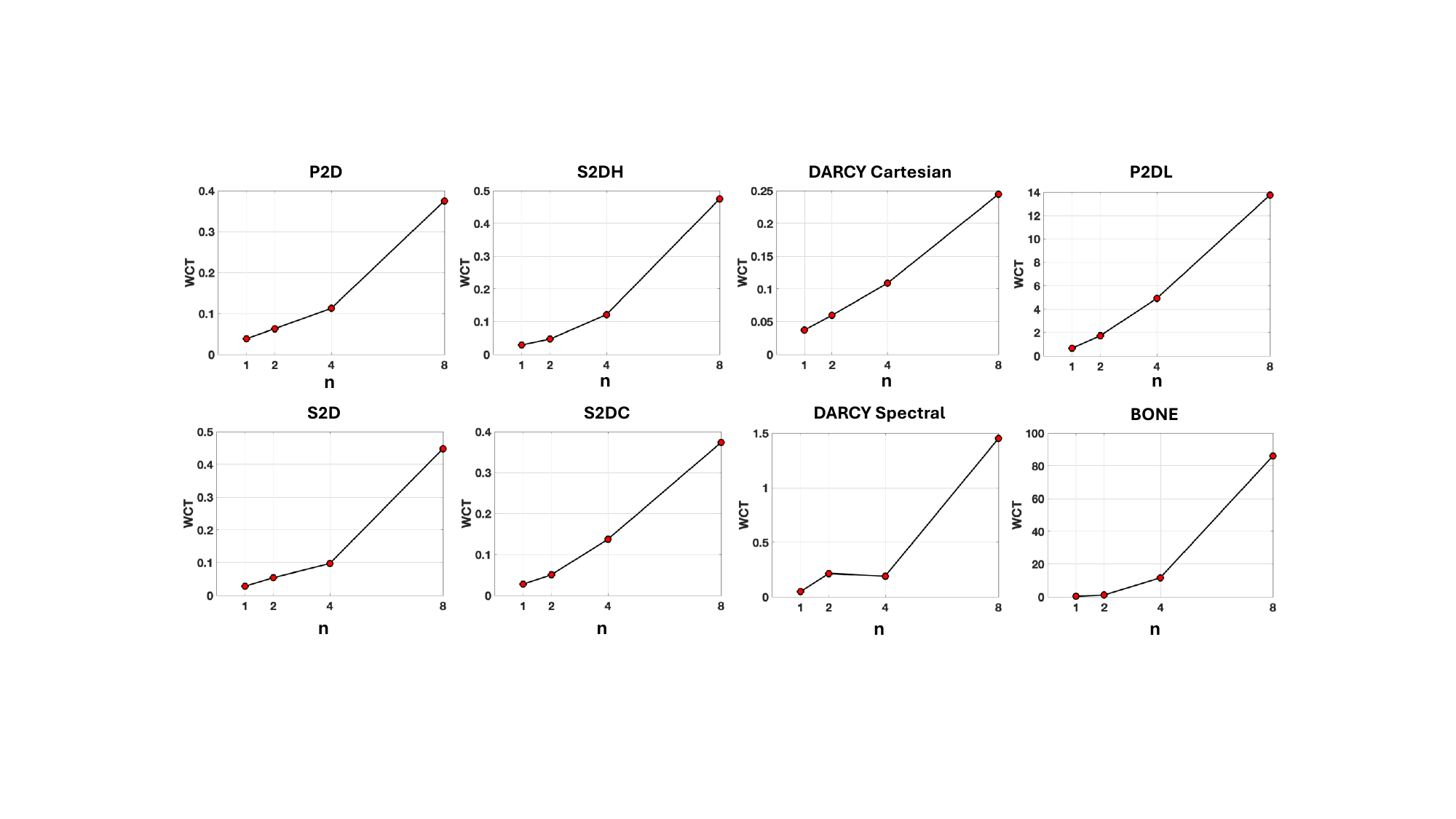}}
	\caption{Wall clock times (WCT) associated with a single application of the coarse preconditioner M\ts{G} of hPLMM (i.e., $\mathrm{M_G^{-1}}\hat{b}$) for different number of mortar nodes per interface $n$ in all of the domains considered in Figs.\ref{fig:problem_set1}-\ref{fig:problem_set2} under shear loading.}
\label{fig:wct_fp}
\end{figure}

We make three key observations. First, in all cases, hPLMM with the M\ts{CG} smoother exhibits the best performance over the benchmarks GDSW and cAMG, as well as the hPLMM preconditioner with the M\ts{ILU($k$)} smoother. The latter indicates the superiority of M\ts{CG} as a compatible smoother for hPLMM. Second, in hPLMM, as $n$ increases, T\ts{MG} also increases, whereas T\ts{sol} decreases up to $n\!=\!4$ then starts to increase for $n\!>\!4$. The reason T\ts{MG} increases is because a larger $n$ entails computing more columns for the reduced prolongation matrix, P, in Eq.\ref{eq:bc_matrix}. Each added mortar node per interface requires $D$ additional shape vectors to be built, via Eq.\ref{eq:col_pk}, on each of the grain grids sharing that interface. The reason T\ts{sol} behaves non-monotonically, as seen in Figs.\ref{fig:wct_1}-\ref{fig:wct_2}, is because as $n$ grows up to $\sim$4, the accuracy of M\ts{G} as a coarse approximator increases, resulting in fewer GMRES iterations (Figs.\ref{fig:iter_tense}-\ref{fig:iter_shear}). But for $n\!>\!4$, the cost of \textit{applying} M\ts{G} within GMRES starts to dominate, i.e., $\mathrm{M_G^{-1}}v_k$ for some $v_k$. This is because the size of the coarse matrix $\!\mathrm{\hat{R}\hat{A}\hat{P}}$ in Eq.\ref{eq:mg_struct} is proportional to $n\!\times\!D$, for which we use a direct solver to apply its inverse to a given RHS vector. Thus, when $n\!>\!4$, the number of GMRES iterations required to converge is low, but each iteration costs more. This is confirmed by Fig.\ref{fig:wct_fp}, where the cost of a single application of M\ts{G} (i.e.,  $\mathrm{M_G^{-1}}\hat{b}$) is seen to grow with $n$.

Our third and final observation is that the total cost of hPLMM, T\ts{tot}, is roughly constant in going from $n\!=\!1$ to 4, except in S2DH and S2DC where a significant drop in cost is observed. The fact that T\ts{tot} for the low-order PLMM ($n\!=\!1$) and the high-order hPLMM ($n\!=\!2-4$) are comparable may give the impression that the added cost of building M\ts{G} in hPLMM is not worthwhile. This is true if the preconditioner is meant to be used to solve only \textit{one} system. However, in numerous problems such as plasticity, wave propagation, and fracture mechanics, linear or linearized systems like Eq.\ref{eq:ls_all} must be solved repeatedly over multiple time steps, load steps, and/or Newton or staggered iterations. In most such cases, the coefficient matrix $\mathrm{\hat{A}}$ remains unchanged (load steps) or is perturbed only slightly (local cracks), in which case the M\ts{G} of hPLMM can be reused. The reusability of M\ts{G} in cracked domains, despite being originally built from an intact domain, has already been successfully demonstrated in PLMM \cite{mehmani2023precond, li2024phase}. Therefore, for solving multiple linear systems with different RHS vectors, T\ts{sol} is more important than T\ts{tot}, as the former dominates T\ts{MG}. In all our domains, except BONE, the drop in T\ts{sol} is at least a factor of two in going from $n\!=\!1$ to 4. For a 100-load-step simulation, the same factor would differentiate the total WCT of PLMM versus hPLMM with $n\!=\!4$. In Section \ref{sec:discussion}, we discuss how multi-level formulations of M\ts{G} can drive down T\ts{sol} even when $n\!>\!4$. 


\section{Discussion}\label{sec:discussion}
\subsection{Need for multi-level formulation}\label{sec:multilevel}
The results of Section \ref{sec:results} indicate that hPLMM is a significant improvement over the low-order PLMM, which itself is superior to existing preconditioners, such as cAMG and GDSW. The latter highlights the importance of embedding information about the physics and geometry of the problem into preconditioners. As an approximate solver, through $\hat{x}_{aprx}\!=\!\mathrm{M_G^{-1}}\hat{b}$, hPLMM is a successful algebraization of the geometric multiscale method of \cite{khan2024high} and yields more accurate first-pass solutions than PLMM under loading conditions that induce substantial bending/torsion moments locally. This is the case for the shear but not tensile loads in Section \ref{sec:results}. The reason for PLMM's lower accuracy is that contact interfaces are assumed to be rigid, a consequence of assigning a single coarse-scale displacement to them. As a preconditioner too, hPLMM is superior to PLMM, cAMG, and GDSW. While the total cost of hPLMM ($n\!=\!2-4$) is comparable to that of PLMM ($n\!=\!1$), this cost is broken down differently into the WCTs of building M\ts{G}, T\ts{MG}, and the self-time of GMRES, T\ts{sol}. In PLMM, a large portion of the cost comes from T\ts{sol}, whereas in hPLMM, T\ts{sol} is lower but the difference is compensated by a higher T\ts{MG}. Hence, if the goal is to solve only one linear system like Eq.\ref{eq:ls_all}, little gains are made by using hPLMM over PLMM. But if the goal is to solve a multitude of systems with differing RHS vectors, or even slightly perturbed coefficient matrices \cite{li2024phase}, then hPLMM would be at least twice as fast.

For $n\!>\!4$ in hPLMM, T\ts{sol} increases due to the larger coarse system involving the matrix $\mathrm{A}^o\!=\!\mathrm{\hat{R}\hat{A}\hat{P}}$ in Eq.\ref{eq:mg_struct} that must be solved with every iteration of GMRES. Simply put, $\mathrm{A}^o$ becomes too large and the only way to reduce the cost of applying its inverse to a vector would be to extend hPLMM from a two-level algorithm to a multi-level one. The procedure would involve coarsening $\mathrm{A}^o$ further, in hierarchical fashion, by grouping collections of nearby grain grids into ``macro-grain grids.'' For large $n$, this is absolutely necessary and should be the subject of future research.

\subsection{Mortar nodes and functions}\label{sec:discuss_mort}
In Section \ref{sec:mortar}, we introduced two mortar functions: (1) Gaussian, and (2) Algebraic. The former was proposed in \cite{khan2024high} and the latter herein. In \cite{khan2024high}, a third type of mortar function, called Fickian, was proposed that involves the initialization of a box function over each mortar node, then smoothing it by solving a $D-1$-dimensional diffusion equation on the corresponding contact interface. We found that none of these choices has a sizable impact on the accuracy of first-pass solutions produced by hPLMM, and by extension, the performance of GMRES. The caveat is that a good user-defined parameter must be set for the Gaussian ($\beta$ in Eq.\ref{eq:gauss}) and Fickian (the extent box functions are smoothed) mortars. By contrast, Algebraic mortars do not have any adjustable parameters, which is an advantage. On the other hand, both Algebraic and Fickian mortars require the solution of a $D-1$-dimensional PDE on contact interfaces, which is more costly than the analytically defined Gaussian mortars in Eq.\ref{eq:gauss}. In \ref{app:gauss_supp}, we show that the $\beta$ parameter in Gaussian mortars, which controls the spread of the functions, cannot be too small or too large. If too small, stresses become oscillatory along contact interfaces and the approximation errors grow (see Fig.D.4 in \cite{khan2024high}). But if too large, the overlap between adjacent mortars increases and the condition number of the coarse matrix $\mathrm{\hat{R}\hat{A}\hat{P}}$ in Eq.\ref{eq:mg_struct} deteriorates. Bad conditioning implies poor approximation. For all the domains herein, $\beta\!=\!4$ was found to be a balanced choice. In \ref{app:alg_vs_gauss}, we further compare Gaussian mortars to Algebraic mortars and find they both produce first-pass solutions with similar accuracy. However, as the number of mortar nodes per interface, $n$, grows, Gaussian mortars converge slightly faster to the exact solution, which is another reason we chose them in Section \ref{sec:results}. 

Lastly, we remark on the placement of mortar nodes on contact interfaces: The algorithm outlined in Section \ref{sec:mortar}, proposed by \cite{khan2024high}, works well for all domains, except S2DH. From Fig.\ref{fig:S2DH_fp}, we see the first-pass solution of S2DH ceases to improve beyond some $n$. We attribute this to the fact that the placement of mortar nodes is not informed by the underlying heterogeneity in stiffness tensor, which exhibits jump discontinuities near interfaces. Accounting for such heterogeneity may require better mortar-node placement, or mortar functions derived from eigenvalue problems \cite{heinlein2019adaptive}.

\subsection{Computational complexity}\label{sec:comp complex}
The computational complexity of the geometric hPLMM was outlined in \cite{khan2024high}. Here, we detail the complexity of building and applying the algebraic hPLMM preconditioner. Our analysis parallels that of \cite{mehmani2023precond, li2024phase} for PLMM, and generalizes it to hPLMM. Let $\Omega_s$ consist of $N^f$ fine grids (FEM nodes), $N^g$ grain grids, $N^\zeta$ contact grids, and $N^m$ mortars. Assuming all fine-scale computations are performed by a linear solver whose WCT scales like $O(N^\vartheta)$, where $N$ is the number of unknowns and $\vartheta \in(1,3)$, then the WCTs of building M in Eq.\ref{eq:multi_precond}, and applying it once are:
\begin{subequations} \label{eq:complexity}
\begin{align}
	&\mathcal{T}^{M}_{build} = 
	\mathcal{T}^{M_G}_{build} + \mathcal{T}^{M_L}_{build}
	\label{eq:Mbuild_wct} \\
	&\mathcal{T}^{M}_{apply} = 
	\mathcal{T}^{M_G}_{apply} + \mathcal{T}^{M_L}_{apply}
	\label{eq:Mapply_wct}
\end{align}
where
\begin{align}
	& \mathcal{T}^{M_G}_{build} = 
	   O(N^f D / N^g)^{\vartheta} \times (2N^m D + N^g) / P_{prc}
	    \label{eq:MGbuild_wct}\\
	&\mathcal{T}^{M_G}_{apply} = 
	   O(N^m D + N^g)^{\vartheta} \label{eq:MGapply_wct}\\
    & \mathcal{T}^{M_L}_{apply} = O\left[
	   (N^f D / N^g)^{\vartheta} \times N^g / P_{prc} + 
	   (f^{\zeta}N^f D / N^\zeta)^{\vartheta} \times N^\zeta / P_{prc}\right]
	   \label{eq:MLapply_wct}
\end{align}
\end{subequations}
In Eq.\ref{eq:complexity}, $P_{prc}$ is the number of parallel processors used, and $f^\zeta$ is the fraction of $\Omega_s$ covered by the union of all contact grids (see Fig.\ref{fig:schem}). Thus, the cost of solving one local system on a grain grid is proportional to $(N^fD/N^g)^\vartheta$, whereas the cost of solving one local system on a contact grid is proportional to $(f^\zeta N^fD/ N^\zeta)^\vartheta$. Recall the number of contact grids, $N^\zeta$, equals the number of interfaces, $N^c$. To build M\ts{G}, $D\!\times\! 2$ grain-grid problems must be solved per mortar node, plus $N^g$ additional grain-grid problems corresponding to the correction vectors in Eq.\ref{eq:bc_matrix}. To apply M\ts{G}, only a coarse system ($\mathrm{\hat{R}\hat{A}\hat{P}}$) must be solved, which has $N^mD + N^g$ rows and columns. As for the smoother, we have assumed the use of M\ts{CG} with the number of stages $n_{st}\!=\!1$, which was identified as the best performing choice in Section \ref{sec:results}. To apply this additive-Schwarz M\ts{L} once, $N^g$ grain-grid problems and $N^\zeta$ contact-grid problems must be solved.

Notice the build-time of M\ts{L}, $\smash{\mathcal{T}^{M_L}_{build}}$, is left unspecified in Eq.\ref{eq:complexity} because, under normal circumstances, this cost is zero. But if the local systems associated with all grain grids and contact grids are LU-decomposed during a preprocessing step, then a significant part of $\smash{\mathcal{T}^{M_L}_{apply}}$ would be front-loaded as a one-time expense, reducing the cost of GMRES iterations (as was done in Section \ref{sec:results}). In that case, the cost of doing LU-decompositions constitutes $\smash{\mathcal{T}^{M_L}_{build}}$.

All computations, except those associated with applying M\ts{G}, are parallelizable since grain-grid and contact-grid problems are fully decoupled. The maximum number of parallel processors, $P_{prc}$, one could use when building M\ts{G} is $2N^mD + N^g$, whereas the maximum one can use when applying (or building) M\ts{L} is $\max(N^g + N^\zeta)$. Given a careful parallel scalability study is outside the scope of this paper, all computations in Section \ref{sec:results} were performed in series.

From Eq.\ref{eq:complexity}, we see that the cost of applying M\ts{G} scales with $\smash{(N^m)^\vartheta}$, where $N^m\!=\!nN^c$ and $n$ is again the number of mortar nodes per interface. In PLMM \cite{li2024phase}, $n\!=\!1$ and M\ts{G} is very cheap to apply, because $N^c\!\ll\! N^f$. However, in hPLMM, as $n$ grows to $>\!4$, the cost of applying M\ts{G} starts to dominate the self-time of GMRES (T\ts{sol} in Table \ref{tab:wct}).

\subsection{Applications beyond linear-elasticity}\label{sec:nonlinear_app}
As noted in Section \ref{sec:multilevel}, PLMM ($n\!=\!1$) and hPLMM ($n\!=\!2-4$) have similar WCTs if the goal is to solve a single system. But because hPLMM front-loads much of its cost as a one-time expense towards building M\ts{G}, which itself is parallelizable, it is faster than PLMM in solving multiple systems by at least a factor of two. Few examples where this is useful includes plasticity, finite-strain mechanics, wave propagation, and phase-field simulations of fracture. Typical algorithms of these problems consist of multiple, often nested, loops for taking load steps, time steps, and performing staggered or Newton iterations \cite{neto2011PlasticFiniteStrain,hughes2012finite,ambati2015review}. In each step or iteration, a linear(ized) system of the form Eq.\ref{eq:ls_all} must be solved, where the RHS vector is altered and the matrix $\mathrm{\hat{A}}$ is potentially perturbed. If the perturbations are localized (i.e., high-frequency) or global but not too severe, then the M\ts{G} built at the start of a simulation can be reused later, as was demonstrated by \cite{li2024phase} for PLMM in the case of phase-field simulations. Otherwise, M\ts{G} must be updated periodically. However, from experience \cite{li2024phase, mehmani2023precond}, the need for such updates is often neither frequent nor demanding of rebuilding M\ts{G} from scratch. Adaptive update criteria can be devised (see \cite{li2024phase}) that identify grain grids with the largest deviation/error and replace the corresponding column in the reduced prolongation matrix, P, in Eq.\ref{eq:bc_matrix}.

\section{Conclusion} \label{sec:conclusion}
In this work, we have successfully algebraized the high-order pore-level multiscale method (hPLMM) of \cite{khan2024high} into a two-level preconditioner for solving linear-elastic deformation in domains with arbitrary geometry and material heterogeneity. It consists of a coarse preconditioner, M\ts{G}, and a fine-scale smoother, M\ts{L}, where the most compatible choice for the latter is M\ts{CG} in Section \ref{sec:local_smoother}, instead of black-box smoothers like M\ts{ILU($k$)}. When applied within GMRES, the hPLMM preconditioner converges faster, both in terms of iterations and wall-clock times (WCT), than other state-of-the-art preconditioners such as the low-order PLMM \cite{li2024phase}, cAMG \cite{gustafsson1998cAMGorig}, and GDSW \cite{heinlein2020fully}. This is achieved through the introduction of mortars in hPLMM, or extra degrees of freedom, at interfaces between subdomains that enable the accurate capturing of local bending/torsion moments, unlike PLMM. As the number of mortar nodes per interface, $n$, grows, the cost of building M\ts{G}, denoted by T\ts{MG}, increases but the self-time of GMRES, T\ts{sol}, is reduced. When $n\!<\!4$, the sum of these two costs is roughly constant, implying that for the same effort, linear systems can be solved at least twice as fast (smaller T\ts{sol}). Of course, if the goal is to solve only one system, there are no benefits to $n\!>\!1$. But if the goal is to solve multiple systems, with different RHS vectors or slightly perturbed coefficient matrices \cite{li2024phase}, then substantial ($>\!2$ times) cost savings are possible with $n\!>\!1$ compared to $n\!=\!1$, where the latter reduces hPLMM to the low-order PLMM. Scenarios where repeated solutions of similar systems are needed abound in the literature, from performing multiple load/time steps to model quasi-static/dynamic deformation of a solid to the staggered or Newton iterations needed to model its failure or plastic yield. For $n\!>\!4$, the total cost (T\ts{MG} + T\ts{sol}) of hPLMM grows with $n$, as the size of the coarse system grows too. Multilevel ($>\!2$) extensions of hPLMM could alleviate this. Finally, the building and application of hPLMM are both amenable to parallelism, which we did not exploit in this work.

\appendix
\setcounter{figure}{0}
\section{Impact of $\beta$ in Gaussian mortars}
\label{app:gauss_supp}
\begin{figure} [t!]
  \centering
  \centerline{\includegraphics[scale=0.25,trim={190 20 215 20},clip]{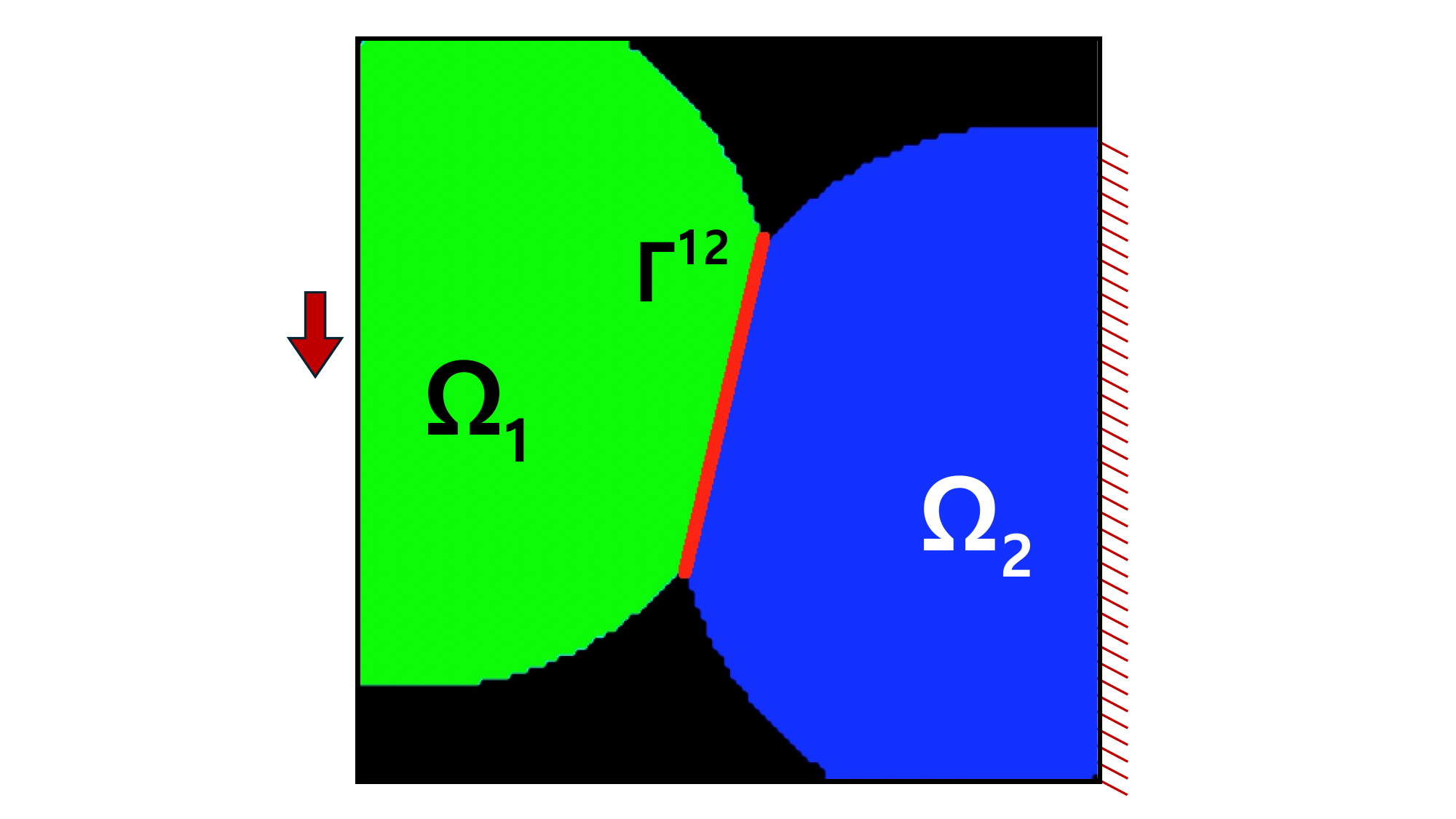}}
	\caption{A two-grain domain comprised of two grain grids and one contact interface, used to quantify the impact of $\beta$ in Eq.\ref{eq:gauss}.}
\label{fig:two_grain}
\end{figure}

\begin{figure} [h!]
  \centering
  \centerline{\includegraphics[scale=0.7,trim={155 120 190 97},clip]{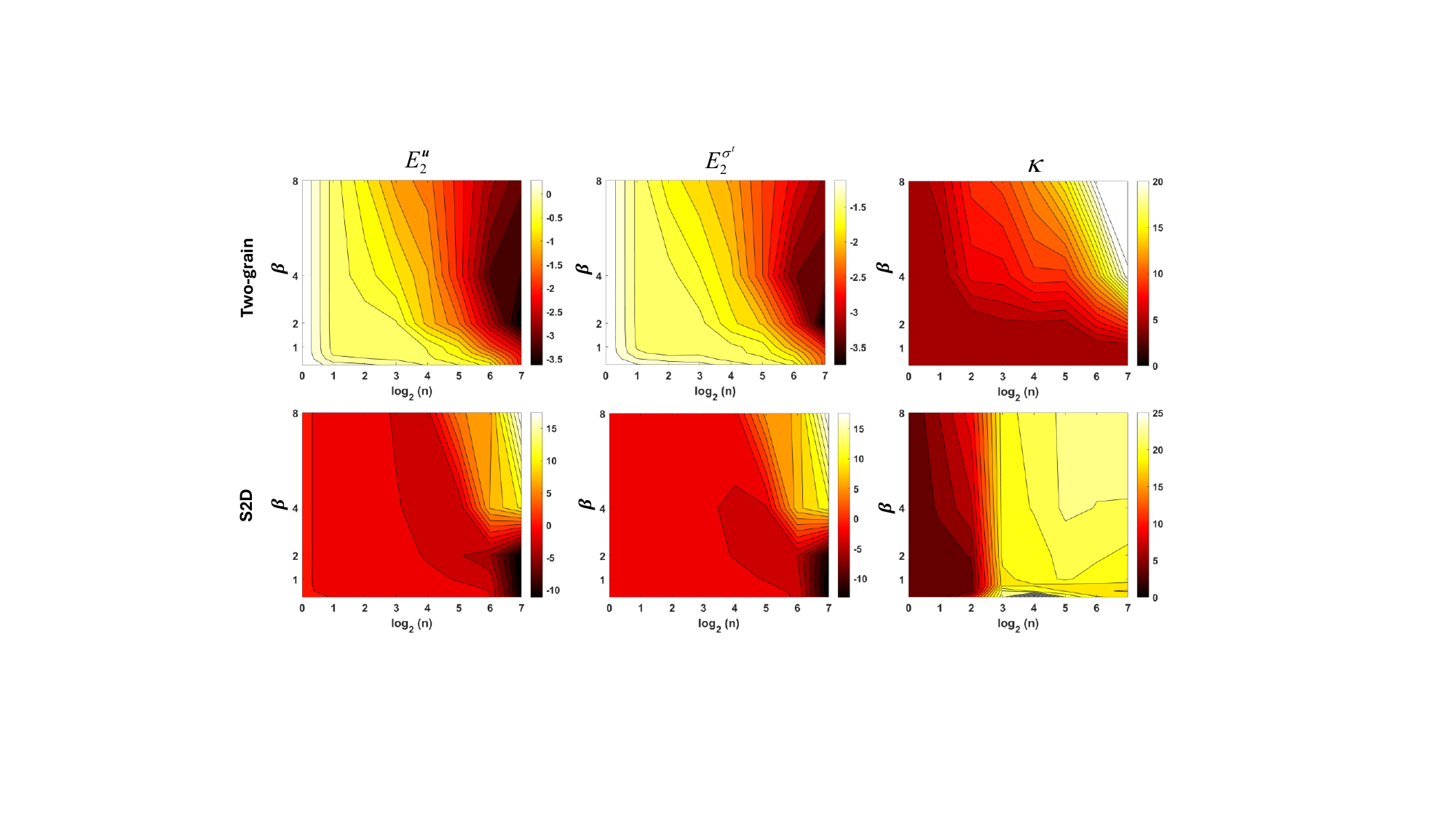}}
	\caption{Contour plots of L$_2$-errors (log-scale) in displacement, $E^{\bs{u}}_2$, and maximum shear stress, $\smash{E^{\sigma^t}_2}$, for the first-pass solutions of hPLMM for different $\beta$ in Eq.\ref{eq:gauss} and number of mortar nodes per interface, $n$. The top row corresponds to the two-grain domain in Fig.\ref{fig:two_grain} and the bottom row to the S2D domain in Fig.\ref{fig:problem_set1}, both subjected to shear loading. The third column depicts the condition number, $\kappa$, of the coarse system $\mathrm{\hat{R}\hat{A}\hat{P}}$ in Eq.\ref{eq:mg_struct}.}
\label{fig:gauss_supp}
\end{figure}

The parameter $\beta$ in Eq.\ref{eq:gauss} controls the spread of the Gaussian mortar function. A large $\beta$ implies more spread, thus overlap with adjacent mortars, whereas a small $\beta$ results in a more localized mortar function around its node. Fig.\ref{fig:gauss_supp} depicts contour plots of L$_2$-error for the first-pass solution of hPLMM for different values of $\beta$ and number of mortar nodes per interface, $n$. Errors are computed via Eq.\ref{eq:err_metric} for displacement, $E^{\bs{u}}_2$, and maximum shear stress, $\smash{E^{\sigma^t}_2}$, on the S2D domain in Fig.\ref{fig:problem_set1} and the simpler two-grain domain in Fig.\ref{fig:two_grain}. Both domains are subjected to shear loading, as detailed in Section \ref{sec:prob_set} and annotated in Fig.\ref{fig:two_grain}. When decomposed, the two-grain domain consists of only two grain grids and one contact interface, providing a simpler setting for quantifying the influence of $\beta$. The general observation from Fig.\ref{fig:gauss_supp} is that at small to moderate $n$ ($<\!2^6$), errors decrease with growing $\beta$. However, at high $n$, errors decrease up to some $\beta$ then increase for larger $\beta$. To understand this behavior, we have included plots of the condition number of the coarse matrix $\mathrm{\hat{R}\hat{A}\hat{P}}$ in Eq.\ref{eq:mg_struct}, $\kappa$, versus $\beta$ and $n$. We see that as $\beta$ is increased, so is $\kappa$, and the increase in $\kappa$ is fastest at large $n$. The growth of $\kappa$ with $\beta$ is due to an increase in the overlap between adjacent mortars and the resulting loss of linear independence (in finite-precision arithmetic). Thus, at large $n$ and $\beta$, the accuracy of the the first-pass solution deteriorates due to numerical inaccuracies in applying the inverse of $\mathrm{\hat{R}\hat{A}\hat{P}}$ in Eq.\ref{eq:mg_struct}. While an optimal value for $\beta$ is geometry dependent, Fig.\ref{fig:gauss_supp} suggests $\beta\!=\!4$ is an acceptable choice, which we use throughout this work.

\section{Comparison between Gaussian and Algebraic mortars}
\label{app:alg_vs_gauss}
\begin{figure} [h!]
  \centering
  \centerline{\includegraphics[scale=0.5,trim={150 23 165 30},clip]{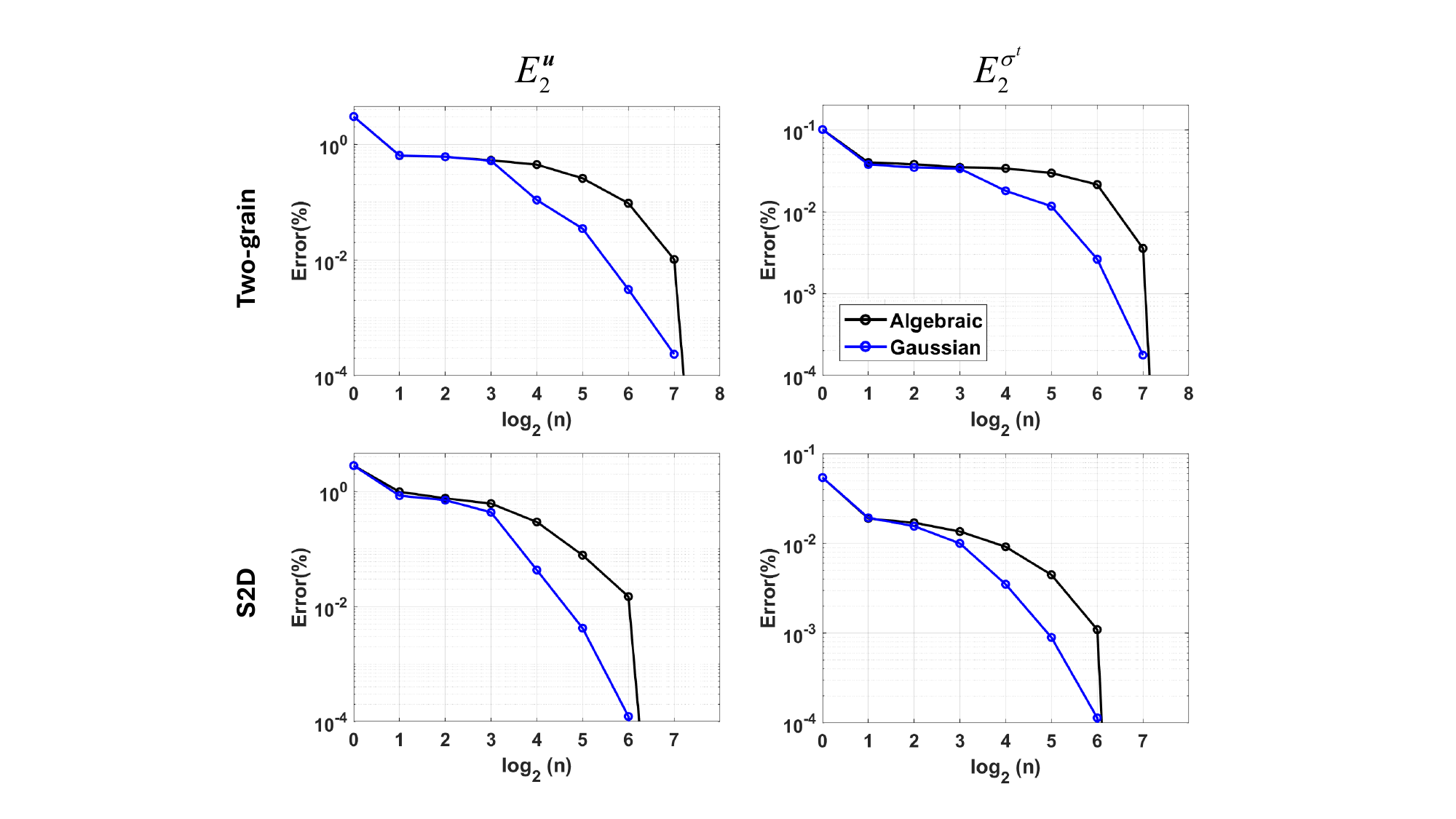}}
	\caption{Comparison between Gaussian and Algebraic mortars in terms of the L$_2$-errors in displacement, $E^{\bs{u}}_2$, and maximum shear stress, $\smash{E^{\sigma^t}_2}$, obtained when they are used in computing the first-pass solutions of hPLMM with different number of mortar nodes per interface, $n$. The plotted errors correspond to (top row) the two-grain domain in Fig.\ref{fig:two_grain} and (bottom row) the S2D domain in Fig.\ref{fig:problem_set1}, both subjected to shear loading.}
\label{fig:alg_vs_gauss}
\end{figure}

In Section \ref{sec:mortar}, we introduced two kinds of mortar functions: Gaussian and Algebraic. While Algebraic mortars, unlike Gaussian mortars, are not burdened by any user-defined parameters like $\beta$, they do require the solution of a $D-1$-dimensional PDE on each contact interface. Fig.\ref{fig:alg_vs_gauss} compares L$_2$-errors of displacement, $E^{\bs{u}}_2$, and maximum shear stress, $\smash{E^{\sigma^t}_2}$, obtained from using each of these mortar types to compute the first-pass solution of hPLMM with different number of mortar nodes per interface, $n$. Results are shown for the two-grain domain in Fig.\ref{fig:two_grain} and the S2D domain in Fig.\ref{fig:problem_set1}, both subjected to shear loading. Fig.\ref{fig:alg_vs_gauss} shows the two mortar functions perform comparably, with Gaussian mortars converging slightly faster when $n\!>\!8$. The steepest drop in error occurs from $n\!=\!1$ to $2$, which is consistent with improvements seen in the first-pass solutions of Fig.\ref{fig:P2D_S2D_fp}. At $n\!>\!2^7$, the number of mortar nodes exceeds the number of fine grids on all contact interfaces, making the coarse matrix $\mathrm{\hat{R}\hat{A}\hat{P}}$ in Eq.\ref{eq:mg_struct} rank deficient. Hence, we set $n$ equal to the number of fine grids on each interface, yielding machine precision errors (out of bounds in Fig.\ref{fig:alg_vs_gauss}). Given Gaussian mortars converge somewhat faster, and are cheaper to construct, we used them throughout Section \ref{sec:results}.

\section{Decomposition for the GDSW preconditioner}
\label{app:GDSW}
\begin{figure} [h!]
  \centering
  \centerline{\includegraphics[scale=0.45,trim={0 130 0 140},clip]{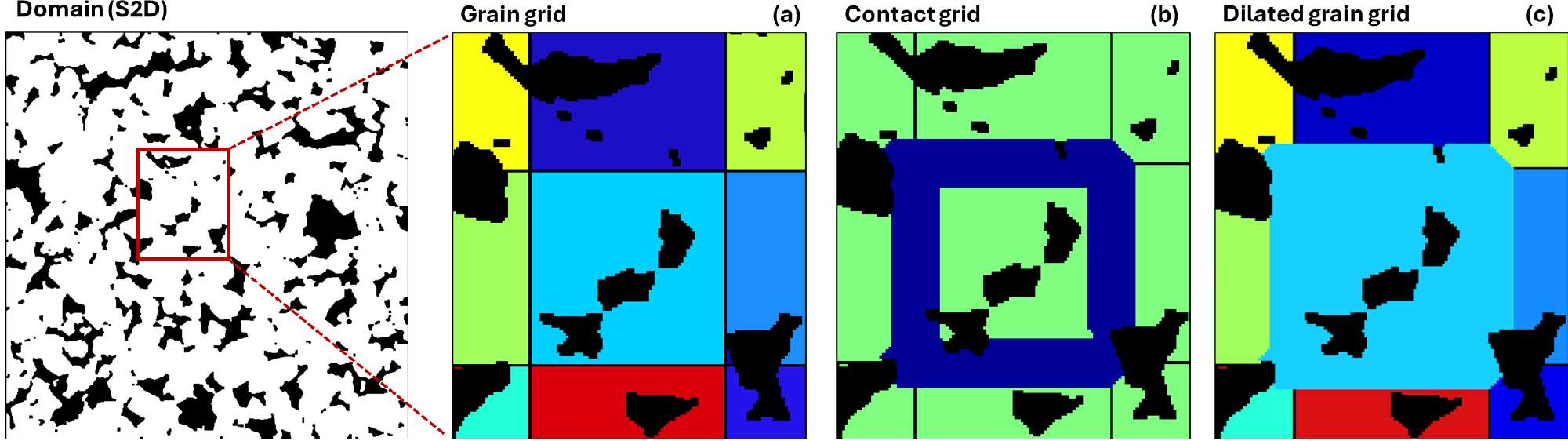}}
	\caption{Decomposition of the S2D domain into Cartesian grids for the GDSW preconditioner. Plot (a) shows a zoom-in of a grain grid in light blue. Plot (b) depicts the corresponding contact grid (16 pixels wide). Plot (c) shows the grain grid in (a) dilated to create overlap with neighboring grain grids. In GDSW, the grain grids in (a) are used to build M\ts{G} (same way as PLMM), and the dilated grain grids in (c) are used to build M\ts{L}.}
\label{fig:GDSW}
\end{figure}

Fig.\ref{fig:GDSW} illustrates the Cartesian decomposition used to construct the GDSW preconditioner for the S2D domain in Fig.\ref{fig:problem_set1}. As noted in Section \ref{sec:prob_set}, the coarse preconditioner M\ts{G} in GDSW can be built using the same algorithm as the low-order PLMM from \cite{li2024phase} (or hPLMM with $n\!=\!1$). Since GDSW does not prescribe a specific decomposition (unlike PLMM), we employ a Cartesian partitioning of $\Omega_s$ into ``grain grids'' (Fig.\ref{fig:GDSW}a). The GDSW smoother M\ts{L} is an additive-Schwarz preconditioner that uses overlapping subdomains derived from the same grain grids, but dilated to create overlap with neighboring regions (Fig.\ref{fig:GDSW}c). The overlap width between adjacent dilated grain grids is 16 pixels (fine-grid FEM elements), comparable to the contact-grid widths used in the hPLMM smoother (Fig.\ref{fig:GDSW}b).

\section{Krylov convergence under tensile loading}
\label{app:Krylov_tensile}
Table \ref{tab:wct_tense} is the counterpart of Table \ref{tab:wct} but when all domains in Figs.\ref{fig:problem_set1}-\ref{fig:problem_set2} are subjected to tensile loading instead of shear. Included in Table \ref{tab:wct_tense} are the number of GMRES iterations and total wall-clock times (WCT), T\ts{tot}, the latter composed of the cost of building the smoother, T\ts{ML}, cost of building the coarse preconditioner, T\ts{MG}, and the self-time of GMRES, T\ts{sol}. The values are almost identical to those in Table \ref{tab:wct} for shear, except convergence is slightly faster under tension (esp. for S2DH and DARCY). The reason is the same as in Section \ref{sec:first_pass} for the higher accuracy of hPLMM's first-pass solution under tension than shear. Namely, tensile loading produces smaller local bending/torsion moments than shear, requiring fewer mortar nodes to capture the displacement field along contact interfaces.

\begin{table}[t!]
\centering
\caption{Summary of the number of iterations and wall-clock times (WCTs) in seconds required by GMRES to converge ($\| \hat{\mathrm{A}} \hat{x} - \hat{b} \| / \| \hat{b} \| < 10^{-9}$) to the solution of Eq.\ref{eq:ls_all} using different preconditioners. This includes hPLMM, where M\ts{G} is combined with either the M\ts{CG} ($n_{st}\!=\!1$) or M\ts{ILU($k$)} ($n_{st}\!=\!6$) smoother, and the number of mortar nodes per interface is set to $n\!=\!1$, 2, 4, and 8. The $n\!=\!1$ case corresponds to the low-order PLMM. Results for GDSW and cAMG are included as benchmarks. The total WCT, T\ts{tot}, consists of the time needed to build the smoother, T\ts{ML}, time to build the coarse preconditioner, T\ts{MG}, and the time spent by GMRES itself, T\ts{sol}. In all domains, the ILU fill-level is $k\!=\!0$ except in S2DH, where it is $k\!=\!1$. All cases correspond to \textit{tensile loading} and \textcolor{red}{red} means ``diverged'' (i.e., not converged within 300 iterations).}
\footnotesize
\setlength{\tabcolsep}{3.5pt}
\begin{tabular}{l l *{4}{r} @{\hspace{12pt}} *{4}{r} @{\hspace{12pt}} r r}
\toprule
  & & \multicolumn{4}{c}{\textbf{Contact-Grain (CG)}} 
  & \multicolumn{4}{c}{\textbf{ILU(0) / ILU(1)}} 
  &  &  \\
\cmidrule(lr){3-6} \cmidrule(lr){7-10} 
 \textbf{Domain} & \textbf{Metric} & $n=1$ & $n=2$ & $n=4$ & $n=8$ & $n=1$ & $n=2$ & $n=4$ & $n=8$ & \textbf{GDSW} & \textbf{cAMG} \\
\midrule

\multirow{5}{*}{P2D} 
 & T\textsubscript{ML}  & 7.49 & 7.49 & 8.77 & 8.71 & 0.11 & 0.11 & 0.11 & 0.11 & 10.57 & -- \\
 & T\textsubscript{MG}  & 10.67 & 14.14 & 22.96 & 40.80 & 10.60 & 13.99 & 22.78 & 40.80 & 9.89 & 1.81 \\
 & T\textsubscript{sol} & 5.39 & 3.83 & 3.05 & 4.21 & 43.70 & 20.35 & 20.21 & 24.35 & 22.18 & 168.4 \\
 & \textbf{T\textsubscript{tot}} & \textbf{23.55} & \textbf{25.46} & \textbf{34.78} & \textbf{53.72} & \textbf{54.41} & \textbf{34.45} & \textbf{43.10} & \textbf{65.26} & \textbf{42.64} & \textbf{170.2} \\
 & \textbf{Iter.}  & \textbf{16} & \textbf{11} & \textbf{8} & \textbf{7} & \textbf{87} & \textbf{47} & \textbf{43} & \textbf{37} & \textbf{45} & \textbf{64} \\
\midrule

\multirow{5}{*}{S2D} 
 & T\textsubscript{ML}  & 4.26 & 4.33 & 5.13 & 5.17 & 0.08 & 0.08 & 0.08 & 0.08 & 6.37 & -- \\
 & T\textsubscript{MG}  & 5.56 & 7.24 & 12.54 & 23.24 & 5.72 & 7.31 & 12.21 & 23.45 & 5.20 & 1.19 \\
 & T\textsubscript{sol} & 6.12 & 2.76 & 2.81 & 4.62 & 30.07 & 19.04 & 13.10 & 21.19 & 16.56 & 207.7 \\
 & \textbf{T\textsubscript{tot}} & \textbf{15.94} & \textbf{14.33} & \textbf{20.48} & \textbf{33.03} & \textbf{35.87} & \textbf{26.43} & \textbf{25.39} & \textbf{44.72} & \textbf{28.13} & \textbf{225.3} \\
 & \textbf{Iter.}  & \textbf{28} & \textbf{13} & \textbf{11} & \textbf{9} & \textbf{87} & \textbf{60} & \textbf{42} & \textbf{36} & \textbf{55} & \textbf{115} \\
\midrule

\multirow{5}{*}{S2DH} 
 & T\textsubscript{ML}  & 4.37 & 4.26 & 5.18 & 4.14 & \textcolor{red}{0.70} & \textcolor{red}{0.73} & \textcolor{red}{0.72} & \textcolor{red}{0.71} & 8.08 & -- \\
 & T\textsubscript{MG}  & 5.82 & 7.58 & 12.66 & 21.81 & \textcolor{red}{5.59} & \textcolor{red}{7.43} & \textcolor{red}{12.52} & \textcolor{red}{22.24} & 6.11 & 1.19 \\
 & T\textsubscript{sol} & 12.33 & 11.71 & 12.23 & 15.55 & \textcolor{red}{346.3} & \textcolor{red}{346.7} & \textcolor{red}{355.7} & \textcolor{red}{435.5} & 125.0 & 365.9 \\
 & \textbf{T\textsubscript{tot}} & \textbf{22.52} & \textbf{23.55} & \textbf{30.07} & \textbf{41.50} & \textcolor{red}{\textbf{352.6}} & \textcolor{red}{\textbf{354.9}} & \textcolor{red}{\textbf{369.0}} & \textcolor{red}{\textbf{458.4}} & \textbf{139.2} & \textbf{375.2} \\
 & \textbf{Iter.}  & \textbf{51} & \textbf{45} & \textbf{38} & \textbf{29} & \textcolor{red}{\textbf{300}} & \textcolor{red}{\textbf{300}} & \textcolor{red}{\textbf{300}} & \textcolor{red}{\textbf{300}} & \textbf{222} & \textbf{179} \\
\midrule

\multirow{5}{*}{S2DC} 
 & T\textsubscript{ML}  & 3.55 & 3.55 & 4.34 & 4.38 & 0.12 & 0.11 & 0.11 & 0.12 & 5.22 & -- \\
 & T\textsubscript{MG}  & 4.93 & 6.71 & 11.53 & 22.00 & 4.99 & 6.66 & 12.42 & 22.20 & 4.46 & 1.12 \\
 & T\textsubscript{sol} & 10.04 & 6.57 & 5.97 & 7.99 & 46.70 & 23.51 & 22.22 & 23.19 & 23.20 & 192.0 \\
 & \textbf{T\textsubscript{tot}} & \textbf{18.52} & \textbf{16.83} & \textbf{21.84} & \textbf{34.37} & \textbf{51.81} & \textbf{30.28} & \textbf{34.75} & \textbf{45.51} & \textbf{32.88} & \textbf{198.3} \\
 & \textbf{Iter.}  & \textbf{42} & \textbf{28} & \textbf{20} & \textbf{15} & \textbf{116} & \textbf{67} & \textbf{52} & \textbf{37} & \textbf{74} & \textbf{105} \\
\midrule

\multirow{5}{*}{DARCY Cart} 
 & T\textsubscript{ML}  & 6.96 & 6.95 & 6.79 & 6.82 & 0.22 & 0.19 & 0.19 & 0.22 & 9.36 & -- \\
 & T\textsubscript{MG}  & 9.71 & 12.62 & 19.35 & 35.60 & 9.60 & 12.57 & 19.33 & 36.28 & 9.93 & 1.98 \\
 & T\textsubscript{sol} & 14.97 & 6.72 & 4.09 & 3.96 & 36.62 & 18.40 & 14.78 & 16.99 & 21.95 & 76.60 \\
 & \textbf{T\textsubscript{tot}} & \textbf{31.64} & \textbf{26.29} & \textbf{30.23} & \textbf{46.38} & \textbf{46.44} & \textbf{31.16} & \textbf{34.30} & \textbf{53.49} & \textbf{41.24} & \textbf{78.58} \\
 & \textbf{Iter.}  & \textbf{36} & \textbf{17} & \textbf{10} & \textbf{8} & \textbf{63} & \textbf{35} & \textbf{28} & \textbf{27} & \textbf{44} & \textbf{27} \\
\midrule

\multirow{5}{*}{DARCY Spec} 
 & T\textsubscript{ML}  & 7.50 & 7.70 & 7.40 & 7.27 & 0.19 & 0.20 & 0.20 & 0.20 & 10.13 & -- \\
 & T\textsubscript{MG}  & 11.97 & 16.49 & 25.83 & 49.90 & 11.72 & 16.46 & 26.18 & 49.77 & 10.33 & 2.41 \\
 & T\textsubscript{sol} & 15.28 & 8.42 & 5.93 & 5.51 & 28.88 & 17.91 & 16.47 & 18.96 & 22.92 & 76.73 \\
 & \textbf{T\textsubscript{tot}} & \textbf{34.75} & \textbf{32.61} & \textbf{39.16} & \textbf{62.68} & \textbf{40.79} & \textbf{34.57} & \textbf{42.85} & \textbf{68.93} & \textbf{43.38} & \textbf{79.14} \\
 & \textbf{Iter.}  & \textbf{36} & \textbf{19} & \textbf{13} & \textbf{10} & \textbf{54} & \textbf{34} & \textbf{29} & \textbf{28} & \textbf{44} & \textbf{27} \\
\midrule

\multirow{5}{*}{P2DL} 
 & T\textsubscript{ML}  & 205.1 & 152.5 & 205.5 & 205.0 & 2.47 & 2.26 & 2.64 & 2.33 & 118.2 & -- \\
 & T\textsubscript{MG}  & 194.1 & 262.3 & 443.5 & 919.8 & 194.0 & 266.1 & 446.0 & 917.9 & 162.1 & 36.47 \\
 & T\textsubscript{sol} & 120.1 & 85.55 & 82.43 & 119.0 & 3998.9 & 436.4 & 510.1 & 739.8 & 613.5 & 2667.2 \\
 & \textbf{T\textsubscript{tot}} & \textbf{519.3} & \textbf{500.3} & \textbf{731.4} & \textbf{1243.8} & \textbf{4195.4} & \textbf{704.7} & \textbf{958.7} & \textbf{1660.0} & \textbf{893.9} & \textbf{3119.3} \\
 & \textbf{Iter.}  & \textbf{17} & \textbf{11} & \textbf{9} & \textbf{7} & \textbf{260} & \textbf{44} & \textbf{43} & \textbf{38} & \textbf{78} & \textbf{59} \\
\midrule

\multirow{5}{*}{\textbf{BONE}} 
 & T\textsubscript{ML}  & 150.6 & 148.1 & 146.4 & 146.4 & 5.66 & 5.97 & 6.07 & 5.80 & 168.2 & -- \\
 & T\textsubscript{MG}  & 198.5 & 357.5 & 761.2 & 1924.2 & 195.5 & 354.5 & 749.5 & 2055.5 & 102.0 & 8.79 \\
 & T\textsubscript{sol} & 195.2 & 152.6 & 420.6 & 1904.3 & 1338.0 & 1448.5 & 4360.3 & 23229.0 & 839.3 & 2713.8 \\
 & \textbf{T\textsubscript{tot}} & \textbf{544.3} & \textbf{658.2} & \textbf{1328.1} & \textbf{3974.9} & \textbf{1539.2} & \textbf{1809.0} & \textbf{5115.9} & \textbf{25290.3} & \textbf{1109.4} & \textbf{2722.6} \\
 & \textbf{Iter.} & \textbf{42} & \textbf{28} & \textbf{22} & \textbf{17} & \textbf{298} & \textbf{280} & \textbf{259} & \textbf{206} & \textbf{67} & \textbf{222} \\

\bottomrule
\end{tabular}
\label{tab:wct_tense}
\end{table} 

\section*{Acknowledgments}
This material is based upon work supported by the National Science Foundation under Grant No. CMMI-2145222. We acknowledge the Institute for Computational and Data Sciences (ICDS) at Penn State University for access to computational resources.

\section*{References}
\bibliographystyle{unsrtnat}
\bibliography{./References.bib}

\end{document}